\documentclass[11pt]{amsart}
\usepackage[utf8]{inputenc}
\usepackage[T1]{fontenc} 
\usepackage{lmodern}
\usepackage{amsmath}
\usepackage{amssymb}
\usepackage{amsthm}
\usepackage{amsfonts}
\usepackage{stmaryrd}
\usepackage{bbm}
\usepackage{enumitem}
\usepackage[dvipsnames]{xcolor}
\usepackage{tikz}
\usetikzlibrary{decorations.markings}
\usepackage{quiver}
\usepackage{hyperref}
\usepackage{adjustbox}
\usepackage{graphicx}
\usepackage{textcomp}
\usepackage[normalem]{ulem}
\usepackage[all]{xy}
\usepackage{cancel}
\usepackage{scalerel}
\usepackage[backend=biber, 
style=alphabetic,
minalphanames=4,
maxalphanames=5,
maxbibnames=5]{biblatex}

\DeclareFontFamily{U}{mathx}{}
\DeclareFontShape{U}{mathx}{m}{n}{<-> mathx10}{}
\DeclareSymbolFont{mathx}{U}{mathx}{m}{n}
\DeclareMathAccent{\widehat}{0}{mathx}{"70}
\DeclareMathAccent{\widecheck}{0}{mathx}{"71}
\DeclareMathAccent{\widebar}{0}{mathx}{"73}

\theoremstyle{plain}
 
\makeatletter
\@namedef{subjclassname@1991}{2020 Mathematics Subject Classification}
\makeatother

\usepackage{chngcntr}

\counterwithin*{equation}{section}

\newtheorem{theorem}{Theorem}[section]
\newtheorem{lemma}[theorem]{Lemma}
\newtheorem{proposition}[theorem]{Proposition}
\newtheorem{corollary}[theorem]{Corollary}
\theoremstyle{definition}
\newtheorem{definition}[theorem]{Definition} 
\newtheorem{example}[theorem]{Example}
\newtheorem{remark}[theorem]{Remark}

\newcommand{\Cat}{\mathcal C}% category C
\newcommand{\Set}{\mathsf{Set}}% category of sets
\newcommand{\Mfd}{\mathsf{Mfd}}% category of smooth manifolds
\newcommand{\Vect}{\mathsf{Vec}}
\newcommand{\SVect}{\mathsf{SVec}}
\newcommand{\SSet}{\mathsf{SSet}}
\newcommand{\Simp}{\mathsf{S}}

\newcommand{\VS}{\mathsf{VS}}
\newcommand{\VB}{\mathsf{VB}}

\newcommand{\Ch}{\mathsf{Ch}}
\newcommand{\hSVect}{\mathsf{hSVec}}
\newcommand{\hCh}{\mathsf{hCh}}
\newcommand{\hChP}{\mathsf{hCh_{\ge 0}}}

\newcommand{\Shuf}{\mathsf{Shuf}} %shuffle
\newcommand{\Sh}{\mathsf{Sh}} %Sheaves cat
\DeclareMathOperator{\IHom}{\underline{Hom}}%internal hom
\newcommand{\stkout}[1]{\ifmmode\text{\sout{\ensuremath{#1}}}\else\sout{#1}\fi} %%%%%%strikeout in math mode, requires amsmath, ulem

\renewcommand{\epsilon}{\varepsilon}

\newcommand{\R}{\mathbb{R}}
\newcommand{\Z}{\mathbb{Z}}

\newcommand{\into}{\hookrightarrow}

\newcommand\rightthreearrows{% Three arrows on top of each other
        \mathrel{\vcenter{\mathsurround0pt
                \ialign{##\crcr
                        \noalign{\nointerlineskip}$\rightarrow$\crcr
                        \noalign{\nointerlineskip}$\rightarrow$\crcr
                        \noalign{\nointerlineskip}$\rightarrow$\crcr
                }%
        }}%
}

\makeatletter
\newcommand{\newrightleftarrows}[2]{%
  \mathrel{\mathop{%
    \vcenter{\offinterlineskip\m@th
      \ialign{\hfil##\hfil\cr
        \hphantom{$\scriptstyle\mspace{8mu}{#1}\mspace{8mu}$}\cr
        \rightarrowfill\cr
        \vrule height0pt width 2em\cr
        \leftarrowfill\cr
        \hphantom{$\scriptstyle\mspace{8mu}{#2}\mspace{8mu}$}\cr
        \noalign{\kern-0.3ex}
      }%
    }%
  }\limits^{#1}_{#2}}%
}
\makeatother

\DeclareMathOperator{\Kan}{Kan}
\DeclareMathOperator{\Hom}{Hom}
\DeclareMathOperator{\Img}{Im}

\DeclareMathOperator{\Ann}{Ann}%Annihilator
\DeclareMathOperator{\id}{id}
\DeclareMathOperator{\sign}{sgn}
\newcommand{\KK}{\mathbb{K}}

\newcommand{\Horn}[2]{\Lambda^{#1}_{#2}}
\DeclareMathOperator{\Cosk}{Cosk}
\DeclareMathOperator{\Sk}{Sk}
\DeclareMathOperator{\cosk}{cosk}
\DeclareMathOperator{\sk}{sk}
\DeclareMathOperator{\tr}{tr}

\newcommand{\huaB}{\mathcal{B}}%{{\mathcal{E}}}%{\mathcal{B}}
\newcommand{\huaS}{\mathcal{S}}
\newcommand{\huaL}{\mathcal{L}}

\newcommand{\huaG}{\mathcal{G}}

\newcommand{\huaU}{\mathcal{U}}
\newcommand{\huaV}{\mathcal{V}}
\newcommand{\huaW}{\mathcal{W}}
\newcommand{\huaX}{\mathcal{X}}
\newcommand{\huaY}{\mathcal{Y}}

\newcommand{\huaC}{{\mathcal{C}}}%{\mathcal{C}}
\newcommand{\huaD}{\mathcal{D}}

\newcommand{\huaK}{\mathcal{K}}

\newcommand{\huaT}{\mathcal{T}}

\newcommand{\unit}{\mathsf{1}}
\newcommand{\ChnV}{\Ch_{\ge 0}}

\title{Duals of Higher Vector Spaces}

\author{Stefano Ronchi}
\address{Mathematics Institute\\Georg-August-University of G\"ottingen\\Bunsenstra{ss}e 3-5\\G\"ottingen 37073\\Germany}
\email{stefano.ronchi@mathematik.uni-goettingen.de}

\author{Chenchang Zhu}
\address{Mathematics Institute\\Georg-August-University of G\"ottingen\\Bunsenstra{ss}e 3-5\\G\"ottingen 37073\\Germany}
\email{chenchang.zhu@mathematik.uni-goettingen.de}

\date{\today}

\begin{document}

\begin{abstract}
    We introduce a notion of ``$n$-dual'' to a simplicial vector space for $n\ge 0$. Coming with it, there is  a canonical pairing,  which we show to be non-degenerate up to homotopy for homotopy $n$-types. As a result this notion of duality is reflexive up to homotopy for $n$-types. 
    In particular the same properties hold for $n$-groupoid objects in vector spaces, whose $n$-duals are again such $n$-groupoid objects. 
    We study this construction in the context of the Dold-Kan correspondence and we reformulate the Eilenberg-Zilber theorem, which classically controls monoidality of the Dold-Kan functors, in terms of internal homs. We compute explicitly the 1-dual of a groupoid object and the 2-dual of a 2-groupoid object in the category of vector spaces. 
    As the 1-dual of a groupoid object, we recover its dual as a $\VB$ groupoid over a point. 
\end{abstract}
    
    
\maketitle

\setcounter{tocdepth}{2}    
\tableofcontents

\section{Introduction}
In this article, we introduce a notion of ``$n$-dual'' for higher vector spaces. This concept arises naturally in our attempt to give an explicit description of cotangent bundles for higher differentiable stacks modeled by higher Lie groupoids.

One of the central ingredients in symplectic geometry and classical mechanics is the {\bf cotangent} bundle $T^*M$ of a manifold $M$. 
Many other symplectic structures appear as its symplectic reduction. When we go higher,  in the algebraic geometric setting, \cite{PaToVaVe13} and \cite{Calaque19} show that the shifted cotangent bundle of a derived higher stack also carries a shifted symplectic structure. Such shifted symplectic structures provide the correct framework for higher symplectic geometry, enabling the use of the AKSZ and BV formalisms to construct sigma models in higher dimensions. They also unify geometric structures arising in Calabi-Yau geometry and derived moduli spaces, making them one of the most conceptually transformative developments of modern symplectic geometry.

However, unlike in algebraic geometry, where the cotangent bundle can be defined directly as the spectrum of the symmetric algebra of the tangent sheaf, the situation for differentiable stacks modeled by Lie $n$-groupoids is considerably more delicate. The corresponding cotangent bundle {\bf cannot} be obtained by simply taking levelwise duals. Already for $n=1$, if $G := G_1 \rightrightarrows G_0$ is a Lie groupoid, the cotangent groupoid, first constructed in \cite{CosteDazordWeinstein1987},  is given by $T^* G_1 \rightrightarrows A^*$, where $A$ is the Lie algebroid of $G$, rather than by the naive levelwise duals $T^*G_i$. 

This observation points to the need for a systematic notion of ``$n$-dual'' for higher vector spaces. Then, as the tangent bundle of $\huaX$ is encoded by the tangent $\VB$ $n$-groupoid $T\huaX$,  the cotangent bundle will be encoded by a corresponding $\VB$ $n$-groupoid by taking the $n$-dual of $T\huaX$ (see the follow-up article  \cite{CuecaRonchi-temp} and the thesis \cite{stefano-thesis}).

Considering the application to higher Lie groupoids,  we use $\VS$ $n$-groupoids as models for higher vector spaces. These are $n$-groupoid objects in $\Vect$, that is,  simplicial vector spaces equipped with certain Kan conditions \cite{BehrendGetzler2017, Duskin2001/02, Henriques2008, Pridham2013, Zhu2009}\footnote{In fact, a simplicial object in a category of groups is automatically a Kan complex \cite[Thm. 3]{Moore1954}. Thus a simplicial vector space is already a $\VS$ $\infty$-groupoid.
}. In particular,  a $\VS$ 1-groupoid is precisely a 2-vector space in the sense of \cite{BaezCrans2004} (see Example \ref{ex:MooreFillers1Gpd}).

The category of simplicial vector spaces $\SVect$ is well known to be closed monoidal \cite{GoerssJardine2009, Riehl2014}\footnote{It is actually so as a category enriched over itself, as we recall in Section \ref{sec:Svect-int-hom}.}. Thus, to take the dual $\huaV^*$ of an object $\huaV \in \SVect$, a natural way is to take $\huaV^*:=\IHom(\huaV, \mathbbm{1})$, where $\IHom$ is the internal hom and $\mathbbm{1}$ is the monoidal unit.  
In $\SVect$, the monoidal unit is $\mathbbm{1} =B^0\R$ (see \eqref{eq:def-BnR}); and
the internal hom is given by the mapping space, which is defined at each level by\footnote{As a convention, we denote the hom functor of a category by the name of that category.}
\begin{equation*}
    \IHom(\huaV, \huaW)_m = \SVect(\huaV \otimes \Delta[m], \huaW).
\end{equation*}
However, if we follow this strategy and define the dual of $\huaV \in \SVect$ as $\IHom(\huaV, B^0\R)$, we recover (see Example \ref{ex:0dual}) only the vector space $\pi_0(\huaV)^*$ with $\pi_0(\huaV)$ the 0-th homotopy group\footnote{By using the Dold-Kan correspondence one has $\pi_0(\huaV)\cong H_0(N(\huaV))$, the 0-th homology of the normalized complex, which in this case is a vector space.} of $\huaV$.
Therefore, this causes a {\bf loss} of higher categorical and homotopy information in general, and the usual expectation of reflexivity of a dual cannot be met. 
Moreover, as in the vector space case, we define a pairing by applying the adjunction 
\begin{equation}\label{eq:0-pairing}
    \SVect(\IHom(\huaV, B^0\R), \IHom(\huaV, B^0\R)) \xrightarrow{\tau} \SVect(\IHom(\huaV, B^0\R) \otimes \huaV, B^0\R)
\end{equation}
to the identity of the dual. 
This pairing is non-degenerate only if $\huaV$ is a 0-groupoid, i.e. a usual vector space;  and it is non-degenerate up to homotopy\footnote{In the same sense of shifted symplectic structures \cite{CuecaZhu2023, Lesdiablerets} (Definition \ref{def:homological-non-deg}).} only if $\huaV$ is a 0-type (Definition \ref{def:order-type}),  i.e. homotopy equivalent to a usual vector space (Theorem \ref{thm:ndual-pairing-hom-nondeg}).

Therefore we need to give a more meaningful definition of a dual $\VS$ $n$-groupoid that preserves the higher information. This is the main focus of this article. 
More precisely, we define an {\bf $n$-dual} of $\huaV \in \SVect$ by 
\begin{equation*}
    \huaV^{n*} := \IHom(\huaV, B^n\R), \quad \forall n\ge 0, 
\end{equation*}
where $B^n\R$ is as in \eqref{eq:def-BnR}.
Then $\huaV^{n*}$ has {\bf all} the desired properties. First of all, $\huaV^{n*}$ is always a $\VS$ $n$-groupoid (Theorem \ref{thm:VSMappingSpace-nGPD}). Secondly, similar to \eqref{eq:0-pairing}, we can define a canonical $n$-shifted dual pairing (Def. \ref{def:IMPairing}), which is proven to be homologically non-degenerate if and only if $\huaV$ is at most an $n$-type, i.e. it is homotopy equivalent to a $\VS$ $m$-groupoid for $m\le n$ (Theorem \ref{thm:ndual-pairing-hom-nondeg}).
Then, in this case, the twice iterated $n$-dual $(\huaV^{n*})^{n*}$ is homotopy equivalent to the original space $\huaV$ (Theorem \ref{eq:double-W.E.}).
Notably, if $\huaV^{n*}$ is an $n$-groupoid on the nose, its $n$-dual is also one. Thus this $n$-dual is internal to the subcategory of $\VS$ $n$-groupoids. However, we cannot push our result further:  the category of $\VS$ $n$-groupoids is not a rigid monoidal category despite the internal dual, because it is not closed under the monoidal product of $\SVect$. This was already noticed in \cite[Remark 8.5]{HoyoTrentinaglia2023} in the $\VB$ groupoid case.

In general, this non-closedness is a consequence of the Eilenberg-Zilber theorem \cite{EilenbergZilber1953,EilenbergMacLane1954}. With this, we show that the order of a monoidal product satisfies a certain additive lower bound, and thus cannot be internal to $\VS$ $n$-groupoids (Theorem \ref{thm:order-of-tensors}). Another consequence of this theorem is the non-degeneracy of the $n$-dual pairing for $n$-types (Theorem \ref{thm:ndual-pairing-hom-nondeg}). This however requires us to prove a version of the Eilenberg-Zilber theorem in terms of the internal hom (Theorem \ref{thm:EZ-thm-hom}).

The {\bf biggest} technical difficulty is to calculate $n$-duals explicitly. This involves solving an exponentially growing amount of linear equations. In this article, we give an explicit solution for $n\le 2$ in Section \ref{sec:computations}. The 0-dual is our first guess for the dual discussed in the previous paragraphs, and it coincides with the usual dual for vector spaces. The 1-dual is consistent with Pradines's construction of dual $\VB$ groupoids \cite{Pradines1988}. The 2-dual is totally {\bf new} and it allows us to construct explicitly the  2-cotangent bundle of a Lie 2-groupoid with a 2-shifted symplectic structure in \cite{stefano-thesis} and the upcoming work \cite{CuecaRonchi-temp}. 

Finally, we remark that, in theory,  there is another possible way of taking  $n$-duals. Given a simplicial vector space $\huaV$, we can take an $n$-shifted dual of its normalized complex and then reassemble it into a simplicial vector space by the Dold-Kan correspondence. The drawback of this approach is that it becomes non-canonical and depends on choices of connections once we pass to the bundle version \cite{stefano-thesis,CuecaRonchi-temp}.  Moreover, the natural 2-shifted paring is generally not multiplicative (Prop. \ref{prop:DK2shiftedDual-counterex}), and thus cannot extend to a simplicial pairing.

\subsection*{Acknowledgements} We would like to thank Anton Alekseev, Miquel Cueca, Ilias Ermeidis, Matias del Hoyo, Camille Laurent-Gengoux, Madeleine Jotz, Kalin Krishna, Jon Pridham, Giorgio Trentinaglia, Luca Vitagliano, and Hao Xu for fruitful discussions, ideas and suggestions on previous iterations of this work. We also thank all the other members of the Higher Structures Seminar in G\"{o}ttingen for many inspiring conversations. This work was partially supported by DFG grant 446784784, DFG grant ZH 274/5-1,  and RTG 2491.

\section{Preliminary: Higher groupoids and simplicial vector spaces} 

In this paper, we will mostly focus on finite-dimensional vector spaces over $\R$, but $\R$ can be replaced any other field $\mathbb{K}$ of characteristic zero. We denote the category of such vector spaces by $\Vect$.   Most of the material we discuss works for simplicial objects in other abelian categories such as abelian groups and modules over a more general ring. 
A notable exception is Proposition \ref{prop:equivalences-in-svect}. There, because in a category of modules over a ring that is not a field, not all sequences split, the notion of weak equivalence is a strictly weaker notion than homotopy equivalence.

\subsection{\texorpdfstring{$n$}{n}-groupoid objects}\label{sec:ngpds-intro}

Here we first recall some fundamental facts and definitions on $n$-groupoids and  $\VS$ $n$-groupoids. In this article,  we assume our category $\Cat$ to be complete, cocomplete, and have a forgetful functor $U: \Cat \to \Set$ which preserves limits. All ``algebraic'' categories as referred to by Goerss-Jardine \cite*[Example II.2.8]{GoerssJardine2009}: groups, rings, modules, algebras or Lie algebras are examples of such categories. In this article we are mainly interested in $\Set$ itself and $\Vect$. We occasionally refer to examples of Lie $n$-groupoids and $\VB$ $n$-groupoids, which we will study in more detail in \cite{CuecaRonchi-temp}. See Remark \ref{rmk:cat-assumptions-2} for further comments on $\Cat$.

Let $\Delta$ be the category of finite ordinals, denoted by $[0] = \{0\}, [1] = \{0,1\}, \dots, [m] = \{0,1,\dots,m\}, \dots$, with order-preserving maps. 
A \textbf{simplicial object} in a category $\Cat$ is a contravariant functor from $\Delta$ to $\Cat$. A \textbf{simplicial map} between two simplicial objects in $\Cat$ is a natural transformation between the corresponding contravariant functors.  We denote a simplicial object by a calligraphic letter, for example $\huaX: \Delta^{op} \to \Cat$. 
More concretely, $\huaX$ is a tower of objects $\huaX_m$ in $\Cat$, with face maps $d^m_j: \huaX_m \to \huaX_{m-1}$ and degeneracy maps $s^m_j: \huaX_{m} \to \huaX_{m+1}$ between them for any $m\ge 0$ and any $0\le j \le m$. These maps satisfy the following simplicial identities:
\begin{equation}\label{eq:face-degen}
    \begin{array}{lll}
        d^{l-1}_i d^{l}_j = & d^{l-1}_{j-1} d^l_i &\text{if}\; i<j,  \\
       s^{l}_i s^{l-1}_j =& s^{l}_{j+1} s^{l-1}_i & \text{if}\; i\leq j,
    \end{array}\qquad  d^l_i s^{l-1}_j =\left\{\begin{array}{ll}
    s^{l-2}_{j-1} d^{l-1}_i  & \text{if}\; i<j, \\
    \id  & \text{if}\; i=j, j+1,\\
    s^{l-2}_j d^{l-1}_{i-1} & \text{if}\; i> j+1.
 \end{array}\right.
\end{equation}
We drop the upper indices on the simplicial structure maps, when they are clear from context. Then a simplicial map between two simplicial objects \(\huaX\) and~\(\huaY\) is a family of maps \(\huaX_k\to \huaY_k\) that intertwine their face and degeneray maps. We denote by $\Simp\Cat$ the category of simplicial objects and simplicial maps in $\Cat$. A simplicial set is a simplicial object in $\Set$, while a simplicial vector space is a simplicial object in $\Vect$. 

\begin{remark}\label{rmk:hom-C}
For simplicial sets $\huaS, \huaT \in \SSet$, the set of simplicial maps $\hom(\huaS, \huaT)=\SSet(\huaS, \huaT)$ is also an equalizer,
\begin{equation}\label{eq:hom-equalizer}
    \begin{tikzcd}
        \hom(\huaS, \huaT) \ar[r] & \prod_{l\ge 0} \SSet(\huaS_l, \huaT_l) \ar[r, shift left=1ex, "a"] \ar[r, shift right=1ex, "b"] & \prod_{g:[m]\to [k]} \SSet(\huaS_{k}, \huaT_m)^g. 
    \end{tikzcd} 
\end{equation}Here for $f=(f_l)_{l\ge 0} \in \prod_{l\ge 0} \SSet(\huaS_l, \huaT_l)$, 
\[
    a(f)^{g:[m]\to [k]} = \huaS(g)\circ f_m: \huaS_k \to \huaT_m, \quad b(f)^{g:[m]\to [k]} = f_k \circ \huaT(g) : \huaS_k \to \huaT_m. 
\]

For $S\in \Set$ and $X\in \Cat$, we shall write $\hom(S, X) := \prod_{s\in \huaS} X^s$, that is copies of $X$ indexed by elements of $S$. Thus we define $\hom(\huaS, \huaX)$ for $\huaS \in \SSet$ and $\huaX \in \Simp\Cat$ as the following equalizer
\begin{equation}
    \begin{tikzcd}
        \hom(\huaS, \huaX) \ar[r] & \prod_{l\ge 0} \hom(\huaS_l, \huaX_l) \ar[r, shift left=1ex, "a' "] \ar[r, shift right=1ex, "b' "] & \prod_{g:[m]\to [k]} \hom(\huaS_{k}, \huaX_m)^g, 
    \end{tikzcd}
\end{equation} with
\begin{equation*}
\begin{aligned}
    &(a')^{g:[m]\to [k]}: ((\huaX_l^{s_l})_{s_l\in \huaS_l})^{l\ge 0} \mapsto (\huaX_m^{\huaS(g)(s_m)})_{s_m \in \huaS_m},\\
    &(b')^{g:[m]\to [k]}: ((\huaX_l^{s_l})_{s_l\in \huaS_l})^{l\ge 0} \mapsto  (\huaX(g)(\huaX_k))^{s_m}_{s_m \in \huaS_m}. 
\end{aligned}
\end{equation*}
Thus if $\Cat$ is complete as we assume in this article, $\hom(\huaS, \huaX)$ is another element in $\Cat$, otherwise, one may view it as a presheaf over $\Cat$.

As shown in \cite[Prop.4.5]{RogersZhu2020}, for a test element $T\in \Cat$, we have
\begin{equation} \label{eq:copower-adjoint}
    \Cat(T,  \hom(\huaS, \huaX) )=\Simp\Cat(\huaS \otimes T, \huaX),
\end{equation} where  $\huaS \otimes T \in \Simp\Cat$ is the copowering of $\huaS$ with the constant simplicial object $T$ (see also \eqref{eq:copowering}), that is each level $(\huaS \otimes T)_k = \coprod_{s\in \huaS_k} T^s$. 

By \eqref{eq:copower-adjoint},  $\hom(-, \huaX): \SSet^{op} \to \Cat$ is a right adjoint functor. Thus it preserves limit, or equivalently as a contravariant functor $\hom(-, \huaX): \SSet \to \Cat$, it turns colimt to limit. 
\end{remark}

There are three classes of simplicial sets which play an important role in our article. These are the $m$-simplex $\Delta[m]$, the boundary of the $m$-simplex $\partial\Delta[m]$ and the $(m,k)$-horn $\Lambda[m,k]$. They are given for each $m \ge 0$, $0 \le l \le m$, as 
\begin{equation}\label{eq:simplex-horn}
\begin{split}
(\Delta[m])_l & = \{ f: [l] \to [m] \mid f(i)\leq
f(j),
\forall i \leq j\}, \\
(\partial\Delta[m])_l &=  \{ f\in (\Delta[m])_l \mid \{0,\dots, m\}
\nsubseteq \{ f(0),\dots, f(l)\} \}\\
(\Lambda[m,k])_l & = \{ f \in (\Delta[m])_l \mid \{0,\dots,\widehat{k},\dots,m\}
\nsubseteq \{ f(0),\dots, f(l)\} \}\\
\end{split}
\end{equation}
Geometrically, the boundary of the $m$-simplex is obtained from it by removing its interior, which is the unique non-degenerate $m$-simplex in the simplicial set $\Delta[m]$. Similarly, the $(m,k)$-horn is obtained from the boundary by further removing the $k$-th face (i.e. the one opposite the vertex $k$). Note that in the simplicial description, removing these non-degenerate simplices also requires removing all of  their higher-dimensional degeneracies. 
Clearly then $\Lambda[m,k] \subseteq \partial\Delta[m] \subseteq \Delta[m]$. 

In these cases, the limit defining $\hom(\huaS, \huaT)$ in \eqref{eq:hom-equalizer} simplifies considerably by choosing a sub-index-category with the same limit. For $\huaS=\Delta[m]$, one can choose the sub-index-category with only one object: the non-degenerate $m$-simplex in $\Delta[m]$ and the identity arrow. Thus for $\huaX \in \Simp\Cat$, $\hom(\Delta[m], \huaX) = \huaX_m$, and we denote the sets of $m$-boundaries and $(m,k)$-horns in $\huaX$ by 
\begin{equation*}
    \begin{array}{cc}
    \partial^m(\huaX) := \hom(\partial\Delta[m], \huaX),
    &\Lambda^m_k(\huaX) := \hom(\Lambda[m,k],\huaX), 
    \end{array}
\end{equation*}
respectively. 
Both of these objects can be written as fiber products of $\huaX_{m-1}$ over lower dimensional horns and boundaries\footnote{For a proof of the general fiber product description of $\Lambda^m_k(\huaX)$ for an $\infty$-groupoid and a formula computing its dimension see \cite{stefano-thesis}. See also \eqref{eq:partial}. }.

For $\Cat=\Set$ or $\Vect$, we have the following explicit elementwise descriptions of the sets of $m$-boundaries and $(m,k)$-horns in $\huaX$:
\begin{equation}\label{eq:DefHornBdrySpaces}
    \begin{split}
    \partial^m(\huaX) &= \{(x_0, \dots, x_m) \in (\huaX_{m-1})^{\times m+1} \mid d_ix_j = d_{j-1}x_i, \text{ for all } i < j\},\\
    \Lambda^m_k(\huaX) &= \{(x_0, \dots, \widehat{x_k}, \dots, x_m) \in (\huaX_{m-1})^{\times m} \mid \\
    &\qquad \qquad d_ix_j = d_{j-1}x_i, \text{ for all } i < j, \text{ such that } i\neq k, j\neq k. \}. 
    \end{split}
\end{equation}

More intuitively, $\partial^m(\huaX)$ is the space of all possible boundaries of $m$-simplices in $\huaX$. In the same way, $\Lambda^m_k(\huaX)$ is the space of all possible configurations of $(m-1)$-simplices in the shape of an $(m,k)$-horn that exist in $\huaX$. That is, pictorially,
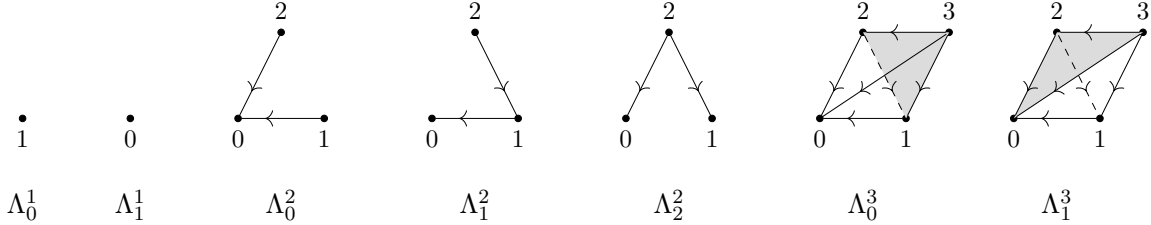
\begin{center}
\begin{adjustbox}{width=\textwidth}
\pgfdeclarelayer{nodelayer}
\pgfdeclarelayer{edgelayer}
\pgfsetlayers{nodelayer,main,edgelayer}
\begin{tikzpicture}[scale=0.6,
    none/.style={}, 
    dot/.style={fill=black, draw=black, shape=circle, scale=0.25}, 
    fill-grey/.style={-, fill={rgb,255: red,220; green,220; blue,220}, draw=none, fill opacity=0.5},
    directed/.style={decoration={markings,mark=at position 0.7 with {\arrow{To[scale=1.25]}}}, postaction={decorate}},
    dashdir/.style={dashed, decoration={markings,mark=at position 0.7 with {\arrow{To[scale=1.25]}}}, postaction={decorate}},
    ]
	\begin{pgfonlayer}{nodelayer}
		\node [style=none] (25) at (20.5, 1) {};
		\node [style=none] (26) at (19.5, 3) {};
		\node [style=none] (27) at (21.5, 3) {};
		\node [style=none] (0) at (0, -1) {$\Lambda^1_0$};
		\node [style=none] (1) at (2.5, -1) {$\Lambda^1_1$};
		\node [style=none] (2) at (6, -1) {$\Lambda^2_0$};
		\node [style=none] (3) at (10.5, -1) {$\Lambda^2_1$};
		\node [style=none] (4) at (15, -1) {$\Lambda^2_2$};
		\node [style=none] (5) at (19.5, -1) {$\Lambda^3_0$};
		\node [style=none] (6) at (24, -1) {$\Lambda^3_1$};
		\node [style=dot, label={below:\small{1}}] (7) at (0, 1) {};
		\node [style=dot, label={below:\small{0}}] (8) at (2.5, 1) {};
		\node [style=dot, label={below:\small{0}}] (9) at (5, 1) {};
		\node [style=dot, label={below:\small{1}}] (10) at (7, 1) {};
		\node [style=dot, label={above:\small{2}}] (11) at (6, 3) {};
		\node [style=dot, label={below:\small{0}}] (12) at (9.5, 1) {};
		\node [style=dot, label={below:\small{1}}] (13) at (11.5, 1) {};
		\node [style=dot, label={above:\small{2}}] (14) at (10.5, 3) {};
		\node [style=dot, label={below:\small{0}}] (15) at (14, 1) {};
		\node [style=dot, label={below:\small{1}}] (16) at (16, 1) {};
		\node [style=dot, label={above:\small{2}}] (17) at (15, 3) {};
		\node [style=dot, label={below:\small{0}}] (18) at (18.5, 1) {};
		\node [style=dot, label={below:\small{1}}] (19) at (20.5, 1) {};
		\node [style=dot, label={above:\small{2}}] (20) at (19.5, 3) {};
		\node [style=dot, label={above:\small{3}}] (21) at (21.5, 3) {};
		\node [style=none] (22) at (23, 1) {};
		\node [style=none] (23) at (24, 3) {};
		\node [style=none] (24) at (26, 3) {};
		\node [style=dot, label={below:\small{0}}] (28) at (23, 1) {};
		\node [style=dot, label={below:\small{1}}] (29) at (25, 1) {};
		\node [style=dot, label={above:\small{2}}] (30) at (24, 3) {};
		\node [style=dot, label={above:\small{3}}] (31) at (26, 3) {};
	\end{pgfonlayer}
	\begin{pgfonlayer}{edgelayer}
		\draw [style=fill-grey] (27.center)
			 to (25.center)
			 to (26.center)
			 to cycle;
		\draw [style=directed] (10) to (9);
		\draw [style=directed] (11) to (9);
		\draw [style=directed] (13) to (12);
		\draw [style=directed] (14) to (13);
		\draw [style=directed] (17) to (15);
		\draw [style=directed] (17) to (16);
		\draw [style=directed] (20) to (18);
		\draw [style=directed] (19) to (18);
		\draw [style=directed] (21) to (20);
		\draw [style=directed] (21) to (19);
		\draw [style=dashdir] (20) to (19);
		\draw [style=fill-grey] (23.center)
			 to (24.center)
			 to (22.center)
			 to cycle;
		\draw [style=directed] (30) to (28);
		\draw [style=directed] (29) to (28);
		\draw [style=directed] (31) to (30);
		\draw [style=directed] (31) to (29);
		\draw [style=dashdir] (30) to (29);
		\draw [style=directed] (31) to (28);
		\draw [style=directed] (21) to (18);
	\end{pgfonlayer}
\end{tikzpicture}
\end{adjustbox}
\end{center}
where our convention is to orient arrows from the higher to the lower vertex number and the shaded faces denote empty triangles. 

\begin{definition} \label{def:n-gpd} 
    For $n\in \Z^{\ge 0}\sqcup \infty$, an \textbf{$n$-groupoid} $\huaX$ is a simplicial set satisfying the Kan conditions $\Kan(m,k)$ for any $m\ge 1$, $0\le k\le m$ and the strict Kan conditions $\Kan!(m,k)$ for any $m \ge n+1$, $0 \le k \le m$. 
    These conditions are
    \begin{itemize}[leftmargin=*]
        \item $\Kan(m,k)$: The horn projection $p^m_k: \hom(\Delta[m],\huaX) \to \hom(\Lambda[m,k],\huaX)$ is surjective.
        \item $\Kan!(m,k)$: The horn projection $p^m_k: \hom(\Delta[m],\huaX) \to \hom(\Lambda[m,k],\huaX)$ is bijective.
    \end{itemize}
    An \textbf{$\infty$-groupoid} is also known as a \textbf{Kan complex}.  An \textbf{$n$-group} is an $n$-groupoid $\huaX$ with $\huaX_0 = pt$.
\end{definition}

Intutively, the Kan condition $\Kan(m,k)$ states the following: 
For any configuration of $(m-1)$-simplices in $\huaX$ forming an $(m,k)$-horn, that is, an $m$-tuple as in \eqref{eq:DefHornBdrySpaces}, there exists a unique \textbf{horn filler}, that is an $m$-simplex $x \in \huaX_m$, such that $d_i x = x_i$ for any $i\neq k$. 

\begin{remark}\label{rmk:cat-assumptions-2}
    These definitions carry over directly to the category $\Vect$ thanks to Remark \ref{rmk:hom-C}. In fact, $n$-groupoid objects make sense for any category $\Cat$ together with a special class of maps called covers\footnote{A choice of such a class of maps gives rise to a Grothendieck pretopology when certain conditions are satisfied.} \cite[Def.1.2]{Zhu2009}, by replacing surjective maps with covers in $\Kan(m,k)$ and bijective maps with isomorphisms in $\Kan!(m,k)$. When treating $n$-groupoid objects in a general category $\Cat$ in this article, we further require that the forgetful functor $U: \Cat\to \Set$ preserves covers and reflects isomorphisms: that is $U$ maps a cover to a surjective map, and if $U(f)$ is an isomorphism then $f$ is an isomorphism.  For instance, $\Set$ with surjective maps, $\Vect$ with surjective maps, and other algebraic categories with surjective maps, are examples of categories with such a choice of covers. Notably, the category of manifolds $\Mfd$ with surjective submersions and the category of vector bundles $\VB$ with levelwise surjective submersions do not satisfy these assumptions: the forgetful functor does not reflect isomorphisms and one needs to use more complex methods, such as described in \cite{Zhu2009, RogersZhu2020}.
\end{remark}

For $n$-groupoid objects in $\Vect$ with surjective maps, we have the following equivalent simplified definition: 

\begin{definition}\label{def:vs-n-gpd}
    A \textbf{$\VS$ $n$-groupoid} $\huaV$ is a simplicial vector space whose underlying set is an $n$-groupoid. A {\bf $\VS$ $n$-group}  $\huaV$ is a $\VS$ $n$-groupoid with $\huaV_0=0$. 
\end{definition}

\begin{remark}[Relation with $\VB$ $n$-groupoids] \label{ex:vb-vs-gpd}
    Simplicial vector bundles and $\VB$ $n$-groupoids were studied in \cite{HoyoTrentinaglia2021, HoyoTrentinaglia2023,HoyoTrentinaglia2025} to generalize $\VB$ groupoids \cite{Pradines1988, Mackenzie2005, GraciaSazMehta2017}.
    As we pointed out just before Definition \ref*{def:vs-n-gpd}, $\VB$ $n$-groupoids are $n$-groupoid objects in the category of vector bundles with surjective submersions as covers in the sense of \cite{Zhu2009}. In particular, this implies that the base of a $\VB$ $n$-groupoid is a Lie $n$-groupoid, which is in turn an $n$-groupoid object in the category of manifolds with surjective submersions.

    Simplicial vector spaces are then simplicial vector bundles over the identity groupoid of a point. Similarly, $\VS$ $n$-groupoids are also just $\VB$ $n$-groupoids over a point. More generally, a $\VB$ $n$-groupoid over a manifold $M$ (viewed as a constant simplicial manifold which is furthermore a Lie $n$-groupoid) is effectively ``a bundle of $\VS$ $n$-groups'' over $M$.  Observe that for any simplicial manifold $\huaG$ there is an embedding via the total unit 
    \begin{equation} \label{eq:total-deg}
       \unit:= s_{i_k}\circ \dots \circ s_{i_0}: \huaG_0 \to \huaG_k,
    \end{equation} 
    for $i_k, \dots, i_0$ arbitrary  indices. In fact,  any combination of indices defines the same map. Then a $\VB$ $n$-groupoid $\huaV\to \huaG$ may be pulled back via $1$ to $\huaG_0$ and this results in a bundle $1^*\huaV$ of $\VS$ $n$-groups over $\huaG_0$. 
\end{remark}

\begin{remark}[Front-to-back duality {\cite[\S 8.2.10]{Weibel1994}}{\cite[\href{https://kerodon.net/tag/003L}{Tag 003L}]{kerodon}}]\label{rem:front-to-back}
    Every simplicial object $\huaX$ in a category $\Cat$ has an opposite $\huaX^{op}$, which consists of the same spaces but has simplicial maps $(d^{op})^n_i = d^n_{n-i}$, $(s^{op})^n_i = s^n_{n-i}$. 
    This can be understood as reversing the order of the elements in each object of the category $\Delta$, or equivalently reversing the direction of all the edges and orienting them from lower to higher vertex number. 
    When the simplicial object is the nerve of a category, its opposite is precisely the nerve of the opposite category \cite[\href{https://kerodon.net/tag/003Q}{Tag 003Q}]{kerodon}. 
    The existence of opposite simplicial objects extends the well-known principle of categorical duality (see e.g. \cite[\S 1.2]{Riehl2016}) to this context.
    In particular, if a property holds for an arbitrary simplicial object in $\Cat$, then it must hold for its opposite as well. This is helpful in certain proofs. 
    We will make use of this property in the computation of 1- and 2-duals in Section \ref{sec:computations}. 
\end{remark}

\subsection{Moore Fillers for horns}\label{sec:MooreFillers}
Simplicial vector spaces are in particular simplicial groups, and it is a general fact that these are always $\infty$-groupoids, because one can construct fillers for all horns by using the group operations. This result is due to Moore \cite[Thm. 3 p. 18-04]{Moore1954}, see also \cite[Thm. 17.1]{May1967}, \cite[Lem. 3.1]{Curtis1971}, \cite[Lem. 8.2.8]{Weibel1994}, \cite[Lem. I.3.4]{GoerssJardine2009}. We recall it here with an explicit proof in the case of vector spaces, adapted from \cite{nlab:simplicial_group}, as it will be used in our later discussion.

\begin{proposition}[Moore]\label{thm:MooreHornFillers}
   Simplicial vector spaces are $\VS$ $\infty$-groupoids. 
\end{proposition}
\begin{proof}
    Let $\huaV$ be a simplicial vector space, and $(x_0, \dots, \widehat{x_{k}}, \dots, x_n) \in \Lambda^n_k(\huaV)$, an $(n,k)$-horn in $\huaV$.
    That is, $d_i x_j = d_j x_{i+1}$ for any $i\ge j$, whenever both sides are defined. 
    We describe an algorithm to construct a horn filler $x \in \huaV_{n+1}$, i.e. $x$ is such that $d_ix = x_i$ for any $i\neq k$. 
    Let us start with the case of $0<k<n$.
    \begin{enumerate}
    \item[\underline{Step 1}] We construct by induction on $i\in [0, k-1]$ a $w_{k-1}$ such that $d_jw_{k-1} = x_j$ for $j\in[0, k-1]$. \\
    We begin with $w_0 = s_0x_0 \in \huaV_{n+1}$. Then by the simplicial identities, $d_0w_0 = x_0$. Now assume that $w_{i-1}$ has been defined in a way such that $d_jw_{i-1} = x_j$ for any $j\in[0, i-1]$. We then define 
    \begin{equation} \label{eq:wi}
        w_i := w_{i-1} - s_id_iw_{i-1} + s_ix_i.
    \end{equation}
    A simple computation using the simplicial identities shows $d_iw_i = x_i$. 
    Moreover, for any $0 \le j \le i-1$, 
    \begin{equation}\label{eq:MooreThmComp1}
        \begin{split} 
            d_jw_i &= d_jw_{i-1} - d_js_id_iw_{i-1} + d_js_ix_i\\
            &= d_jw_{i-1} - s_{i-1}d_jd_iw_{i-1} + s_{i-1}d_jx_i\\
            &= d_jw_{i-1} - s_{i-1}d_{i-1}d_jw_{i-1} + s_{i-1}d_{i-1}x_j = d_jw_{i-1} = x_j.
        \end{split}
    \end{equation}
Thus the inductive step is proven. Hence $w_{k-1}$ has the property that $d_iw_{k-1}=x_i$ for $0\le i \le k-1$.
    \item[\underline{Step 2}] We define $w_n = w_{k-1} - s_{n-1}d_nw_{k-1} + s_{n-1}x_n$. 
    With a similar calculation to the one after \eqref{eq:wi} and that in \eqref{eq:MooreThmComp1}, $d_n w_n = x_n$ and $d_iw_n = x_i$ for any $0\le i \le k-1$.
    \item[\underline{Step 3}] We proceed by induction backwards, that is for all $i\in [k+1, n]$ we construct a $w_i$ such that $d_jw_i = x_j$ for $j\in[0, k-1] \cup [i, n]$. \\
    The initial case ($i=n$) was verified in the last step. Now assume that $w_{i+1}$ has been defined in a way such that $d_j w_{i+1} = x_j$ for any $j$ such that $0 \le j \le k-1$ or $i+1 \le j \le n$. 
    Then we define 
    \begin{equation}
        w_i := w_{i+1} - s_{i-1}d_iw_{i+1} + s_{i-1}x_i.
    \end{equation}
    Again, we have that $d_iw_i = x_i$ and, for any $j$ such that $0 \le j \le k-1$ or $i+1 \le j \le n$, a computation similar to \eqref{eq:MooreThmComp1}, shows that $d_jw_{i+1}=x_j$.
    Then, by induction,  $w_{k+1}$ is our desired horn filler, because by this construction $d_jw_{k+1} = x_j$ for any $j\neq k$. 
\end{enumerate}
    
Now, in the case of $k=n$, we start with $w_0= s_0x_0 \in \huaV_{n+1}$ as above and repeat the first induction procedure until $w_{n-1}$. Then $w_{n-1}$  is already the desired horn filler. 
    
Similarly,  for $k=0$, we start with $w_n=w_0-s_{n-1}d_nw_0+s_{n-1}x_n$. Then $d_nw_n=x_n$ and we proceed by the downward induction as in Step 3 until $w_1$, which is the desired horn filler. 
\end{proof}

\begin{remark}\label{rmk:MooreFillers}
    Given  a simplicial vector space $\huaV$, we define a linear map $\mu^m_k: \Lambda^m_k(\huaV) \to \huaV_m$ by the algorithm in the the proof of Proposition \ref{thm:MooreHornFillers}, and we call $\mu^m_k(\lambda_k) \in \huaV_m$ the \textbf{Moore filler} of $\lambda_k \in \Lambda^m_k(\huaV)$. In general, $\mu^m_k$ is a right inverse of the horn projection $p^m_k: \huaV_m \to \Lambda^m_k(\huaV)$. If $\huaV$ is in particular a $\VS$ $n$-groupoid, then $\mu^m_k$ is \textit{the} inverse of $p^m_k$ when $m\ge n+1$. This makes the structure of a $\VS$ $n$-groupoid particularly rigid, as we will see in Theorem \ref{thm:FiniteDataVS}. The following example illustrates this situation for $n=1$. 
\end{remark}

\begin{example}\label{ex:MooreFillers1Gpd}
    Let $\huaV_1 \rightrightarrows  \huaV_0$ be a category object internal to the category of vector spaces $\Vect$. As proven in \cite{BaezCrans2004}, $\huaV_1 \rightrightarrows  \huaV_0$ is automatically a groupoid internal to $\Vect$, that is, a $\VS$ (1-)groupoid. Such an object is called a {\em 2-vector space} in \cite{BaezCrans2004}.
    
    By Proposition \ref{thm:MooreHornFillers}, since $\mu^2_1$ provides the inverse of  $p^2_1$, we obtain a linear multiplication on a 2-vector space
    \[ m_1=d_1 \mu^2_1  : \Lambda^2_1(\huaV) \to \huaV_1,  \]
    when identifying $\huaV_1 \rightrightarrows  \huaV_0$ with its nerve $\huaV$, which is a simplicial vector space.

More explicitly,  for $u,w \in \huaV_1$ meeting at $x = d_1u = d_0w \in \huaV_0$,
    \begin{equation} \label{eq:mul-Gpd}
        w \cdot u := m_1(u,w) = d_1\mu^2_1(u,w) = u + w - \unit x.
    \end{equation}
    This was also observed in \cite[Lemma 3.2]{BaezCrans2004}. In fact, using \eqref{eq:mul-Gpd}, one easily obtains the inverse of $u$ as $u^{-1} = -u + \unit d_0u + \unit d_1u$.

    We thus notice that a 2-vector space is completely determined by a pair of vector spaces with source, target and unit maps between them, as the multiplication and inversion can be inferred from these. We state a general version of this fact in Theorem \ref{thm:FiniteDataVS}. 

    This inversion also appears in the involution exchanging the right and left core of a $\VB$ groupoid in \cite[\S 3.2]{GraciaSazMehta2017}, \cite[\S 11.2]{Mackenzie2005}, which is $\unit^*\ker \widetilde{d}_0 \to \unit^*\ker \widetilde{d}_1$, $c \mapsto -c^{-1} = c - \widetilde{\unit}\widetilde{d}_0c$, where $(\widetilde{d}_i, d_i)$ and $(\widetilde{\unit}, \unit)$ denote the source, target maps and the unit of the $\VB$ groupoid respectively. 

\end{example}

For a groupoid in general,  the degeneracy $s_0$ serves as the groupoid unit for the groupoid multiplication. This is true for higher $\VS$ groupoids as well. That is, the degeneracy maps give higher ``units'' for the higher multiplications in the language of \cite[\S 2.3]{Zhu2009}.  This holds thanks to   the fact that Moore fillers provide explicit formulas for higher multiplications and   $\mu^m_k p^m_k$ preserves degenerate simplices.  We summarize these in the following lemma. 

\begin{lemma}\label{lem:MooreFillers-degeneracies}
    Let $\huaV$ be a simplicial vector space. The space of degenerate $m$-simplices $D_m\huaV$ is isomorphic to the horn space $\Lambda^m_k(\huaV)$ for any $0 \le k \le m$. In particular the Moore fillers $\mu^m_k$ are compatible with the degeneracy maps in the sense that 
    \begin{equation*}
        \mu^m_k p^m_k s_j = s_j, \qquad  \forall k\in [0, m],\,  \forall j\in [0, m-1].   
    \end{equation*}
\end{lemma}
\begin{proof}
    Consider the following short exact sequences: 
    \begin{equation}
        0\to D_m\huaV \into \huaV_m \to \huaV_m/D_m\huaV \to 0, \qquad 0\to  \ker p^m_m \to \huaV_m  \xrightarrow{p^m_m} \Lambda^m_m(\huaV) \to 0. 
    \end{equation}
    By the isomorphism giving equivalent definitions of normalized chains in the Dold-Kan correspondence (see e.g. \cite[Thm. III.2.1]{GoerssJardine2009}), $\huaV_m/D_m\huaV \cong \ker p^m_m$, so $\dim \Horn{m}{m}(\huaV) = \dim D_m(\huaV)$. In addition, by the formula computing the dimension of horn spaces in \cite[Cor. 1.2.34]{stefano-thesis}, the space of horns $\Horn{m}{k}(\huaV)$ has the same dimension for all $k\in [0, m]$. Thus  $\dim \Horn{m}{k}(\huaV) = \dim D_m\huaV$ for all $k\in [0, m]$.

    On the other hand, by the construction of $\mu^m_k$ in the proof of Proposition \ref{thm:MooreHornFillers}, the image of $\mu^m_k$ is contained in $D_m\huaV$, therefore $p^m_k|_{D_m\huaV}\circ \mu^m_k = \id$. In particular,  $p^m_k|_{D_m\huaV}: D_m \huaV \to \Horn{m}{k}(\huaV) $ is surjective. Since the two vector spaces have the same dimension, $p^m_k|_{D_m\huaV}$ is an isomorphism with inverse $\mu^m_k$. Clearly, then, for any degenerate simplex $s_j v$, $\mu^m_k p^m_k s_j v = s_j v$.
\end{proof}

\subsection{\texorpdfstring{$\Cosk$}{Cosk} functor and finite data for \texorpdfstring{$\VS$}{VS} \texorpdfstring{$n$}{n}-groupoids}\label{sec:FiniteDataVS}

We demonstrate a general way to represent an $n$-groupoid object $\huaX$ with finite data, which is expanded upon in \cite{Zhu2009, Duskin2001/02} with a particular focus on the $n=2$ case. This construction is based on the coskeleton construction, which we now briefly recall.

Consider the \textbf{$m$-truncation} of a simplicial object $\huaX$, which is obtained by ``forgetting'' all simplicial levels higher than $\huaX_m$, and it consists of the objects $\huaX_0, \dots, \huaX_m$ and simplicial maps between them. We denote it by $\tr^m\huaX$. 
This defines a functor $\tr^m$ from the category $\Simp\Cat$ of simplicial objects in $\Cat$ to the category $\Simp\Cat_{\le m}$ of $m$-truncated simplicial objects in $\Cat$, which is defined as the category of functors from the full subcategory $\Delta_{\le m}$ of $\Delta$ generated by the objects $\{[0],[1],\dots [m]\}$ to $\huaC$.
The truncation $tr^m$ admits a left adjoint $\sk^m$ called the \textbf{skeleton} functor and a right adjoint $\cosk^m$ called the \textbf{coskeleton} functor. By composing each of them with the truncation, we obtain a pair of adjoint endofunctors on $\Simp\Cat$, $\Sk^m \vdash \Cosk^m$, where $\Sk^{m}(\huaX) := \sk^m(\tr^m(\huaX))$ is called the $m$-skeleton of $\huaX$ and $\Cosk^{m}(\huaX):= \cosk^{m}(\tr^m(\huaX))$ is the $m$-coskeleton of $\huaX$.
The skeleton functor is defined in such a way that $\Sk^m(\huaX)$ is the sub-simplicial object of $\huaX$ that coincides with it on all levels $l$ for $l\le m$ and it only consists of the degenerate simplices at higher levels.
Meanwhile, the coskeleton can be described as a limit \cite[\href{https://stacks.math.columbia.edu/tag/0AMA}{Section 0AMA}]{stacks-project}. We have the following equivalent explicit inductive formula\footnote{This formula is taken from \cite[\S 2.2-2.3]{Duskin2001/02} which treats the case of simplicial sets. However it is not hard to verify that the general formula in  \cite[\href{https://stacks.math.columbia.edu/tag/0183}{Lemma 0183}]{stacks-project} coincides with this inductive one.} for the coskeleton of a simplicial object $\huaX$:  
\begin{equation*}
    \Cosk^m(\huaX)_l := \begin{cases}
        \huaX_l &\text{for } l \le m,\\
        \partial^{m+1}(\huaX) &\text{for } l = m+1,\\
        \partial^{l}(\tr^{l-1}\Cosk^m(\huaX)) &\text{for } l > m+1.
    \end{cases}
\end{equation*}
Here $\partial^{l}(\huaX)$ denotes the $l$-th boundary space of $\huaX$. Notice that $\partial^{l}(\huaX)=\hom (\partial \Delta^l, \huaX)$ is an equalizer induced by the coequalizer diagram
\begin{equation}\label{eq:partial}
    \coprod_{0\le i, j \le l } \Delta[l-2] \rightrightarrows \coprod_{i=0}^l \Delta[l-1] \to \partial \Delta[l],
\end{equation} given by the relations $d^jd^i=d^id^{j-1}$ for $i<j$ (see \cite[Prop. 2.3]{GoerssJardine2009}).
Thus $\partial^{l}(\huaX)$ is a limit involving only $\huaX_0, \dots, \huaX_{l-1}$ and face maps between them. Therefore it depends only on the $(l-1)$-truncation of $\huaX$ and it is also defined for any $(l-1)$-truncated simplicial object.

\begin{lemma} \label{lem:Cosk}
    For an $n$-groupoid object $\huaX$ in $\Cat$ (as in Remark \ref{rmk:cat-assumptions-2}):
    \begin{enumerate}[label={\arabic*)}]
        \item $U(\huaX)$ is an $n$-groupoid, that is an $n$-groupoid object in $\Set$.
        \item $\Cosk^{n+1}(\huaX) \cong \huaX$.
    \end{enumerate}
\end{lemma}
\begin{proof}
    The first statement is immediate from the fact that the forgetful functor $U$ preserves covers\footnote{Notice that all functors preserve isomorphisms.}.  
    Then since  $\Cosk^{n+1}(\huaY) \cong \huaY$ \cite[\S 2.3]{Duskin2001/02} for an $n$-groupoid object $\huaY$ in $\Set$, we have
    \[ U(\Cosk^{n+1}(\huaX)) \cong \Cosk^{n+1}(U(\huaX)) \cong U(\huaX). 
    \]
    So the second statement follows from the fact that $U$ reflects isomorphisms. 
\end{proof}

In particular, given an $n$-groupoid object $\huaX$, statement (2) above reduces the data of $\huaX$ to its $(n+1)$-truncation $\tr^{n+1}(\huaX)$. Now since $\huaX_{n+1} \cong \Lambda^{n+1}_k(\huaX)$ for any $0 \le k \le n+1$, we can further reduce the data of $\huaX$ to the $n$-truncation $\tr^n(\huaX)$ together with $n+2$ many $(n+1)$-ary \textbf{multiplication maps} $m_k: \Lambda^{n+1}_k(\huaX) \to \huaX_{n}$. These $m_k$'s are the higher counterparts of the multiplication and inverse maps of a groupoid. The compatibility conditions they satisfy shall be understood as higher counterparts of associativity and other laws obeyed by groupoid structure maps. This is discussed in more detail in \cite[\S 2.3]{Zhu2009}.

In the case of $\VS$ $n$-groupoids, we have an even stronger result: a $\VS$ $n$-groupoid $\huaV$ is already determined by its $n$-truncation $\tr^n(\huaV)$ without further information on multiplication maps. This is because the multiplications are automatically determined by the Moore fillers $\mu^{n+1}_k$. In fact, 
\begin{equation} \label{eq:vs-ngpd-mul}
    m_k := d^{n+1}_k\mu^{n+1}_k, \quad \forall 0\le k \le n+1. 
\end{equation} are the multiplications of a $\VS$ $n$-groupoid. We summarize this in the following theorem.

\begin{theorem}[Finite data for $\VS$ $n$-groupoids]\label{thm:FiniteDataVS}
    The data of a $\VS$ $n$-groupoid is equivalent to the data of an $n$-truncated simplicial vector space.  That is, a $\VS$ $n$-groupoid $\huaV$ can be recovered from its $n$-truncation $\tr^n(\huaV)$. 
    
    The data of a simplicial linear map $f: \huaW \to \huaV$, between a simplicial vector space $\huaW$ and a $\VS$ $n$-groupoid $\huaV$ is equivalent to the data of an $n$-truncated simplicial map $f': \tr^n(\huaW) \to \tr^n(\huaV)$ such that
\begin{equation} \label{eq:mult}
    \forall w \in \huaW_{n+1},  \quad f'(d_k w) = m_k((f'd_i w)_{i \neq k}), \qquad  \text{ for a certain}\; k \in [0, n+1]. 
\end{equation}
 We call such an $f'$ \textbf{multiplicative}. 
\end{theorem}

\begin{proof}For the first part, one direction is obvious. Let us start with an $n$-truncated simplicial vector space $\widebar{\huaV}$. We choose a $k \in [0, n+1]$ and define $m_k:= d^{n+1}_k\mu^{n+1}_k$. Then, 
    \begin{equation}\label{eq:n+1truncationVSnGpd}
        \Lambda^{n+1}_k(\widebar{\huaV}) \overset{d_0, \dots, d_{n+1}}{\rightthreearrows} {\widebar{\huaV}}, 
    \end{equation}
    is an $n+1$-truncated simplicial vector space with face maps $d_{j\neq k}$ the natural projections,  $d_k = m_k$,  and degeneracy maps $s_j$ given by $(s_jx)_i = d_is_j x$ componentwise. Then $\cosk^{n+1}(\Lambda^{n+1}_k(\widebar{\huaV}) \rightthreearrows \widebar{\huaV})$ is a $\VS$ $n$-groupoid.

    For the second part, we have that, by the adjunction $\tr^{n+1} \dashv \cosk^{n+1}$ and the fact that $\huaV$ is a $\VS$ $n$-groupoid,
    \begin{equation*}
        \begin{split}
        \SVect(\huaW, \huaV) 
        &\cong \SVect (\huaW, \cosk^{n+1}\tr^{n+1}\huaV) 
        \cong \SVect_{\le n+1}(\tr^{n+1}\huaW, \tr^{n+1}\huaV)\\
        &\cong \SVect_{\le n+1}(\tr^{n+1}\huaW, \Lambda^{n+1}_k(\huaV) \rightthreearrows \tr^n{\huaV}),
        \end{split}
    \end{equation*}
    where the last isomorphism is given by writing $\tr^{n+1}\huaV$ as in \eqref{eq:n+1truncationVSnGpd} by applying $p^{n+1}_k$ at level $n+1$. With this, all that is left to show is that the data of an $n$-truncated multiplicative simplicial map $f': \tr^n(\huaW) \to \tr^n(\huaV)$ is the same as the data of an $(n+1)$-truncated simplicial map $f: \tr^{n+1}\huaW \to (\Lambda^{n+1}_k(\huaV) \rightthreearrows \tr^n{\huaV})$. 
    
    Given $f$, we define $f'$ by $f': = \tr^nf$. Given $f'$, we define $f$ by $f'$ on level $0, \dots, n$ and $f_{n+1}$ by
\begin{equation}\label{eq:f-from-f'}
    f_{n+1}(w) := (f'(d_0 w), \dots, \widehat{f'(d_k w)}, \dots, f'(d_n w)) = (f'd_i w)_{i \neq k}, \quad \forall w \in \huaW_{n+1}.
\end{equation}
Then the commutativity of $f$ with the simplicial maps between levels $n$ and $n+1$ is equivalent to the multiplicativity of $f'$. We show this in one direction, as the other direction is similar. Given $f$ constructed from $f'$ by \eqref{eq:f-from-f'}, we automatically have $d_j f_{n+1} = f'd_j$ for $j \neq k$. When $j=k$,  $d_k f_{n+1} = m_k((f'd_i)_{i \neq k} )$, which is the multiplicativity condition of $f'$. The commutativity of $f$ with the $s_j$'s follows automatically by the construction of the latter maps in \eqref{eq:n+1truncationVSnGpd}. 
\end{proof}

\begin{remark}
In the above proof, since $k$ is arbitrarily chosen, we see that a multiplicative map $f'$ even satisfies
\begin{equation}\label{eq:mult2}
    \forall w \in \huaW_{n+1}, \quad f'(d_k w) = m_k((f'd_i w)_{i \neq k}),  \qquad \forall k \in [0, n+1],
\end{equation}
which is apparently stronger than \eqref{eq:mult}, but actually equivalent to it. 
\end{remark}

\begin{remark}\label{ref:VSnGpd-is-n-skeletal}
    Let $\huaV$ be a simplicial vector space. By Lemma \ref{lem:MooreFillers-degeneracies}, since $D_m\huaV \cong \Horn{m}{k}(\huaV)$ for all $0\le k \le m$ and $m\ge 0$, the $n$-skeleton $\Sk^n(\huaV)$ is always canonically a $\VS$ $n$-groupoid for any $n\ge 0$. Moreover,  if $\huaV$ is a $\VS$ $n$-groupoid,  then $\huaV  \cong \Sk^n(\huaV) \cong  \Cosk^{n+1}(\Sk^{n} (\huaV)) $. 
\end{remark}

\begin{remark}\label{rem:category-of-n-gpds}
    Let us denote $\SVect_{\le n}^{k-mult}$ the category of $n$-truncated simplicial vector spaces with $n$-truncated simplicial maps $f': \widebar{\huaW} \to \widebar{\huaV}$ that are multiplicative,  which in this case means that $f' \circ m_k^\huaW = m_k^\huaV \circ f'$ for a choice of $0 \le k \le n+1$. As a consequence of the theorem above, this category is equivalent to the category of $\VS$ $n$-groupoids and consequently a model for $(n+1)$-vector spaces generalizing the 2-vector spaces in \cite{BaezCrans2004}. 
    
    However, $\SVect_{\le n}^{k-mult}$ is not closed under the monoidal product in $\SVect$, as we show in Theorem \ref{thm:order-of-tensors}, so it is not a monoidal subcategory. Thus we will not restrict our discussion to this category. 
\end{remark}

\begin{example}\label{ex:MooreFillers2Gpd}
    In analogy with Example \ref{ex:MooreFillers1Gpd}, we can describe the canonical multiplications for a $\VS$ $2$-groupoid $\huaV_2 \rightthreearrows \huaV_1 \rightrightarrows \huaV_0$. For our convenience, we introduce a new notation for the ternary multiplications of a $2$-groupoid object $\huaX$. That is, we denote each $m_i$ of three triangles fitting in a $(3,i)$-horn by writing the missing face  with an empty box $\square$. More precisely,  given $(X,Y,Z) \in \Lambda^3_0(\huaX) \cong \huaX_3$, let $W= d_0(X,Y,Z) $, then the four ternary multiplications are
    \begin{equation*}
        \begin{array}{ccc}
            W = m_0(X,Y,Z) =: \square X Y Z,
            &\qquad
            &Y = m_2(W,X,Z) =: W X \square Z,\\
            X = m_1(W,Y,Z) =: W \square Y Z,
            &\qquad
            &Z = m_3(W,X,Y) =: W X Y \square.
        \end{array}
    \end{equation*} Notice that we automatically have $(W, X, Y) \in \Lambda^3_3(\huaX)$, $(W, Y, Z) \in \Lambda^3_1(\huaX)$, and $(W, X, Z) \in \Lambda^3_2(\huaX)$. 
    In the case of a $\VS$ 2-groupoid $\huaV$, the Moore fillers give
	\begin{equation}\label{eq:TriMultOverPointAll}
    \begin{split}
    	W = m_0(X,Y,Z) = \square X Y Z &= X - Y + Z + s_0d_0 Y - s_0d_0 Z + s_1d_0 Z - s_1d_1 Z,\\
    	X = m_1(W,Y,Z) = W \square Y Z &= W + Y - Z - s_0d_1 W + s_0d_0 Z - s_1d_0 Z + s_1d_1 Z,\\
    	Y = m_2(W,X,Z) = W X \square Z &= - W + X + Z + s_0d_1 W - s_0d_0 Z + s_1d_0 Z - s_1d_1 Z,\\
    	Z = m_3(W,X,Y) = W X Y \square &= W - X + Y - s_0d_1 W + s_0d_2 W - s_1d_2 W + s_1d_2 X,
    \end{split}
	\end{equation}
    for $(W,X,Y,Z) \in \huaV_3$.
\end{example}

\subsection{The simplicial category of simplicial vector spaces and its homotopy category}\label{sec:Svect-int-hom}

In this subsection, we recall and sort out some neccessary preliminary knowledge on the simplicial category $\SVect$. See e.g. \cite{GoerssJardine2009,Riehl2014,kerodon}. 
This is a closed monoidal category, which by definition is a category equipped with a monoidal product and an internal hom which is right adjoint to it.

The monoidal product $\otimes$ in $\SVect$ is given levelwise by the tensor product of vector spaces. 
For any two simplicial vector spaces $\huaV$ and $\huaW$, the face and degeneracy maps of $\huaV \otimes \huaW$ are $d_i^\huaV \otimes d_i^\huaW$ and $s_j^\huaV \otimes s_j^\huaW$, respectively. 
This product is symmetric with respect to the symmetry isomorphism $v \otimes w \mapsto w \otimes v$.
In addition, there is a way to ``tensor'' a simplicial vector space and a simplicial set, known as \textbf{copowering}, which is defined at each level $m$ and for any $\huaV \in \Simp\Vect$ and $\huaK \in \Simp\Set$, by
\begin{equation}\label{eq:copowering}
(\huaV \otimes \huaK)_m
:= \huaV_m \otimes \KK [\huaK_m]
:= \bigoplus_{u \in \huaK_m} \huaV_m^u,
\end{equation}
The simplicial structure of $\huaV \otimes \huaK$ is similar to that of the tensor product. More explicitly, for any $v^r \in \huaV_m^r \subseteq (\huaV \otimes \huaK)_m$,  the face and degeneracy maps of $\huaV \otimes \huaK$ are given in components by 
\begin{equation*}
    d_i(v^r) 
    = (d_iv)^{d_i r}
    \in \huaV_{m-1}^{d_i r} 
    \subseteq (\huaV \otimes \huaK)_{m-1},
    \qquad s_i(v^r) 
    = (s_iv)^{s_i r}
    \in \huaV_{m+1}^{s_i r} 
    \subseteq (\huaV \otimes \huaK)_{m+1}. 
\end{equation*}

The \textbf{internal hom} is given by the mapping space $\IHom(\huaV, \huaW)$, which is defined at each level by
\begin{equation}\label{eq:HomSpaceDef}
    \IHom(\huaV, \huaW)_m
    := \SVect(\huaV \otimes \Delta[m], \huaW),
\end{equation}
where the right-hand side is the usual hom vector space in $\SVect$.
The face and degeneracy maps of the mapping space are induced on each level by the coface and codegeneracy simplicial morphisms between $\Delta[m]$ and $\Delta[m-1]$ (or $\Delta[m+1]$). 
More explicitly, given a simplicial map $f \in \IHom(\huaV, \huaW)_m$, its $l$-th level is a map $f_l: (\huaV \otimes \Delta[m])_l \to \huaW_l$. Then, face and degeneracy maps are defined in components by
\begin{equation}\label{eq:HomSpaceMaps}
    \begin{split}
        (d_if)_l (x^u) &:= f_l (x^{\delta^i(u)}), \quad \forall x \in \huaV_{l}, u \in \Delta[m-1]_{l},\\
        (s_if)_l (x^u) &:= f_l (x^{\sigma^i(u)}), \quad \forall x \in \huaV_{l}, u \in \Delta[m+1]_{l},
    \end{split}
\end{equation}
where $\delta^i: \Delta[m-1]_{l} \to  \Delta[m]_{l}$ are the coface maps, and $\sigma^i: \Delta[m+1]_{l} \to \Delta[m]_{l}$ are the codegeneracy maps\footnote{If we write simplices in $\Delta[m]_l$ as non-decreasing sequences of $l$ numbers in $\{0,1,2,\dots,n\}$, then the coface map $\delta^i$ is the postcomposition with the unique injection of ordinals that skips $i$, and the codegeneracy map $\sigma^i$ is by postcomposition with the unique surjection of ordinals that repeats $i$.}. It is easy to see that they are linear and satisfy the simplicial identities, so that $\IHom(\huaV, \huaW)$ is indeed a simplicial vector space. The fact that $\IHom$ is right adjoint to the monoidal product $\otimes$ follows from the same proof as in \cite[Prop. I.5.1]{GoerssJardine2009}. By the discussion in \cite[Remark 3.3.9, \S 3.7]{Riehl2014}, the tensor-hom adjunction upgrades to an $\SVect$-adjunction. We prove this explicitly in Proposition \ref{prop:tensor-hom-svect} later as it is an important result for the rest of the article. 

We now describe more precisely the simplicial category structure of $\SVect$ in order to introduce its homotopy category and show that the tensor-hom adjunction descends to it. 
The internal hom $\IHom$ is a functor $\SVect \times \SVect \to \SVect$, whose level 0 coincides with the hom space in $\SVect$. Going higher, the level 1 of $\IHom(\huaV, \huaW)$ is the space of simplicial homotopies between simplicial linear maps in $\IHom(\huaV, \huaW)_0$: a \textbf{simplicial homotopy} between two simplicial linear maps $f,g:\huaV \to \huaW$ is a simplicial linear map $h: \huaV \otimes \Delta[1] \to \huaW$, such that 
\begin{equation*}
    \begin{split}
        &d_0 h = h \circ (id \otimes \delta_0) = f: \huaV \otimes \Delta[0] \cong \huaV \to \huaW, \text{ and }\\
        &d_1 h = h \circ (id \otimes \delta_1) = g: \huaV \otimes \Delta[0] \cong \huaV \to \huaW.
    \end{split}
\end{equation*}
Analogously, the higher levels of the internal hom encode increasingly higher homotopies.
These higher homotopies can be composed across different internal homs according to the composition defined at each level $m$ by
\begin{equation*}
    \begin{split}
        &\cdot \circ \cdot : \IHom(\huaV, \huaW)_m \otimes \IHom(\huaU, \huaV)_m \longrightarrow \IHom(\huaU, \huaW)_m\\
        &k \circ h: \huaU \otimes \Delta[m] \xrightarrow{id \times diag} \huaU \otimes (\Delta[m] \times \Delta[m]) \cong (\huaU \otimes \Delta[m]) \otimes \Delta[m] \xrightarrow{h \otimes id} \huaV \otimes \Delta[m] \xrightarrow{k} \huaW,
    \end{split}
\end{equation*}
where $diag$ is the diagonal inclusion. 
Note that for $m=1$, if $h$ is a homotopy between $f$ and $f'$ and $k$ is a homotopy between $g$ and $g'$, then $h \circ k$ is a homotopy between $g\circ f$ and $g'\circ f'$.
It is straightforward to check that with this composition and the obvious identity element, $\SVect$ is an $\SVect$-category, i.e. a category enriched in $\SVect$,  in the sense of \cite[Def. 3.3.1]{Riehl2014}. In particular it is a simplicial category, i.e. a category enriched in $\SSet$, in the sense of \cite[\S 3.6]{Riehl2014}, \cite[\href{https://kerodon.net/tag/00JQ}{Tag 00JQ}]{kerodon}.

We can now define the \textbf{homotopy category} $\hSVect$ as the category with objects simplicial vector spaces and with morphisms simplicial homotopy classes of simplicial linear maps. In other words, the hom spaces in $\hSVect$ are given by $\pi_0\IHom$. For a general construction of this starting from a simplicial category see for example \cite[\href{https://kerodon.net/tag/00LW}{Tag 00LW}]{kerodon}. In this case, the construction is simplified by the fact that all simplicial vector spaces  are Kan by Theorem \ref{thm:MooreHornFillers}, so homotopy of maps is an equivalence relation. We denote it by $\sim$. See for example \cite[Lemma I.6.1]{GoerssJardine2009}, \cite[\href{https://kerodon.net/tag/00HC}{Tag 00HC}]{kerodon} and \cite[\href{https://kerodon.net/tag/00M0}{Tag 00M0}]{kerodon} for more details. An isomorphism in the homotopy category is precisely a \textbf{homotopy equivalence}. That is, a pair of maps $f: \huaV \to \huaW$ and $g: \huaW \to \huaV$ such that $fg \sim id_{\huaW}$ and $gf \sim id_{\huaV}$. In this case, we write $\huaV \simeq \huaW$. The tensor product also descends in a straightforward way to the homotopy category.

\begin{proposition}[Enriched tensor-hom adjunction for $\SVect$]\label{prop:tensor-hom-svect}
    For any $\huaU, \huaV, \huaW$ in $\SVect$, there are natural isomorphisms of simplicial vector spaces
    \begin{equation}\label{eq:tensor-hom-svect}
        \IHom(\huaU \otimes \huaV, \huaW) \cong \IHom(\huaU, \IHom(\huaV, \huaW)) \cong
        \IHom(\huaV, \IHom(\huaU, \huaW)).
    \end{equation}
\end{proposition}
\begin{proof}
    We construct a pair of $\SVect$-natural morphisms which are inverse to each other:
    \begin{equation}\label{eq:tensor-hom-rho-tau}
    \IHom(\huaU \otimes \huaV, \huaW) \overset{\xrightarrow{\rho_{\huaU, \huaV, \huaW}}}{\xleftarrow[\tau_{\huaU, \huaV, \huaW}]{}}  \IHom(\huaV, \IHom(\huaU, \huaW)). 
    \end{equation}
    The construction to show $\IHom(\huaU \otimes \huaV, \huaW) \cong \IHom(\huaU, \IHom(\huaV, \huaW))$ is similar. 

    We write  $f \in \IHom(\huaU \otimes \huaV, \huaW)_m$ in components $f_i: \huaU_i \otimes \huaV_i \otimes \Delta[m]_i \to \huaW_i$. We view an element $(v,r) \in \huaV_l \otimes \Delta[m]_l$ as a simplicial linear map $(v,r) : \Delta[l] \to \huaV \otimes \Delta[m]$. Then we define $\rho_{\huaU, \huaV, \huaW}(f) \in \IHom(\huaV, \IHom(\huaU, \huaW))_m$ in components $\rho_{\huaU, \huaV, \huaW}(f)_l: \huaV_l \otimes \Delta[m]_l \longrightarrow \SVect(\huaU \otimes \Delta[l], \huaW)$ for any $(v,r)$ as the composition 
    \begin{equation*}
    \rho_{\huaU, \huaV, \huaW}(f)_l(v,r):\; \huaU \otimes \Delta[l] \xrightarrow{id_\huaU \otimes(v, r)} \huaU \otimes \huaV \otimes \Delta[m] \xrightarrow{f} \huaW. 
    \end{equation*}
    For a more explicit expression of $\rho_{\huaU, \huaV, \huaW}(f)_l$, we observe that any $i$-simplex $t \in \Delta[l]_i$ can be expressed in terms of the unique non-degenerate $l$-simplex $E_l = 012\dots l \in \Delta[l]_l$, as $t = s_I d_J E_l$, for some multi-indices $I,J$ such that $|I| - |J| = i - l$. 
    Then
    \begin{equation}\label{eq:tensor-hom-svect-def-rho}
        \begin{split}
            \rho_{\huaU, \huaV, \huaW}(f)_l: &\huaV_l \otimes \Delta[m]_l \longrightarrow \SVect(\huaU \otimes \Delta[l], \huaW)\\
            &\rho_{\huaU, \huaV, \huaW}(f)_l(v,r)(u,t) = f_i(u, s_I d_Jv, s_I d_J r) \in \huaW_i,
        \end{split}
    \end{equation}
    for any $(v,r) \in \huaV_l \otimes \Delta[m]_l$, and any $(u,t) \in \huaU_i \otimes \Delta[l]_i$ with $t=s_Id_JE_l$. Notice that the choice of the multiindices $I,J$ may not be unique. But different choices give rise to the same result.

    For the other direction,  we start with maps with components $g_l: \huaV_l \otimes \Delta[m]_l \to \SVect(\huaU \otimes \Delta[l], \huaW)$, and we define 
    \begin{equation}\label{eq:tensor-hom-svect-def-tau}
            \tau_{\huaU, \huaV, \huaW}(g)_l: \huaU_l \otimes \huaV_l \otimes \Delta[m]_l \longrightarrow \huaW_l, \quad \text{by} \quad \tau_{\huaU, \huaV, \huaW}(g)_l(u, v, r) = g_l(v,r)(u, E_l) \in \huaW_l,
    \end{equation}
    where $E_l \in \Delta[l]_l$ is the unique non-degenerate simplex therein, as before.
    Equivalently, we can write $\tau(g) = ev \circ (id_{\huaU} \otimes g)$, where $ev: \huaU \otimes \Delta[l] \otimes \SVect(\huaU \otimes \Delta[l], \huaW) \to \huaW$ is the canonical evaluation map, which is simplicial \cite[\S I.5]{GoerssJardine2009} and linear. 
\end{proof}

From this, we recover the usual unenriched hom-tensor adjunction by taking the simplicial level 0 in \eqref{eq:tensor-hom-svect}. Analogously, the adjunction descends to the homotopy category by applying the functor $\pi_0$ to \eqref{eq:tensor-hom-svect}.

\begin{corollary}[Tensor-hom adjunction for $\hSVect$]\label{cor:tensor-hom-hsvect}
    For any $\huaU, \huaV, \huaW$ in $\SVect$, there are natural isomorphisms
    \begin{equation}\label{eq:tensor-hom-hsvect}
        \hSVect(\huaU \otimes \huaV, \huaW) \cong \hSVect(\huaU, \IHom(\huaV, \huaW)) \cong
        \hSVect(\huaV, \IHom(\huaU, \huaW)).
    \end{equation}
\end{corollary}

\subsection{The DG-category of non-negative chain complexes and its homotopy category}\label{sec:monoidal-struct-chains}

The category of chain complexes of vector spaces $\Ch$ is also a closed symmetric monoidal category which is enriched over itself. Its symmetric monoidal structure is given by the tensor product defined degreewise as
\begin{equation*}
    (A \otimes B)_m = \bigoplus_{p+q=m} A_p \otimes B_q,
\end{equation*}
for chain complexes $(A, \partial)$, $(B, \partial)$, with the differential defined on homogeneous elements by $\partial(a \otimes b) = \partial a \otimes b + (-1)^{|a|} a \otimes \partial b$, where $|a|$ denotes the degree of $a$. 
It is important to note that with these definitions, the symmetry isomorphism $A \otimes B \to B \otimes A$ has a sign: $a \otimes b \mapsto (-1)^{|a||b|}b \otimes a$.

The \textbf{internal hom} is given by the mapping complex $\IHom(A,B)$, which is the chain complex with degree $m$ the space of degree $m$ maps of the underlying graded vector spaces 
\begin{equation} \label{eq:chain-iHom}
    \IHom(A,B)_m = \{ f: A_\bullet \to B_{\bullet + m} \},
\end{equation}
and differential $\partial f = \partial \circ f - (-1)^{|f|} f \circ \partial$. 
Clearly this is a chain complex of vector spaces, as every level is a vector space and the differential is linear. With this, and the obvious composition rule, $\Ch$ becomes a category enriched over itself, also commonly known as a DG-category \cite[\href{https://kerodon.net/tag/00ND}{Tag 00ND}]{kerodon}. 

Note that the internal hom is generally non-zero in both positive and negative degrees and that the usual space of chain maps is given by the space of 0-cycles:
\begin{equation*}
    \Ch(A, B) = \{f: A \to B | \partial f = f \partial \} = \ker \partial_0.
\end{equation*}
Furthermore, two chain maps $f, g:A \to B$ are homotopic if and only if they differ by a 1-boundary $\partial_1 h \in \Img(\partial_1)$ of the internal hom, i.e. $f - g = \partial h + h \partial$. This has the immediate consequence that homotopy of chain maps is an equivalence relation which is compatible with composition of chain maps. However, we need to be careful that chain maps and chain homotopies do {\it not} make chain complexes into a bicategory because interchange laws do not hold on the nose \cite{mo:interchange1,mo:interchange2}.\footnote{One can obtain a 2-category only by considering \textit{homotopy classes} of homotopies as 2-morphisms as in \cite[V.8.4]{GabrielZisman1967}, which amounts to truncating the $\infty$-category at 2-morphisms, while the homotopy category truncates it at 1-morphisms. This is enough for us.}

A more convenient solution for us is to consider the \textbf{homotopy category} of chain complexes $\hCh$ (see e.g. \cite[\href{https://kerodon.net/tag/00NM}{Tag 00NM}]{kerodon}) with objects chain complexes of vector spaces and morphisms the chain maps up to homotopy, i.e. $\hCh(A,B) = H_0(\IHom(A,B))$, for any two chain complexes $A$ and $B$. 
In particular, an isomorphism in the homotopy category $\hCh$ is a \textbf{chain homotopy equivalence}. That is, a pair of maps $f:A \to B$ and $g:B \to A$ such that $fg \sim id_B$ and $gf \sim id_A$, where $\sim$ denotes the homotopy relation. In this case, we write $A \simeq B$. The monoidal structure also descends to the homotopy category, because given $f,g: A \to B$ homotopic via $h$ and $f',g':A' \to B'$ homotopic via $k$, there is a homotopy between $f \otimes f'$ and $g \otimes g'$ given by $h \otimes f' + g \otimes k$.\footnote{Observe here that by the Koszul sign rule $(g \otimes k)(a \otimes a')=(-1)^{|k||a|}g(a) \otimes k(a')=(-1)^{|a|}g(a) \otimes k(a')$.}

Consider now the full subcategory of non-negative chain complexes $\ChnV$. This is a DG-subcategory, but we can also see it as a category enriched over itself by truncating the internal hom at 0 and replacing its 0-chains by its 0-cycles. That is, we define, for any non-negative $A$ and $B$,
\begin{equation}\label{eq:non-neg-chain-maps}
    \IHom_{\ge 0}(A, B)_i := tr_{\ge 0}\IHom(A, B)_i = \begin{cases}
        \IHom(A,B)_i &\text{ for } i > 0\\
        \ker\partial_0 &\text{ for } i = 0\\
        0 &\text{ for } i < 0.
    \end{cases}
\end{equation}
With this, $\ChnV$ is a category enriched over itself. Since $H_0= \ker\partial_0 / \Img \partial_1 $, the homotopy category $\hChP$ of $\ChnV$  is the full subcategory of $\hCh$ generated by the non-negative complexes. More precisely, for $A,B \in \ChnV$, we have 
\begin{equation*}
    \hChP(A,B) = H_0(\IHom_{\ge 0}(A,B)) = H_0(\IHom(A,B)) = \hCh(A,B).
\end{equation*}

We now recall that the tensor-hom adjunction in $\Ch$ upgrades to an adjunction of $\Ch$-categories. It is in fact more practical to show the enriched adjunction first. The usual adjunction then follows by taking the 0-cycles of the internal homs. 
This is for example the content of \cite[Exercise 10.8]{Rotman2009} and the incorrectly stated \cite[Exercise 2.7.3]{Weibel1994}, but we did not otherwise find a complete statement in the literature. The proof is quite straightforward, the main difficulty being choosing the right sign conventions so that the signs of the differentials match. 
A proof using triple complexes can be found at \cite{SE:4803304}.
We write here the natural morphisms composing the adjunction, which we will use later in the article and leave it to the reader to check that they form a $\Ch$-natural isomorphism.

\begin{proposition}\label{prop:tensor-hom-chains}
    For any $A,B,C \in \Ch$, there are natural isomorphism of chain complexes
    \begin{equation}\label{eq:tensor-hom-chains}
        \IHom(A \otimes B, C) \cong \IHom(A, \IHom(B, C)) \cong
        \IHom(B, \IHom(A, C)).
    \end{equation}
\end{proposition}

\begin{proof}
The natural isomorphism $\rho_{A,B,C}: \IHom(A \otimes B, C) \to \IHom(A, \IHom(B, C))$ is given in each degree by 
\begin{equation}\label{eq:tensor-hom-chvect-def-rho}
    \begin{split}
        \rho_{A, B, C}(f): &A_{\bullet} \to \IHom_{\bullet+|f|}(B,C)\\
        &\rho_{A, B, C}(f)(a)(b) = f(a \otimes b) \in C_{|a|+|b|+|f|},
    \end{split}
\end{equation}
for a homogeneous element $f \in \IHom_{|f|}(A \otimes B, C)$. 
Its inverse $\tau_{A,B,C}: \IHom(A, \IHom(B, C)) \to \IHom(A \otimes B, C)$ is given by
\begin{equation}\label{eq:tensor-hom-chvect-def-tau}
    \begin{split}
        \tau_{A, B, C}(g): &(A \otimes B)_\bullet \to C_{\bullet+|g|}\\
        &\tau_{A, B, C}(g)(a \otimes b) = g(a)(b) \in C_{|a|+|b|+|g|},
    \end{split}
\end{equation}
for a homogeneous element $g \in \IHom_{|g|}(A, \IHom(B, C))$.
It is worth noting that the other natural isomorphism $\rho^L_{A,B,C}: \IHom(A \otimes B, C) \to \IHom(B, \IHom(A, C))$ obtained by precomposing $\rho$ with the symmetry isomorphism has a sign. In fact, $\rho^L$ is given by
\begin{equation}\label{eq:tensor-hom-chvect-def-rho-left}
    \begin{split}
        \rho^L_{A, B, C}(f): &B_{\bullet} \to \IHom_{\bullet+|f|}(A,C)\\
        &\rho^L_{A, B, C}(f)(b)(a) = (-1)^{|a||b|} f(a \otimes b) \in C_{|a|+|b|+|f|},
    \end{split}
\end{equation}
for any homogeneous $f \in \IHom_{|f|}(A \otimes B, C)$. 
\end{proof}

This restricts to an adjunction of $\ChnV$-categories on $\ChnV$, and it induces a tensor-hom adjunction on the homotopy categories $\hCh$ and $\hChP$. 

\begin{corollary}\label{cor:tensor-hom-positive-chains}
    For any $A,B,C \in \ChnV$, there are natural isomorphism of chain complexes
    \begin{equation}\label{eq:tensor-hom-non-neg-chains}
        \IHom_{\ge 0}(A \otimes B, C) \cong \IHom_{\ge 0}(A, \IHom_{\ge 0}(B, C)) \cong
        \IHom_{\ge 0}(B, \IHom_{\ge 0}(A, C)).
    \end{equation}
\end{corollary}

\begin{proof}
    It follows from the previous proposition that if the complexes in \eqref{eq:tensor-hom-chains} are isomorphic, their $0$-truncations obtained by \eqref{eq:non-neg-chain-maps} are isomorphic. That is
    \begin{equation*}
        \IHom_{\ge 0}(A \otimes B, C) \cong \IHom_{\ge 0}(A, \IHom(B, C)) \cong
        \IHom_{\ge 0}(B, \IHom(A, C)).
    \end{equation*}
    To obtain \eqref{eq:tensor-hom-non-neg-chains} it is enough to observe that for any non-negative chain complexes $A, B, C$,
    \begin{equation*}
        \IHom_{\ge 0}(A, \IHom(B, C)) = \IHom_{\ge 0}(A, \IHom_{\ge 0}(B, C)). 
    \end{equation*}
This is    because   a chain map from a non-negative chain complex does not detect anything in negative degree, and it sends 0-chains to 0-cocycles, as any 0-chain in a non-negative is a 0-cocycle.
\end{proof}

\begin{corollary}\label{cor:tensor-hom-hchp}
    For any $A,B,C \in \ChnV$, there are natural isomorphisms
    \begin{equation*}
        \hCh(A \otimes B, C) \cong \hCh(A, \IHom_{\ge 0}(B, C)) \cong
        \hCh(B, \IHom_{\ge 0}(A, C)).
    \end{equation*}
\end{corollary}

\section{Eilenberg-Zilber for internal homs, orders, amplitudes and equivalences}

We first recall the Dold-Kan correspondence, which holds in general for simplicial objects in abelian categories. For more details we refer to \cite[\S 8.4]{Weibel1994}, \cite[\S 22]{May1967}, \cite[\S III.2]{GoerssJardine2009}, \cite[\href{https://kerodon.net/tag/00QQ}{Tag 00QQ}]{kerodon}, for example.

\begin{definition}\label{def:Moore-normalized}
    Let $\huaV$ be a simplicial vector space. The \textbf{Moore complex} of $\huaV$ is the non-negative chain complex $C(\huaV)$ with $C(\huaV)_n =\huaV_n$ for any $n\ge 0$. The differential is given at each level $n$ by the boundary map
    \begin{equation*}
        \partial_n = \sum_{i=0}^n (-1)^{i} d_i^n.
    \end{equation*}

    The \textbf{normalized complex} of $\huaV$ is the subcomplex $N(\huaV)$ of the Moore complex with
    \begin{equation}\label{eq:normalized}
        N(\huaV)_n = \ker p^n_n = \bigcap_{i=0}^{n-1} \ker d_i^n.
    \end{equation}
    The Moore complex differential restricted to this is $\partial_n = (-1)^n d_n^n$ at each level $n$.
\end{definition}

\begin{proposition}[Dold-Kan correspondence]\label{prop:DK-correspondence}
    The normalized complex functor $N: \Simp\Vect \to \ChnV$ admits an inverse functor $DK: \ChnV \to \Simp\Vect$, and the two form an equivalence of categories between the category of simplicial vector spaces and that of non-negative chain complexes. Under this equivalence, simplicial homotopies $h: \huaV \otimes \Delta[1] \to \huaW$ between two simplicial maps $f,g:\huaV \to \huaW$, correspond to chain homotopies $N(h): N(\huaV) \to N(\huaW)[-1]$ between $N(f)$ and $N(g)$. 
\end{proposition}

Since the Dold-Kan correspondence sends homotopic maps in $\SVect$ to homotopic maps in $\ChnV$, it descends to an equivalence of the homotopy categories as well.

\begin{corollary}\label{cor:DK-htpy-cats}
    The Dold-Kan correspondence induces an equivalence of categories between the homotopy categories $\hSVect$ and $\hChP$.
\end{corollary}

As a first consequence, the normalized complex of a simplicial vector space contains information about its ``order'' as a higher groupoid. We now give a precise definition of this in the following discussion. 

\begin{definition}
    Let $(E,\partial)$ be a chain complex. We say $(E,\partial)$ has \textbf{amplitude} $(0,n)$ if it is concentrated in degrees 0 to $n$, that is $E_i = 0$ for $i< 0$ or $i>n$.
\end{definition}

\begin{proposition}\label{prop:VS-order-amplitude}
    A simplicial vector space $\huaV$ is a $\VS$ $n$-groupoid if and only if $N(\huaV)$ has amplitude $(0,n)$. 
\end{proposition}

\begin{proof}
    By definition, all normalized complexes are already non-negative. 
    If $\huaV$ is a $\VS$ $n$-groupoid, then clearly $N(\huaV)_i = \ker p_i^i = 0$ for any $i>n$. Hence $N(\huaV)$ has amplitude $(0,n)$.
    
    For the converse assume that $N(\huaV)$ has amplitude $(0,n)$. Then $\ker p^m_m = N(\huaV)_m = 0$ for $m\ge n+1$. Thus $\huaV_m \cong \Lambda^m_m(\huaV)$ for $m\ge n+1$. 
    By Lemma \ref{lem:MooreFillers-degeneracies}, $ \Lambda^m_m(\huaV) \cong D_m(\huaV) \cong \Lambda^m_k(\huaV)$, for all $0 \le k \le m$. Thus $\Kan!(m,k)$ holds for all $0\le k \le m \ge n+1$.
    Hence $\huaV$ is a $\VS$ $n$-groupoid.
\end{proof}

\begin{definition}\label{def:order-type}
    We define the \textbf{order} of a simplicial vector space as the maximal degree for which its normalized chain complex is non-zero. In other words, $\huaV$ has order $n$ if $\ker p^m_m = 0$ for all $m \ge n$ and $\ker p^n_n \neq 0$. 
    We define the \textbf{homotopy type} of a simplicial vector space as the maximal degree for which the homology of the normalized complex is non-zero. That is, we say $\huaV$ is of type $n$, or an \textbf{$n$-type},  if $H_m(N(\huaV)) = 0$ for all $m \ge n$ and $H_n(N(\huaV)) \neq 0$. If $H_i(N(\huaV)) = 0$ for all $i$, then we call $\huaV$  \textbf{acyclic}.
\end{definition}

\begin{remark}\label{rem:H-isom-Pi}
    It follows from Definition \ref{def:order-type} and Proposition \ref{prop:VS-order-amplitude} that a $\VS$ $n$-groupoid can only have homotopy type lower or equal to $n$. 
    Recall (e.g. from \cite[Thm. 22.1]{May1967}, \cite[III.2.7]{GoerssJardine2009}) that for a $\VS$ $n$-groupoid $\huaV$,  $H_i(N(\huaV)) \cong \pi_i(\huaV, 0)$.  This is where our terminology of homotopy type comes from.  
\end{remark}

\begin{example}
    A $\VS$ 1-groupoid of order 1 can be seen as a $\VS$ $n$-groupoid for any $n \ge 1$ but it is not an $n$-groupoid of order $n$ for any $n > 1$. 
    $\VS$ 1-groupoids can be of homotopy type 1, 0 or be acyclic. $\VS$ 2-groupoids can in addition be of type 2.
\end{example}

To be able to discuss in what sense the Dold-Kan correspondence is monoidal up to homotopy, we need to describe in what sense simplicial vector spaces can be equivalent up to homotopy. 
In many cases, the notion of homotopy equivalence is too strong, and one introduces weak equivalences and model structures instead (see e.g. \cite{GoerssJardine2009}). 
For $\SVect$, the situation is especially simple, as weak equivalences are the same as homotopy equivalences. We summarize this in the following proposition. 

\begin{proposition}\label{prop:equivalences-in-svect}
    In $\SVect$, the category of simplicial vector spaces, the following are equivalent:
    \begin{enumerate}
        \item $f:\huaV \to \huaW$ is a \textbf{weak equivalence} in the sense that for each point $p \in \huaV_0$, the induced map $\pi_i(\huaV,p) \to \pi_i(\huaW, f(p))$ on the homotopy groups is an isomorphism for $i \ge 0$. 
        \item $f:\huaV \to \huaW$ is a \textbf{homotopy equivalence} in the sense that it admits a homotopy inverse, that is there exists $g:\huaW \to \huaV$ such that $fg \sim id_{\huaW}$ and $gf \sim id_{\huaV}$.
        \item $N(f): N(\huaV) \to N(\huaW)$ is a \textbf{quasi-isomorphism} of chain complexes, in the sense that it induces an isomorphism between the respective homologies. 
        \item $N(f): N(\huaV) \to N(\huaW)$ is a \textbf{chain homotopy equivalence}, in the sense that it admits a chain homotopy inverse, that is there exists a chain map $g:N(\huaW) \to N(\huaV)$ such that $fg \sim id_{N(\huaW)}$ and $gf \sim id_{N(\huaV)}$. 
    \end{enumerate}
\end{proposition}

\begin{proof}
    By the Dold-Kan correspondence being an equivalence between the homotopy categories (Corollary \ref{cor:DK-htpy-cats}), (2) and (4) are equivalent.

    By the discussion in \cite[\S 1.4]{Weibel1994}, (3) and (4) are equivalent. In one direction,  it is clear that chain homotopy equivalences are quasi-isomorphisms. The other direction follows from the observation that a chain complex $(A,\partial)$ splits as $A_n \cong H_n(A) \oplus \Img\partial_{n+1} \oplus \Img\partial_{n}$. This gives a chain homotopy equivalence between $A$ and its homology $H(A)$. 

    Finally (1) and (3) are equivalent because $\pi_i(\huaV, 0) \cong H_i(N(\huaV))$ for any $i\ge 0$. The change of basepoint can be accounted for by using the isomorphism $\pi_i(\huaV, 0) \to \pi_i(\huaV, p)$ induced by multiplication by the unit $1p$ of each basepoint $p\in \huaV_0$.
\end{proof}

\begin{remark}\label{rem:n-type-model}
    It is important to note that all of the above equivalences preserve the homotopy type but not the order. 
\end{remark}

\begin{remark}
    A reader familiar with the theory of higher Lie groupoids might wonder what a Morita equivalence is for $\SVect$. 
    Morita equivalences between $n$-groupoid objects were defined in \cite[\S 2]{Zhu2009} as spans of hypercovers (see also \cite{RogersZhu2020}). 
    In the case of simplicial vector spaces, which are $\infty$-groupoid objects,  a hypercover is a trivial (or acyclic) fibration in the standard model category. By a classical result (e.g. \cite[Lemma 3.2]{Curtis1971}, \cite[Lemma III.2.10-11]{GoerssJardine2009}) these are the levelwise surjective weak equivalences.

    By Brown's factorization Lemma, a weak equivalence between fibrant objects (i.e. $\infty$-groupoids in this case) can be written as a span of trivial fibrations.  Thus a Morita equivalence between simplicial vector spaces is equivalently one of the following
    \begin{enumerate}
        \item a span of levelwise surjective weak equivalences,
        \item a span of weak equivalences, 
        \item a span of  homotopy equivalences.
    \end{enumerate}
\end{remark}

\subsection{Eilenberg-Zilber theorem for internal homs}

The Eilenberg-Zilber theorem describes how the Dold-Kan correspondence acts on the tensor products on each side. This originally appeared in \cite{EilenbergZilber1953}, \cite[Thm. 2.1]{EilenbergMacLane1954}. See also \cite[\S 29]{May1967}, \cite[Tags \href{https://kerodon.net/tag/00RF}{00RF},\href{https://kerodon.net/tag/00S0}{00S0}]{kerodon}. We summarize neccessary definitions and results in this section. Then we give a version of Eilenberg-Zilber theorem for internal homs in Theorem \ref{thm:EZ-thm-hom}.

\begin{definition}\label{def:EilenbergZilberMap}
    Let $\huaV$ and $\huaW$ be two simplicial vector spaces. The (normalized) \textbf{Eilenberg-Zilber map} is the map 
    \begin{equation*}
        EZ : N(\huaV) \otimes N(\huaW) \longrightarrow N(\huaV \otimes \huaW), 
    \end{equation*}
    which is natural both in $\huaV$ and $\huaW$, and defined on elementary tensors $v\otimes w \in N(\huaV)_p \otimes N(\huaW)_q$ as 
    \begin{equation}\label{eq:EZ}
        EZ(v\otimes w)=\sum_{(\mu,\nu)\in\Shuf(p,q)} \sign(\mu,\nu) (s_{\nu_q}\dots s_{\nu_1} v) \otimes (s_{\mu_p}\dots s_{\mu_1} w),
    \end{equation}
    where $\Shuf(p,q)$ is the set of $(p,q)$-shuffles, which we write as 
    \begin{equation*}
        \sigma = (\mu,\nu) = (\mu_1\dots\mu_p,\nu_1\dots\nu_q) = (\sigma(0)\dots\sigma(p-1),\sigma(p)\dots\sigma(p+q -1)),
    \end{equation*}
    with $\sign(\mu,\nu)$ the signature of the corresponding permutation $\sigma$. 
\end{definition}

This definition makes use of what is sometimes known as the shuffle product \cite[Tag \href{https://kerodon.net/tag/00RF}{00RF}]{kerodon} between chains. Additionally, this map can also be defined in the same way at the level of the Moore complexes and referred to as the \textit{unnormalized} Eilenberg-Zilber map. It can then be observed that it preserves normalized chains by the simplicial identities. 

\begin{proposition}[Normalized Eilenberg-Zilber Theorem {\cite[Thm. 2.1a]{EilenbergMacLane1954}}]\label{prop:EZ-thm-tensor}
    Let $\huaV$ and $\huaW$ be two simplicial vector spaces. The Eilenberg-Zilber map $EZ$ admits a left inverse $AW$, called Alexander-Whitney map,
    \begin{equation*}
        N(\huaV) \otimes N(\huaW) \newrightleftarrows{EZ}{AW} N(\huaV \otimes \huaW),
    \end{equation*}
    such that 
    \begin{equation} \label{eq:EZ-AW}
        AW \circ EZ = id_{N(\huaV) \otimes N(\huaW)}, 
        \qquad 
        EZ \circ AW \sim id_{N(\huaV \otimes \huaW)}.
    \end{equation}
    Both of these maps and the homotopy are natural in $\huaV$ and $\huaW$.
    In other words, $AW$ and $EZ$ form a deformation retraction and they provide a natural chain homotopy equivalence between $N(\huaV)\otimes N(\huaW)$ and $N(\huaV \otimes \huaW)$.
\end{proposition}

Since chain homotopy equivalences are precisely the isomorphisms in $\hChP$, we have the following result.

\begin{corollary}\label{cor:EZ-implies-DK-homotopy-monoidal}
    The functor $N: \hSVect \to \hChP$ is monoidal, that is, for any simplicial vector spaces $\huaV$ and $\huaW$, 
    \begin{equation*}
        N(\huaV) \otimes N(\huaW) \simeq N(\huaV \otimes \huaW). 
    \end{equation*}
    Moreover, by the Yoneda embedding, for any simplicial vector spaces $\huaU$, $\huaV$ and $\huaW$,
    \begin{equation*}
        \hChP(N(\huaU) \otimes N(\huaV), N(\huaW)) \cong \hChP(N(\huaU \otimes \huaV) , N(\huaW)).
    \end{equation*}
\end{corollary}

We now prove the main result of this section.

\begin{theorem}[Eilenberg-Zilber theorem for internal homs]\label{thm:EZ-thm-hom}
    There is a natural chain homotopy equivalence
    \begin{equation*}
        N(\IHom(\huaV, \huaW)) \newrightleftarrows{EZ^H}{AW^H} \IHom_{\ge 0}(N(\huaV), N(\huaW)).
    \end{equation*}
    such that
    \begin{equation*}
        EZ^H \circ AW^H = id_{\IHom_{\ge 0}(N(\huaV), N(\huaW))}, 
        \qquad 
        AW^H \circ EZ^H \sim id_{N(\IHom(\huaV, \huaW))}.
    \end{equation*}
\end{theorem}

\begin{proof}
    We define the two maps by using the inverse of the Yoneda embedding. 
    Hence we first define natural transformations
    \begin{equation*}
        \Ch_{\ge 0}(N(\_), N(\IHom(\_,\_))) \newrightleftarrows{EZ^H_*}{AW^H_*}
        \Ch_{\ge 0}(N(\_), N(\IHom_{\ge 0}(N(\_), N(\_)))), 
    \end{equation*}
    between these two functors $(\Simp\Vect^{op})^2 \times \Simp\Vect \to \Set$.

    We obtain these in two steps by composing natural transformations in the 2-category of categories. As the first step we define two natural isomorphisms $T$ and $R$
    \begin{equation}\label{eq:RT}
        \mathsf{Ch}_{\ge 0}(N(\_), N(\IHom(\_,\_)))
        \newrightleftarrows{T}{R}
        \mathsf{Ch}_{\ge 0}(N(\_\otimes\_), N(\_)),
    \end{equation}
    which are inverse to each other, as illustrated by the following commutative diagram:
    \begin{equation*}
    \begin{adjustbox}{width=\textwidth}
    \begin{tikzcd}[ampersand replacement=\&]
	\& {\mathsf{SVec}^{op}\times \mathsf{SVec}} \&\& {\mathsf{Ch}_{\ge 0}^{op}\times \mathsf{Ch_{\ge 0}}} \\
	\&\&\& {\mathsf{SVec}^{op}\times \mathsf{SVec}} \\
	{(\mathsf{SVec}^{op})^{\times 2}\times \mathsf{SVec}} \&\&\&\&\&\& {\mathsf{Set}} \\
	\&\&\& {\mathsf{SVec}^{op}\times \mathsf{SVec}} \\
	\& {\mathsf{SVec}^{op}\times \mathsf{SVec}} \&\& {\mathsf{Ch}_{\ge 0}^{op}\times \mathsf{Ch_{\ge 0}}}
	\arrow[""{name=0, anchor=center, inner sep=0}, "{N\times N}", from=1-2, to=1-4]
	\arrow["{DK \times DK}"{description}, from=1-4, to=2-4]
	\arrow["{\mathsf{Ch_{\ge 0}}}", curve={height=-18pt}, from=1-4, to=3-7]
	\arrow["{\mathsf{SVec}}"{description}, from=2-4, to=3-7]
	\arrow["\tau"{description}, shift right=2, between={0.2}{0.8}, Rightarrow, from=2-4, to=4-4]
	\arrow["{id \times \underline{\mathrm{Hom}}}", from=3-1, to=1-2]
	\arrow[""{name=1, anchor=center, inner sep=0}, "{id \times \underline{\mathrm{Hom}}}"{description, pos=0.8}, from=3-1, to=2-4]
	\arrow[""{name=2, anchor=center, inner sep=0}, "{\otimes \times id}"{description, pos=0.7}, from=3-1, to=4-4]
	\arrow["{\otimes \times id}"', from=3-1, to=5-2]
	\arrow["\rho"{description}, shift right=2, between={0.2}{0.8}, Rightarrow, from=4-4, to=2-4]
	\arrow["{\mathsf{SVec}}"{description}, from=4-4, to=3-7]
	\arrow[""{name=3, anchor=center, inner sep=0}, "{N\times N}"', from=5-2, to=5-4]
	\arrow["{\mathsf{Ch_{\ge 0}}}"', curve={height=18pt}, from=5-4, to=3-7]
	\arrow["{DK \times DK}"{description}, from=5-4, to=4-4]
	\arrow["\cong"{description}, shift right, between={0.2}{0.8}, Rightarrow, 2tail reversed, from=0, to=1]
	\arrow["\cong"{description}, between={0.2}{0.8}, Rightarrow, 2tail reversed, from=2, to=3]
    \end{tikzcd}
    \end{adjustbox}
    \end{equation*}
    Here the composition from the left corner to the right corner along the top path is the functor on the left hand side of \eqref{eq:RT}, while the one along the bottom path is the functor on its right hand side. Then $T$ is the top-to-bottom composition of natural isomorphisms and $R$ is the bottom-to-top composition.
    The natural isomorphism in the central diamond is given by the tensor-hom adjunction of $\SVect$ (Prop. \ref{prop:tensor-hom-svect}), while the ones in the upper and lower left quadrilaterals are given by the natural isomorphisms establishing $N \circ DK \cong id$ and $DK \circ N \cong id$ in the Dold-Kan correspondence. The triangles to the right commute because the $DK$ functor is full and faithful, by virtue of being a categorical equivalence, so naturally 
    $\SVect(DK(\_), DK(\_)) \cong \Ch_{\ge 0}(\_, \_)$.
    With this, $T\circ R= id$. 

    In the second step we define the other two natural transformations $E$ and $A$
    \begin{equation}\label{eq:EA}
        \mathsf{Ch}_{\ge 0}(N(\_\otimes\_), N(\_)) \newrightleftarrows{E}{A} \mathsf{Ch}_{\ge 0}(N(\_), \IHom_{\ge 0}(N(\_), N(\_))),
    \end{equation}as illustrated by the following diagram: 
    \begin{equation*}
    \begin{adjustbox}{width=\textwidth}
    % https://q.uiver.app/#q=WzAsNixbMCwyLCIoXFxtYXRoc2Z7U1ZlY31ee29wfSlee1xcdGltZXMgMn1cXHRpbWVzIFxcbWF0aHNme1NWZWN9Il0sWzAsMCwiXFxtYXRoc2Z7U1ZlY31ee29wfVxcdGltZXMgXFxtYXRoc2Z7U1ZlY30iXSxbMiwwLCJcXG1hdGhzZntDaH1fe1xcZ2UgMH1ee29wfVxcdGltZXMgXFxtYXRoc2Z7Q2hfe1xcZ2UgMH19Il0sWzQsMCwiXFxtYXRoc2Z7U2V0fSJdLFsyLDIsIihcXG1hdGhzZntDaH1fe1xcZ2UgMH1ee29wfSlee1xcdGltZXMgMn1cXHRpbWVzIFxcbWF0aHNme0NoX3tcXGdlIDB9fSJdLFs0LDIsIlxcbWF0aHNme0NofV97XFxnZSAwfV57b3B9XFx0aW1lcyBcXG1hdGhzZntDaF97XFxnZSAwfX0iXSxbMiwzLCJcXG1hdGhzZntDaF97XFxnZSAwfX0iXSxbMSwyLCJOXFx0aW1lcyBOIl0sWzAsMSwiXFxvdGltZXMgXFx0aW1lcyBpZCJdLFswLDQsIk5cXHRpbWVzIE4gXFx0aW1lcyBOIiwyXSxbNCwyLCJcXG90aW1lcyBcXHRpbWVzIGlkIiwxXSxbNCw1LCJpZCBcXHRpbWVzIFxcdW5kZXJsaW5le1xcbWF0aHJte0hvbX19X3tcXGdlIDB9IiwyXSxbNSwzLCJcXG1hdGhzZntDaF97XFxnZSAwfX0iLDJdLFsyLDUsIlxccmhvIiwxLHsib2Zmc2V0IjoyLCJzaG9ydGVuIjp7InNvdXJjZSI6MjAsInRhcmdldCI6MjB9LCJsZXZlbCI6Mn1dLFsxLDQsIkVaXioiLDEseyJvZmZzZXQiOjQsInNob3J0ZW4iOnsic291cmNlIjoyMCwidGFyZ2V0IjoyMH0sImxldmVsIjoyfV0sWzQsMSwiQVdeKiIsMSx7Im9mZnNldCI6Mywic2hvcnRlbiI6eyJzb3VyY2UiOjIwLCJ0YXJnZXQiOjIwfSwibGV2ZWwiOjJ9XSxbNSwyLCJcXHRhdSIsMSx7Im9mZnNldCI6Miwic2hvcnRlbiI6eyJzb3VyY2UiOjIwLCJ0YXJnZXQiOjIwfSwibGV2ZWwiOjJ9XV0=
    \begin{tikzcd}
	{\mathsf{SVec}^{op}\times \mathsf{SVec}} && {\mathsf{Ch}_{\ge 0}^{op}\times \mathsf{Ch_{\ge 0}}} && {\mathsf{Set}} \\
	\\
	{(\mathsf{SVec}^{op})^{\times 2}\times \mathsf{SVec}} && {(\mathsf{Ch}_{\ge 0}^{op})^{\times 2}\times \mathsf{Ch_{\ge 0}}} && {\mathsf{Ch}_{\ge 0}^{op}\times \mathsf{Ch_{\ge 0}}}
	\arrow["{N\times N}", from=1-1, to=1-3]
	\arrow["{EZ^*}"{description}, shift right=4, shorten <=16pt, shorten >=16pt, Rightarrow, from=1-1, to=3-3]
	\arrow["{\mathsf{Ch_{\ge 0}}}", from=1-3, to=1-5]
	\arrow["\rho_{\ge 0}"{description}, shift right=2, shorten <=15pt, shorten >=15pt, Rightarrow, from=1-3, to=3-5]
	\arrow["{\otimes \times id}", from=3-1, to=1-1]
	\arrow["{N\times N \times N}"', from=3-1, to=3-3]
	\arrow["{AW^*}"{description}, shift right=3, shorten <=16pt, shorten >=16pt, Rightarrow, from=3-3, to=1-1]
	\arrow["{\otimes \times id}"{description}, from=3-3, to=1-3]
	\arrow["{id \times \IHom_{\ge 0}}"', from=3-3, to=3-5]
	\arrow["\tau_{\ge 0}"{description}, shift right=2, shorten <=15pt, shorten >=15pt, Rightarrow, from=3-5, to=1-3]
	\arrow["{\mathsf{Ch_{\ge 0}}}"', from=3-5, to=1-5]
    \end{tikzcd}
    \end{adjustbox}
    \end{equation*}
    Here the composition from the bottom left corner to the top right corner along the top path is the functor on the left hand side of \eqref{eq:EA}, while the one along the bottom path is the functor on its right hand side.  Then $E$ is the top-to-bottom composition of natural transformations and $A$ is the bottom-to-top composition. The natural transformations $EZ^*$ and $AW^*$ are the precompositions of the natural maps in Eilenberg-Zilber Theorem (Prop. \ref{prop:EZ-thm-tensor}). The natural isomorphisms $\tau_{\ge 0}$ and $\rho_{\ge 0}$ come from the tensor-hom adjunction in Cor. \ref{cor:tensor-hom-positive-chains}. Clearly $E \circ A = id$ by the Eilenberg-Zilber theorem.

    Now we define 
    \begin{equation*}
    \begin{split}
        EZ^H_* &= E \circ T:
        \\
        &\Ch_{\ge 0}(N(\_), N(\IHom(\_,\_))) \rightarrow 
        \Ch_{\ge 0}(N(\_), \IHom_{\ge 0}(N(\_), N(\_))), \\
        AW^H_* &= R \circ A:\\
        &\Ch_{\ge 0}(N(\_), \IHom_{\ge 0}(N(\_), N(\_))) 
        \rightarrow 
        \Ch_{\ge 0}(N(\_), N(\IHom(\_,\_))), \\
    \end{split}
    \end{equation*}
    and 
    \begin{equation}\label{eq:EZ-AW-internal-hom}
        EZ^H = EZ^H_*(id_{N(\IHom(\huaV, \huaW))}), \qquad 
        AW^H = AW^H_*(id_{\IHom_{\ge 0}(N(\huaV), N(\huaW))}). 
    \end{equation}
    By replacing all functors with the ones in the respective homotopy categories, and using Corollaries \ref{cor:tensor-hom-hsvect}, \ref{cor:tensor-hom-hchp}, \ref{cor:DK-htpy-cats} and \ref{cor:EZ-implies-DK-homotopy-monoidal}, the above diagrams define natural \textit{isomorphisms}
    \begin{equation*}
        {\hChP}(\_, N(\IHom(\_, \_))) \newrightleftarrows{EZ^H_*}{AW^H_*} {\hChP}(\_, \IHom_{\ge 0}(N(\_), N(\_))). 
    \end{equation*}
    Because the Yoneda embedding is fully faithful, and thus reflects isomorphisms, $EZ^H$ and $AW^H$ give rise to isomorphisms in the homotopy category. Hence $EZ^H$ and $AW^H$ form a chain homotopy equivalence. Additionally, since $T\circ R=id$ and $E\circ A=id$,  $EZ^H_* \circ AW^H_* = id$. Therefore $EZ^H \circ AW^H = id$. 
\end{proof}

\subsection{Order of tensor products and internal homs between higher groupoids}

As we see in Proposition \ref{prop:EZ-thm-tensor}, the chain homotopy equivalences in the Eilenberg-Zilber Theorem are deformation retractions. This affects the order and homotopy type of tensor products of $\VS$ $n$-groupoids. The following result expands on \cite[Remark 8.5]{HoyoTrentinaglia2023}, which observes that the tensor product of two $\VS$ 1-groupoid is not generally a $\VS$ 1-groupoid.  
Similarly,  $\mathsf{2Vec}$ studied in \cite{BaezCrans2004} is also not a monoidal category, because it is not closed under the induced monoidal product from $\SVect$.

\begin{theorem}\label{thm:order-of-tensors}
Let $\huaV$ be an $n$-groupoid of order $n$ and $\huaW$ be an $m$-groupoid of order $m$. 
Then $\huaV \otimes \huaW$ has order at least $n+m$ and homotopy type at most $n+m$.
\end{theorem}

\begin{proof}
This is because $(N(\huaV) \otimes N(\huaW))_{n+m} = N(\huaV)_n \otimes N(\huaW)_m$, since the other summands vanish by the hypothesis on the order of $\huaV$ and $\huaW$. For the same reason, this is the maximal non-zero degree of the tensor product of the normalized complexes.
By the Eilenberg-Zilber Theorem \ref{prop:EZ-thm-tensor}, $EZ$ admits a left inverse, hence it must be degreewise injective. Therefore  $EZ((N(\huaV) \otimes N(\huaW))_{n+m})$ is a non-zero subspace of $N(\huaV \otimes \huaW)_{m+n}$, which must then be non-zero. 
Since $N(\huaV) \otimes N(\huaW)$ has amplitude $(0,n+m)$, its homology has at most the same amplitude. This is isomorphic to the homology of $N(\huaV \otimes \huaW)$ by the Eilenberg-Zilber homotopy equivalence, hence the homotopy type of $\huaV \otimes \huaW$ is at most $n+m$. 
\end{proof}

\begin{remark}
    Theorem \ref{thm:order-of-tensors} effectively states that, a priori, the tensor product of a $\VS$ $n$-groupoid with a $\VS$ $m$-groupoid is a groupoid of order $\ge n+m$, but it is homotopy equivalent to an $(n+m)$-groupoid. Models of this $(n+m)$-groupoid can be obtained by $DK(H(N(\huaV) \otimes N(\huaW)))$, as mentioned in Remark \ref{rem:n-type-model}. 
\end{remark}

Now we give a precise statement about the order of an internal hom $\IHom(\huaB, \huaX)$ being \textit{exactly} the order of $\huaX$.  This proof is based on the classical result in the theory of simplicial sets stating that an internal hom with target a Kan complex is Kan \cite[Thm. 6.9]{May1967}, \cite[Cor. I.5.3]{GoerssJardine2009}, \cite[\href{https://kerodon.net/tag/00TJ}{Tag 00TJ}]{kerodon}. 
We add here a new proof to show that \textit{strict} Kan conditions are also inherited from the target, and then adapt this to the case of simplicial vector spaces. 

\begin{lemma}\label{lem:hom-ngpd}
    Let $\huaB$ be a simplicial set, and $\huaX$ be an $n$-groupoid. Then $\IHom(\huaB, \huaX)$ is an $n$-groupoid.
\end{lemma}
\begin{proof}
    First of all, by the aforementioned \cite[Thm. 6.9]{May1967}, $\IHom(\huaB, \huaX)$ is an $\infty$-groupoid. Thus we only need to prove the strict Kan condition $\Kan(m, k)!$ for $m\ge n+1$. 
    
    To prove this, we use the theory of \textit{anodyne extensions} and the standard simplicial model structure on the category of simplicial sets, reviewed in e.g. \cite[Ch. I]{GoerssJardine2009}. 
    A fibration (or Kan fibration) of simplicial sets is a map with the right lifting property with respect to trivial cofibrations. In simplicial sets,  trivial cofibrations are precisely the anodyne extensions, which are the saturated class of morphisms generated by the horn inclusions $i_{m,k} : \Lambda[m,k] \to \Delta[m]$. 
    Then \cite[Cor. 4.6]{GoerssJardine2009} states that if $i:\huaK \to \huaL$ is an anodyne extension and $\huaC\to \huaD$ is an inclusion, the map between the pushout $(\huaK \times \huaD) \cup (\huaL \times \huaC)$ and $\huaL \times \huaD$ induced by the universal property is an anodyne extension. 
    If we take $\huaK = \Lambda[m,k]$, $\huaL= \Delta[m]$, $i = i_{m,k}$, and $\huaC=\huaD=\huaB$, then the map $id \times i_{m,k}: \huaB \times \Lambda[m,k] \to \huaB \times \Delta[m]$ is an anodyne extension, hence a trivial cofibration. 
    Note that for $m \ge n+1$ and any $0\le k\le m$, this is an isomorphism between the $(n-1)$-truncations $\Lambda[m,k]_{\le n-1} \to \Delta[m]_{\le n-1}$. Hence we can apply \cite[Lemma 2.14]{Pridham2013} to the $n$-groupoid $\huaX$, and $id \times i_{m,k}: \huaB \times \Lambda[m,k] \to \huaB \times \Delta[m]$, which is a trivial cofibration and an isomorphism on the $(n-1)$-truncations. Then the induced maps 
    \begin{equation*}
    \hom(\huaB \times \Delta[m], \huaX) \to \hom(\huaB \times \Lambda[m,k], \huaX)    
    \end{equation*}
    are isomorphisms. But these are precisely the horn projections of the internal hom $\IHom(\huaB, \huaX)$ for $m\ge n+1$.  
\end{proof}

\begin{theorem}\label{thm:VSMappingSpace-nGPD}
    Let $\huaU$ be a simplicial vector space, and $\huaV$ be a $\VS$ $n$-groupoid. Then the internal hom $\IHom(\huaU, \huaV)$ is a $\VS$ $n$-groupoid.
\end{theorem}
\begin{proof}
A choice of bases at each level of the simplicial vector space $\huaU$ gives a simplicial set $\huaB$ such that $\huaU \cong \R [\huaB]$.
Since a linear morphism is determined by where the base vectors go, we have
$\SVect(\huaU, \huaV) = \hom(\huaB, \huaV)$.
Similarly, we have further  $\SVect(\huaU\otimes \Delta[m], \huaV) = \hom(\huaB\times \Delta[m], \huaV)$ because $\huaU \otimes \Delta[m] \cong \R [\huaB \otimes \Delta[m]]$.   Thus we have
\begin{equation}
    \IHom_{\SVect}(\huaU, \huaV) = \IHom_{\SSet}(\huaB, \huaV). 
\end{equation}
Then Lemma \ref{lem:hom-ngpd} implies that the underlying simplicial set of $\IHom_{\SVect}(\huaU, \huaV)$ is a $n$-groupoid. Therefore $\IHom_{\SVect}(\huaU, \huaV)$ is a $\VS$ $n$-groupoid. 
\end{proof}

\section{\texorpdfstring{$n$}{n}-duals and pairings}\label{sec:ndual-and-pairings}

Since the dual space is usually defined with the help of  a canonical pairing, we now  discuss pairings of simplicial vector spaces. We consider $n$-shifted pairings, which means they take value in the $\VS$ $n$-groupoid $B^n\R$. This is the $\VS$ $n$-groupoid consisting of $\R$ on the $n$-th level and $0$ on the levels lower than $n$. Its full simplicial data is 
\begin{equation}\label{eq:def-BnR}
    (B^n\R)_m= \left\{\begin{array}{rl}
         0, & \text{for } m \leq n-1,\\
         \R^{(^m_n)}, & \text{for } m\ge n,  
    \end{array}\right. 
\end{equation} 
and its canonical multiplication maps given by the Moore fillers are 
\begin{equation}\label{eq:defMultiplicationBnR}
    m^{B^n\R}_k((a_i)_{0\le i \neq k \le n}) 
    = \sum_{i=0, \, i\neq k }^{n} (-1)^{i-k+1} a_i,
\end{equation} 
Recall that each of these multiplications is the $k$-th face map at level $n+1$ when writing $(B^n \R)_{n+1}=\Lambda^n_k(B^n\R)$. In this case, the other face maps $d^\R_{i\neq k}$ are simply projections towards the $i$-th face. The normalized chain complex of $B^n\R$ is precisely the chain complex with $\R$ in degree $n$ and $0$ in all other degrees: $N(B^n\R) = \R[-n]$. 
Topologically, $B^n\R$ can be seen as the Eilenberg-MacLane space $K(\R,n)$.

\begin{definition}\label{def:simp-pairing}
    Let $\huaV$ and $\huaW$ be simplicial vector spaces. 
    We call a simplicial linear map $\alpha: \huaV \otimes \huaW \to B^n\R$ an \textbf{$n$-shifted simplicial pairing}.
    We call a linear map $\alpha_n: \huaV_n \otimes \huaW_n \to \R$ an \textbf{$n$-shifted pairing} of $\huaV$ with $\huaW$. 
    Additionally, we say $\alpha_n$ is \textbf{multiplicative} if 
    \begin{equation}\label{eq:nShiftedPairingMultiplicative}
        \alpha_n d_0 - \alpha_n d_1 + \alpha_n d_2 - \dots + (-1)^n \alpha_n d_n = 0, 
    \end{equation}
    with $d_i = d_i^\huaV \otimes d_i^\huaW$.
    We also say $\alpha_n$ is \textbf{normalized} if
    \begin{equation}\label{eq:nShiftedPairingNormalized}
        \alpha_n s_i = 0, \quad \forall 0 \le i < n,
    \end{equation}
    with $s_i = s_i^{\huaV} \otimes s_i^{\huaW}$.
\end{definition}

\begin{remark}
    Thanks to the Finite Data Theorem \ref{thm:FiniteDataVS},
    there is a one-to-one correspondence between $n$-shifted multiplicative and normalized pairings $\alpha_n$ and $n$-shifted simplicial pairings $\alpha$. In fact, starting with the simplicial linear map $\alpha:\huaV \otimes \huaW \to B^n\R$, this is determined by its $n$-th level $\alpha_n$, whose compatibility with the face and degeneracy maps is equivalent to \eqref{eq:nShiftedPairingMultiplicative} and \eqref{eq:nShiftedPairingNormalized} respectively. Since all levels $\alpha_i$ for $i<n$ vanish, the commutativity of  $\alpha_{\le n}$ with the face maps and the commutativity of $\alpha_{\le n-1}$ with the degeneracy maps are automatic. Meanwhile, commutativity with the face maps between levels $n+1$ and $n$ translates directly to the multiplicativity condition \eqref{eq:nShiftedPairingMultiplicative} as in Theorem \ref{thm:FiniteDataVS}, and commutativity with the degeneracy maps between levels $n$ and $n-1$ translates directly to the normalization conditions \eqref{eq:nShiftedPairingNormalized}.
\end{remark}

\begin{remark}
    In the language of simplicial cohomology of $\huaV \otimes \huaW$, $\alpha_n \in C^n(\huaV \otimes \huaW)$ is multiplicative if and only if it is closed with respect to the simplicial differential $\delta_n := \sum_{i=0}^n (-1)^i (d_i^n)^*$. In other words, \eqref{eq:nShiftedPairingMultiplicative} can be written simply as $\delta \alpha_n = 0$, as it appears for $n$-shifted symplectic forms in \cite{Lesdiablerets, CuecaZhu2023}.
\end{remark} 

We now introduce the central concept of this article: the  $n$-duals. 
\begin{definition}\label{def:ndual}
    Let $\huaV$ be a simplicial vector space. For each $n\ge 0$, the \textbf{$n$-dual}  of $\huaV$, is defined by
    \begin{equation}\label{eq:vs-n-dual}
        \huaV^{n*} := \IHom(\huaV, B^n\R). 
    \end{equation}
    By Theorem \ref{thm:VSMappingSpace-nGPD} and the fact that $B^n\R$ is a $\VS$ $n$-groupoid, $\huaV^{n*}$ is a $\VS$ $n$-groupoid. Thus we also call it the \textbf{dual $n$-groupoid} of $\huaV$.

    The \textbf{$n$-dual pairing} is $\langle, \rangle = \tau_{\huaV^{n*}, \huaV, B^n \R}(id_{\huaV^{n*}}): \huaV^{n*} \otimes \huaV \to B^n \R$, where $\tau$ is the adjunction isomorphism in \eqref{eq:tensor-hom-rho-tau}. 
\end{definition}

\begin{example}[0-dual]\label{ex:0dual}
    The 0-dual $\huaV^{0*}$ of any simplicial vector space $\huaV$ is exactly the identity groupoid of $\pi_0(\huaV)^* = H^0(N(\huaV))^* = (\huaV_0/\partial(\huaV_1))^*$. In fact, since $\huaV^{0*}$ is a 0-groupoid, it is determined entirely by its level 0. Additionally, the data of a simplicial map $f:\huaV \to B^0\R$ reduces to that of an element $f\in \huaV_0^*$ such that for any $v \in \huaV_1$, $fd_0v = fd_1v$. This is the same as saying that $f$ is constant on each connected component of $\huaV$. Thus it descends to an element $[f] \in \pi_0(\huaV)^*$.
    
    In the particular case that $\huaV$ is a 0-groupoid, i.e. the identity groupoid of the vector space $\huaV_0$, $\huaV^{0*} = \huaV_0^*$, the identity groupoid of the dual of $\huaV_0$. 
    Here the $0$-dual pairing coincides with the usual dual pairing of vector spaces and it is even non-degenerate on the nose. 
    If $\huaV$ is not a 0-groupoid on the nose, but at most a 0-type, then the 0-dual is still homotopy equivalent to $\huaV$, but the 0-dual pairing is only non-degenerate up to homotopy, as we will show in Theorem \ref{thm:ndual-pairing-hom-nondeg}. 
\end{example}

The first thing that the $n$-dual allows us to do is to use the tensor-hom adjunction in Prop. \ref{prop:tensor-hom-svect} to see that the data of any simplicial $n$-shifted pairing $\alpha: \huaV \otimes \huaW \to B^n\R$ is equivalent to either of two induced simplicial linear maps
\begin{equation*}
    \alpha^l: \huaW \to \huaV^{n*}, \qquad \alpha^{r}: \huaV \to \huaW^{n*}, 
\end{equation*}
which we call the \textbf{left} and \textbf{right induced map}, respectively. While $\alpha^r$ is defined by applying the $\rho$ map in \eqref{eq:tensor-hom-rho-tau}, $\alpha^l$ is defined similarly using the other isomorphism in the flipped version of the adjunction.  In the notation of \eqref{eq:n-dual-element-in-comp}, where we consider each $n$-simplex $q = s_I d_J E_k \in \Delta[k]_n$ as indexing a component of the map at each level $k$, we write for any $w \in \huaW_k$, $v \in \huaV_n$, and any $v' \in \huaV_k$, $w' \in \huaW_n$,
\begin{equation}\label{eq:nShiftedPairingIndMapsDef}
    (\alpha^l(w))^q (v) = \alpha(v, s^{\huaW}_I d^{\huaW}_J w), \qquad (\alpha^r(v'))^q (w') = \alpha(s^{\huaV}_I d^{\huaV}_J v', w').
\end{equation}
Here we recall $I, J$ are a pair of multi-indices such that $|I|-|J| = k - n$, so that $s^{\huaW}_I d^{\huaW}_J w\in \huaW_n$ and $s^{\huaV}_I d^{\huaV}_J v' \in \huaV_n$, as expected. 

\begin{example}\label{ex:ndual-pairing-induced-maps}
    The simplicial maps induced by the $n$-dual pairing are
    \begin{equation*}
    \langle, \rangle^l: \huaV \to (\huaV^{n*})^{n*}, \qquad \langle, \rangle^r = id: \huaV^{n*} \to \huaV^{n*}.
    \end{equation*}
    Here, the right induced map is of course the identity by definition of the $n$-dual pairing as $\tau(id)$. 
    On the other hand, if $\huaV$ is at most an $n$-type, the left induced map provides the homotopy equivalence between $\huaV$ and its double $n$-dual, as we will show in Theorem \ref{thm:ndual-reflexive-uth}. 
\end{example}

In Example \ref{ex:0dual},  we see that, despite the fact that every simplicial vector space admits an $n$-dual for any $n\ge 0$, this might a priori come with a loss of information. Our aim in the following discussion is to make such a statement precise and explain why this happens. The main tools at our disposal are the Dold-Kan correspondence and the Eilenberg-Zilber Theorem for internal homs. 
Roughly speaking the problem lies in the fact that the $n$-dual pairing is related to the homotopy equivalence $EZ^H$, and this sees only a truncated version of the $n$-shifted dual in chain complexes. 
We define the latter object in Definition \ref{def:nshifted-dual-chain-cplx} and make this statement precise in Theorem \ref{thm:ndual-pairing-hom-nondeg}.

On the other side of the Dold-Kan correspondence from simplicial pairings, we have shifted pairings of chain complexes, which we call IM-pairings following \cite{CuecaZhu2023}, where the associated IM-form to a shifted 2-form on a Lie $n$-groupoid was defined. This is in reference to the theory of infinitesimally multiplicative forms and tensors appearing in \cite{BursztynCrainicWeinsteinZhu2004}, \cite{BursztynCabrera2012}, \cite{BursztynDrummond2019}.

\begin{definition}\label{def:IMPairing}
    Let $(A, \partial^A)$, $(B, \partial^B)$ be chain complexes concentrated in non-negative degrees. 
    An \textbf{$n$-shifted IM-pairing} $\lambda$ between $A$ and $B$ is a chain map $\lambda: A\otimes B \to \R[-n]$. In other words, $\lambda$ is a linear map
    \begin{equation*}
        \lambda: (A \otimes B)_n = \bigoplus_{i=0}^n A_i \otimes B_{n-i} \to \R,
    \end{equation*}
    that is \textbf{infinitesimally multiplicative}, i.e. it satisfies 
    \begin{equation}\label{eq:IMPairingDef}
        \lambda(\partial u, w) + (-1)^{i+1} \lambda(u, \partial w) = 0,
    \end{equation}
    for any $u\in A_{i+1}$ and $w \in B_{n-i}$.
\end{definition}

In chain complexes we have a natural notion of shifted duality, because of the existence of the $m$-shift operator $[m]$, which is defined for any chain complex $(A, \partial)$, by $A[m]_i = A_{i+m}$ with differential $\partial^{[m]} = (-1)^m \partial$. 
Because the $n$-shifted dual we are about to define is generally not concentrated in non-negative degrees, we also introduce its truncation as the $n$-shifted dual inside $\ChnV$.

\begin{definition}\label{def:nshifted-dual-chain-cplx}
    Let $(A, \partial)$ be a non-negative chain complex. The \textbf{$n$-shifted dual} $A^*[-n]$ of $A$ is 
    \begin{equation}\label{eq:def-nshifted-dual-chain-cplx}
        A^*[-n] := \IHom(A, \R[-n]) = \IHom(A, \R[0])[-n],
    \end{equation}
    with differential $\partial^*_i = (-1)^{i+1}(\partial_{n-i+1})^t$ with $\partial_{n-i+1}: A_{n-i+1} \to A_{n-i}$ and
    \begin{equation*}
        (\partial_{n-i+1})^t: (A^*[-n])_i = A^*_{n-i}\to A^*_{n-i+1}=(A^*[-n])_{i-1}
    \end{equation*}
    its dual.  

    The \textbf{non-negative $n$-shifted dual} $A^*[-n]_{\ge 0}$ is the truncation 
    \begin{equation}\label{eq:def-nonneg-nshifted-dual-chain-cplx}
        A^*[-n]_{\ge 0} := \IHom_{\ge 0}(A, \R[-n]) = tr_{\ge 0}(\IHom(A, \R[0])[-n]),
    \end{equation}
    with the same differential. Here $tr_{\ge 0}$ is the truncation functor defined in \eqref{eq:non-neg-chain-maps}.
\end{definition}

\begin{remark}
    The sign in the differential of the $n$-shifted dual comes from the definition of the internal hom in Section \ref{sec:monoidal-struct-chains}. 

\end{remark}

\begin{remark}\label{rem:DK-n-shift-dual-appendix}
    Using Theorem \ref{thm:EZ-thm-hom}, the non-negative $n$-shifted duals above make the $\VS$ $n$-groupoid $DK(N(\huaV)^*[-n]_{\ge 0})$ also a possibly reasonable notion of $n$-dual of $\huaV$. We discuss this in Appendix \ref{sec:appendix-DK-N-nshifted-dual} and compare it to our definition of $n$-dual in Definition \ref{def:ndual}.
\end{remark}

By the hom-tensor adjunction in chain complexes (Proposition \ref{prop:tensor-hom-chains}), any IM-pairing $\lambda: A \otimes B \to \R[-n]$ also induces two chain maps 
\begin{equation}\label{eq:IMPairingDef-induced-maps}
    \begin{split}
        &\lambda^l: B \to A^*[-n], \qquad \lambda^r: A \to B^*[-n]\\
        &\lambda^l(w)(v) := \rho^L\lambda(w)(v) = (-1)^{(n-j)j} \lambda(v, w), \text{ for } w\in B_{j}, v \in A_{n-j},\\
        &\lambda^r(v)(w) := \rho\lambda(w)(v) = \lambda(v, w), \text{ for } v \in A_i, w\in B_{n-i},
    \end{split}
\end{equation}
which we call \textbf{left} and \textbf{right induced map} respectively.

As for usual vector spaces, the induced maps above are isomorphisms precisely when the IM-pairing is non-degenerate. For chain complexes, we can additionally consider non-degeneracy of pairing in homology, and this is equivalent to the induced maps being quasi-isomorphisms. 
By the following proposition associating an IM-pairing to each simplicial pairing, this allows to define a notion of non-degeneracy up to homotopy for simplicial pairings which we summarize in Defintion \ref{def:homological-non-deg}. 
This result appeared previously in \cite{CuecaZhu2023} for $n$-shifted 2-forms, where it was interpreted in \cite[Remark 2.15]{CuecaZhu2023} as an instance of the Van Est map discussed in \cite[\S 6]{AriasAbadCrainic2011} in relation with IM-forms. The formula \eqref{eq:IMPairingAssociatedDef} originally comes from \cite{Lesdiablerets}. Here, we reformulate it for pairings in $\SVect$ in terms of the Eilenberg-Zilber map. 

\begin{proposition}\label{prop:associated-IM-pairing}
    A multiplicative normalized $n$-shifted pairing $\alpha: \huaV_n \otimes \huaW_n \to \R$ induces an associated $n$-shifted IM-pairing $\lambda_{\alpha}: (N(\huaV) \otimes N(\huaW))_n \to \R$ between the respective normalized complexes. As a chain map,  $\lambda_\alpha$ is defined by the composition
    \begin{equation}
        \lambda_\alpha: N(\huaV) \otimes N(\huaW) \overset{EZ}{\longrightarrow}
        N(\huaV \otimes \huaW) \overset{N(\alpha)}{\longrightarrow} N(B^n\R) = \R[-n],
    \end{equation}
    where $EZ$ is the Eilenberg-Zilber map from Definition \ref{def:EilenbergZilberMap}. More explicitly, $\lambda_\alpha$ is given, for any $v \in N(\huaV)_i$ and $w \in N(\huaW)_{n-i}$, by
    \begin{equation}\label{eq:IMPairingAssociatedDef}
        \lambda_\alpha(v,w) = \sum_{(\mu,\nu)\in\Sh(i,n-i)} \sign(\mu,\nu) \alpha(s_{\nu_{n-i}}\dots s_{\nu_1} v, s_{\mu_{i}}\dots s_{\mu_1} w),
    \end{equation}
    where $\Sh(i,n-i)$ is the set of $(i,n-i)$-shuffles.
\end{proposition}

\begin{proof}
   Equation \eqref{eq:IMPairingAssociatedDef} follows immediately from \eqref{eq:EZ}. However, if we use the explicit formula \eqref{eq:IMPairingAssociatedDef} as the definition of $\lambda_\alpha$, to verify that $\lambda_\alpha$ is multiplicative directly is not easy. One needs a long combinatorial calculation (see \cite[Lemma E.1]{CuecaZhu2023}).\footnote{In that article, this fact is proven for the IM-pairing associated to an $n$-shifted 2-form on a simplicial manifold. However, the combinatorics are the same.}
\end{proof}

\begin{definition}\label{def:homological-non-deg}
    An $n$-shifted pairing $\alpha: \huaV_n \otimes \huaW_n \to \R$ is \textbf{homologically non-degenerate} if its associated IM-pairing $\lambda_\alpha$ descends to a non-degenerate pairing 
    \begin{equation*}
        \lambda_\alpha: H(N(\huaV))_i \otimes H(N(\huaW))_{n-i} \to \R,
    \end{equation*}
    for any $i \in \Z$. That is, it induces an isomorphism between the homologies of the normalized complexes, up to a degree shift of $n$.
    Equivalently $\alpha$ is homologically non-degenerate if either $\lambda_\alpha^r$ or $\lambda_\alpha^l$ is a quasi-isomorphism.\footnote{If one is a quasi-isomorphism, they both are. This follows from the fact that $n$-shifted duality in unbounded chain complexes is reflexive.}
\end{definition}

Due to the fact that homological non-degeneracy of an $n$-shifted pairing requires the normalized complex of $\huaV$ to be quasi-isomorphic to the \textit{total} $n$-shifted dual, and not the truncated one, the $n$-dual pairing is only homologically non-degenerate in certain cases, which we now discuss.

\begin{theorem}\label{thm:ndual-pairing-hom-nondeg}
    Let $\huaV$ be a simplicial vector space.
    The $n$-dual pairing $\langle, \rangle: \huaV^{n*} \otimes \huaV \to B^n\R$ is homologically non-degenerate if and only if $\huaV$ is at most an $n$-type.
\end{theorem}
\begin{proof}
    By definition 
    \begin{equation*}
        \lambda_{\langle, \rangle}^r = \rho (N(\tau(\id_{\huaV^{n*}}))\circ EZ): N(\huaV^{n*}) \to N(\huaV)^*[-n].
    \end{equation*}
    From the definition in \eqref{eq:EZ-AW-internal-hom}, with $\huaW=B^n\R$, we can write
    \begin{equation*}
        EZ^H = \rho_{\ge 0}(N (\tau(\id_{\huaV^{n*}})) \circ EZ): N(\huaV^{n*}) \to N(\huaV)^*[-n]_{\ge 0}.
    \end{equation*}

    If $\huaV$ is at most an $n$-type, then $H(N(\huaV)^*[-n]_{\ge 0}) \cong H(N(\huaV)^*[-n])$. Thus $ \lambda_{\langle, \rangle}^r$ and $ EZ^H $ descend to the same map. Since $EZ^H$ is a quasi-isomorphism by Theorem \ref{thm:EZ-thm-hom},  $\lambda_{\langle, \rangle}^r$ is a quasi-isomorphism.

    On the contrary, if $\huaV$ is an $m$-type for $m > n$, the homology of $N(\huaV)^*[-n]$ must be non-zero at the negative degree $n-m$. However, the homology of $N(\huaV^{n*})$ concentrates in non-negative degrees. Therefore $\lambda^r_{\langle, \rangle}$ cannot be a quasi-isomorphism and consequently $\langle , \rangle$ cannot be homologically non-degenerate.
\end{proof}

\begin{remark}\label{rem:loss-of-info}
We can see immediately from the proof that if $\huaV$ is a $\VS$ $n$-groupoid, then the right induced map $\lambda_{\langle, \rangle}^r$ is exactly $EZ^H$. 
In a sense, this is the best case scenario, where the information of $N(\huaV)^*[-n]$ does not get lost at all, not just up to homotopy as in the case of an $n$-type. 
This is because $n$-dualization corresponds to shifted dualization on the normalized complexes but the Dold-Kan correspondence only sees non-negative degrees.

In fact, if $\huaV$ has order $m$, $N(\huaV)$ has amplitude $(0, m)$, and the $n$-shifted dual $N(\huaV)^*[-n]$ has amplitude $(n-m, n)$.
If $m\le n$, $N(\huaV)^*[-n]$ is non-negative and it coincides with the truncation $N(\huaV)^*[-n] = N(\huaV)^*[-n]_{\ge 0}$. In this case $n$-dualization causes no truncation. 
Otherwise, if $m > n$, the amplitude of $N(\huaV)^*[-n]_{\ge 0}$ is strictly smaller than that of $N(\huaV)^*[-n]$, and in particular the non-negative dual contains less information than $N(\huaV)$. But if $\huaV$ is at most an $n$-type, this process causes no loss of information up to homotopy as the homology of $N(\huaV)$ has amplitude smaller than $(0, n)$.
Otherwise, the truncation forgets some of the homology as well. Thus, for $m > n$, there is generally a loss of information, which is only avoided up to homotopy for $m$-groupoids with homotopy type at most $n$.

This phenomenon is consistent with the situation for $n$-shifted symplectic structures explained in Remark 2.16, Example 2.20, and Example 2.23 of \cite{CuecaZhu2023}. 
\end{remark}

\begin{remark}
    When extending this construction to the category of simplicial vector bundles, which has its own version of the Dold-Kan correspondence \cite{HoyoTrentinaglia2021}, extra care must be taken. Whatever the $n$-dual may be in this category, the above discussion still applies and the truncation might additionally cause problems because quotients of vector bundles may not be vector bundles again. These might cause the $n$-dual to not be representable as a simplicial vector bundle and ``fall out of the category''. We plan to discuss this in detail in \cite{CuecaRonchi-temp}. See also \cite{stefano-thesis}.
\end{remark}

\begin{example}
    As in Example \ref{ex:0dual}, the 0-dual pairing for a $\VS$ 0-groupoid is even non-degenerate at the level of spaces, and its IM-pairing is non-degenerate on chains. As we will show in Remark \ref{rem:VS1Dual-DK}, the 1-dual pairing for a $\VS$ 1-groupoid is also non-degenerate at the level of spaces, and its IM-pairing is non-degenerate on chains before taking homology. The situation changes entirely for the 2-dual pairing, as we will see in Remark \ref{rem:VS2Dual-DK}. The 2-dual paring becomes only homologically non-degenerate even for $\VS$ 2-groupoids. 
    This can be interpreted as the 2-dual containing additional information that is redundant up to homotopy.
\end{example}

\begin{theorem}\label{thm:hom-nondeg-homotopy-equiv}
    Let $\huaV$ and $\huaW$ be simplicial vector spaces and $\alpha : \huaV \otimes \huaW \to B^n\R$ an $n$-shifted simplicial pairing between them. The following diagrams commute:
    \begin{equation}\label{diag:nShiftedPairingsInducedMapsTriangle}
    % https://q.uiver.app/#q=WzAsNixbMywwLCJOKFxcaHVhV157bip9KSJdLFsyLDAsIk4oXFxodWFWKSJdLFszLDEsIk4oXFxodWFXKV4qWy1uXSJdLFswLDAsIk4oXFxodWFXKSJdLFsxLDAsIk4oXFxodWFWXntuKn0pIl0sWzEsMSwiTihcXGh1YVYpXipbLW5dIl0sWzEsMCwiTihcXGFscGhhXnIpIl0sWzEsMiwiXFxsYW1iZGFfXFxhbHBoYV5yIiwyXSxbMCwyLCJcXGxhbWJkYV97XFxsYW5nbGUsIFxccmFuZ2xlX3t0b3R9fV5yIl0sWzMsNSwiXFxsYW1iZGFfXFxhbHBoYV5sIiwyXSxbNCw1LCJcXGxhbWJkYV97XFxsYW5nbGUsIFxccmFuZ2xlX3t0b3R9fV5yIl0sWzMsNCwiTihcXGFscGhhXmwpIl1d
    \begin{tikzcd}[ampersand replacement=\&,cramped]
	{N(\huaW)} \& {N(\huaV^{n*})} \& {N(\huaV)} \& {N(\huaW^{n*})} \\
	\& {N(\huaV)^*[-n]} \&\& {N(\huaW)^*[-n]}
	\arrow["{N(\alpha^l)}", from=1-1, to=1-2]
	\arrow["{\lambda_\alpha^l}"', from=1-1, to=2-2]
	\arrow["{\lambda_{\langle, \rangle_{}}^r}", from=1-2, to=2-2]
	\arrow["{N(\alpha^r)}", from=1-3, to=1-4]
	\arrow["{\lambda_\alpha^r}"', from=1-3, to=2-4]
	\arrow["{\lambda_{\langle, \rangle_{}}^r}", from=1-4, to=2-4]
    \end{tikzcd}
    \end{equation}
    Additionally, if $\huaW$ is at most an $n$-type, then $\alpha$ is homologically non-degenerate if and only if $\alpha^r$ is a weak equivalence. 
    Analogously, if $\huaV$ is at most an $n$-type, then $\alpha$ is homologically non-degenerate if and only if $\alpha^l$ is a weak equivalence. 
\end{theorem}

\begin{proof}
    The commutativity of the right-hand diagram follows from the fact that $\rho EZ^* N \tau$ can be seen as a natural transformation between the functors $\SVect^{op} \times \SVect^{op} \times \SVect \to \Set$ given by
    \begin{equation*}
        \SVect(\_, \IHom(\_, \_)) \to \Ch(N(\_), \IHom(N(\_),N(\_))).
    \end{equation*}
    Evaluating the third argument at $B^n\R$ gives a natural transformation 
    \begin{equation*}
        \SVect(\_, \_^{n*}) \to \Ch(N(\_), N(\_)^*[-n]),
    \end{equation*}
    where both are functors $\SVect^{op} \times \SVect^{op} \to \Set$. 
    The naturality square at the map $((\alpha^r)^*, id)$ in $\SVect^{op} \times \SVect^{op}$ is 
    % https://q.uiver.app/#q=WzAsNCxbMCwwLCJcXFNWZWN0KFxcaHVhV157bip9LCBcXGh1YVdee24qfSkiXSxbMSwwLCJcXENoKFxcVmVjdCkoTihcXGh1YVdee24qfSksIE4oXFxodWFXKV4qWy1uXSkiXSxbMCwxLCJcXFNWZWN0KFxcaHVhViwgXFxodWFXXntuKn0pIl0sWzEsMSwiXFxDaChcXFZlY3QpKE4oXFxodWFWKSwgTihcXGh1YVcpXipbLW5dKSJdLFsyLDMsIlxccmhvIEVaXipOXFx0YXUiLDJdLFswLDEsIlxccmhvIEVaXipOXFx0YXUiXSxbMCwyLCIoKFxccmhvXFxhbHBoYSleKiwgaWReKikiLDJdLFsxLDMsIigoTihcXHJob1xcYWxwaGEpKV4qLCBpZF4qKSJdXQ==
    \[\begin{tikzcd}[ampersand replacement=\&,cramped,row sep=scriptsize]
	{\SVect(\huaW^{n*}, \huaW^{n*})} \& {\Ch(N(\huaW^{n*}), N(\huaW)^*[-n])} \\
	{\SVect(\huaV, \huaW^{n*})} \& {\Ch(N(\huaV), N(\huaW)^*[-n])}
	\arrow["{\rho EZ^*N\tau}", from=1-1, to=1-2]
	\arrow["{((\alpha^r)^*, id)}"', from=1-1, to=2-1]
	\arrow["{((N(\alpha^r)^*, id)}", from=1-2, to=2-2]
	\arrow["{\rho EZ^*N\tau}"', from=2-1, to=2-2]
    \end{tikzcd}\]
   Applying this to  $id_{\huaW^{n*}} \in \SVect(\huaW^{n*}, \huaW^{n*})$ and together with \eqref{eq:IMPairingDef-induced-maps}, we have
    \begin{equation*}
        \lambda_{\langle, \rangle}^r \circ N(\alpha^r) = (\rho EZ^* N \tau)(id_{\huaW^{n*}}) \circ N(\alpha^r) = \rho EZ^* N \tau ((\alpha^r)^*, id) (id_{\huaW^{n*}}) = \rho EZ^* N \tau ( \alpha^r) = \lambda_{\alpha}^r.
    \end{equation*}

    Commutativity of the left-hand diagram follows analogously.

    The other two statements follow from Theorem \ref{thm:ndual-pairing-hom-nondeg}, Definition \ref{def:homological-non-deg}, the two-out-of-three property for quasi-isomorphisms, and Proposition \ref{prop:equivalences-in-svect}.
\end{proof}

\begin{remark}\label{rem:hndg-pairing-admitted}
    In the setting of the theorem, if either $\huaV$ or $\huaW$ is at most an $n$-type and $\alpha$ is homologically non-degenerate, then the other is also at most an $n$-type because it is weak equivalent to a $\VS$ $n$-groupoid. This puts restrictions on which simplicial vector spaces admit homologically non-degenerate pairings with respect to a certain shift. 
    This is analogous to the fact discussed in \cite[Remark 2.16]{CuecaZhu2023} that an $m$-shifted symplectic Lie $n$-groupoid with $m < n$ must have certain vanishing homology groups, which make it interpretable as an $m$-shifted symplectic Lie $m$-groupoid with added singularities. 
\end{remark}

\begin{remark}
    A trivial but perhaps illustrative fact to observe is that in the case of the $n$-dual pairing $\alpha= \langle, \rangle: \huaV^{n*} \otimes \huaV \to B^n\R$ the right-hand diagram in \eqref{diag:nShiftedPairingsInducedMapsTriangle} becomes tautological and the theorem only reiterates the fact that if $\huaV$ is an $n$-type to begin with, then $\lambda^r_{\langle, \rangle}$ is a quasi-isomorphism and the $n$-dual pairing is homologically non-degenerate. 
    
    On the contrary, the left-hand diagram yields an interesting property which deserves to be a theorem of its own.
\end{remark}

\begin{theorem}[$n$-duality is reflexive up to homotopy]\label{thm:ndual-reflexive-uth}
    Let $\huaV$ be at most an $n$-type. The double $n$-dual $(\huaV^{n*})^{n*}$ is weak equivalent to $\huaV$ itself via the weak equivalence given by
    \begin{equation} \label{eq:double-W.E.}
        \langle, \rangle^l: \huaV \xrightarrow{\simeq} (\huaV^{n*})^{n*}. 
    \end{equation}
\end{theorem}

\begin{proof}
    By Theorem \ref{thm:hom-nondeg-homotopy-equiv} and Theorem \ref{thm:ndual-pairing-hom-nondeg}, the left induced map in \ref{ex:ndual-pairing-induced-maps} is a weak equivalence.
\end{proof}

\begin{remark}
    By Proposition \ref{prop:equivalences-in-svect}, in the above situation, $\langle, \rangle^l: \huaV \to (\huaV^{n*})^{n*}$ is also a homotopy equivalence. 
\end{remark}

\section{Computing \texorpdfstring{$n$}{n}-duals of \texorpdfstring{$\VS$}{VS} \texorpdfstring{$n$}{n}-groupoids}\label{sec:computations}

To calculate explicitly what the $n$-dual \eqref{eq:vs-n-dual} is in general involves solving many linear equations. 
Determining the solution space of these linear equations for a general $n$ is not a trival task due to the large number of equations involved. 
In this section we give an overview of the equations for a general $n$ and set out to solve them for $n=1$ and $n=2$. For $n=1$, we rediscover Pradines's dual of $\VB$ groupoids \cite{Pradines1988} (see also \cite[\S 11.2]{Mackenzie2005} \cite[\S 4]{GraciaSazMehta2017}) applied to the case where the base is a point. For $n=2$, we make a completely new discovery. 

\subsection{Overview of the general computation for the \texorpdfstring{$n$}{n}-dual of a \texorpdfstring{$\VS$}{VS} \texorpdfstring{$n$}{n}-groupoid}\label{sec:ndual-overview}

Let $\huaV$ be a $\VS$ $n$-groupoid. We now write down the linear equations to compute $\huaV^{n*}$ for a general $n$. First of all, since $\huaV^{n*}$ is a $\VS$ $n$-groupoid, its full data is contained in its $n$-truncation by Theorem \ref{thm:FiniteDataVS}. Even though the multiplications of its $n$-simplices can be recovered by computing Moore fillers, these have a nice description in terms of the $n$-dual pairing, which we later discuss.

At any fixed level $l$, since $B^n\R$ is a $\VS$ $n$-group, by Theorem \ref{thm:FiniteDataVS}, a simplicial map $f \in \huaV^{n*}_l = \Simp\Vect(\huaV \otimes\Delta[l], B^n\R)$ is determined by its first $n+1$ levels and a multiplicativity condition.  More precisely, for any choice of $0\le k \le n+1$,
we have
% https://q.uiver.app/?q=WzAsNixbMCwwLCJcXExhbWJkYV57bisxfV9rKFxcaHVhVl9cXGJ1bGxldClcXG90aW1lcyBcXERlbHRhW21dX3tuKzF9Il0sWzEsMCwiXFxodWFWX25cXG90aW1lcyBcXERlbHRhW21dX24iXSxbMSwxLCJcXFIiXSxbMCwxLCJcXFJee24rMX0iXSxbMiwxLCIwIl0sWzIsMCwiKFxcaHVhVl9cXGJ1bGxldFxcb3RpbWVzIFxcRGVsdGFbbV0pX3tcXGxlIG4tMX0iXSxbMCwzLCIoZl9uLFxcZG90cyAsZl9uKSJdLFswLDEsIiIsMSx7Im9mZnNldCI6LTJ9XSxbMCwxLCJcXGRvdHMiLDFdLFswLDEsIiIsMSx7Im9mZnNldCI6Mn1dLFszLDIsIlxcZG90cyIsMV0sWzMsMiwiIiwxLHsib2Zmc2V0IjotMn1dLFszLDIsIiIsMSx7Im9mZnNldCI6Mn1dLFsxLDIsImZfbiJdLFsyLDQsIlxcZG90cyIsMV0sWzEsNSwiXFxkb3RzIiwxXSxbNSw0LCJmX3tcXGxlIG4tMX0gPTAiXSxbMiw0LCIiLDAseyJvZmZzZXQiOi0yfV0sWzIsNCwiIiwxLHsib2Zmc2V0IjoyfV0sWzEsNSwiIiwxLHsib2Zmc2V0IjotMn1dLFsxLDUsIiIsMSx7Im9mZnNldCI6Mn1dXQ==
\begin{equation}\label{diag:VSnDualMultiplicativity}
    \begin{tikzcd}[ampersand replacement=\&]
	{\Lambda^{n+1}_k(\huaV)\otimes \Delta[l]_{n+1}} \& {\huaV_n\otimes \Delta[l]_n} \& {(\huaV\otimes \Delta[l])_{\le n-1}} \\
	{\R^{n+1}} \& \R \& 0
	\arrow["{f_{n+1}=(f_n d_0,\dots, \widehat{f_n d_k}, \dots ,f_n d_n)}"', from=1-1, to=2-1]
	\arrow[shift left=2, from=1-1, to=1-2]
	\arrow["\dots"{description}, from=1-1, to=1-2]
	\arrow[shift right=2, from=1-1, to=1-2]
	\arrow["\dots"{description}, from=2-1, to=2-2]
	\arrow[shift left=2, from=2-1, to=2-2]
	\arrow[shift right=2, from=2-1, to=2-2]
	\arrow["{f_n}", from=1-2, to=2-2]
	\arrow["\dots"{description}, from=2-2, to=2-3]
	\arrow["\dots"{description}, from=1-2, to=1-3]
	\arrow["{f_{\le n-1} =0}", from=1-3, to=2-3]
	\arrow[shift left=2, from=2-2, to=2-3]
	\arrow[shift right=2, from=2-2, to=2-3]
	\arrow[shift left=2, from=1-2, to=1-3]
	\arrow[shift right=2, from=1-2, to=1-3]
    \end{tikzcd}
\end{equation}
where the components $f_{\le n-1}$ vanish because $B^n\R_{\le n-1}= 0$. Thus $f$ is uniquely determined by $f_n$.
As in Section \ref{sec:Svect-int-hom},  we write $f_n$ in the components
\begin{equation}\label{eq:n-dual-element-in-comp}
    f_n^q (v) := f_n(v^q),
    \text{ for all $v^q \in \huaV_n^q \subset \huaV_n \otimes \Delta[l]_n$ and $q\in \Delta[l]_n$ }.
\end{equation}
With this, $f$ is determined by the $\binom{l+n+1}{l}$-tuple of linear functions $(f_n^q: \huaV_n \to \R)_{q \in \Delta[l]_n}$ which satisfies the following two sets of equations coming from commutativity of \eqref{diag:VSnDualMultiplicativity}:
\begin{equation}\label{eq:MultiplicativityMostGeneral}
    f_n^{d_k r}(m_k((v_i)_{0\le i \neq k \le n})) = \sum_{i=0, \, i\neq k }^{n} (-1)^{i-k+1} f^{d_i r}(v_i), \quad \forall r\in \Delta[l]_{n+1}, \; \forall (v_i)_{0\le i \neq k \le n} \in \Lambda^{n+1}_k(\huaV), 
\end{equation}
and
\begin{equation}\label{eq:NormalizationMostGeneral}
    f_n^{s_i p}(s_i v) = f_n((s_i v)^{s_i p})= f_n(s_i(v^p)) = 0 , \quad 0\le i \le n-1,  \forall v \in \huaV_{n-1}, \; \forall p \in \Delta[l]_{n-1}.
\end{equation}

The set of $\binom{l+n+2}{l}$ equations in \eqref{eq:MultiplicativityMostGeneral}, one for each $r\in \Delta[l]_{n+1}$, encodes precisely the multiplicativity condition of Theorem \ref{thm:FiniteDataVS}. Thus we call these equations \textbf{multiplicativity conditions}.
They correspond to the fact that $f$ commutes with face maps in \eqref{diag:VSnDualMultiplicativity}.\footnote{For a fixed $k$, as we see in diagram \eqref{diag:VSnDualMultiplicativity}, $f_{n+1}$ automatically commutes with $d_{i\neq k}$ by definition. Thus the only nontrivial condition to impose is to ask $f_{n+1}$ to commute with $d_k$. Writing this down in $\binom{l+n+2}{l}$ components, we obtain eq. \eqref{eq:MultiplicativityMostGeneral}.}
In particular because $\huaV$ is a $\VS$ $n$-groupoid, by \eqref{eq:n+1truncationVSnGpd}, if \eqref{eq:MultiplicativityMostGeneral} holds for a certain $k\in [0,n]$, then it also holds for all other $k$.\footnote{This comes up in the computation of the 2-dual in \ref{sec:computation-2dual}, compare e.g. equations \eqref{eq:VS2DualMult0012with0} and \eqref{eq:VS2DualMult0122with0}}
The set of $n\cdot\binom{l+n}{l}$ equations in \eqref{eq:NormalizationMostGeneral}, one for each $q\in \Delta[l]_{n-1}$ and $i\in \{0, \dots, n\}$,  is the set of \textbf{normalization conditions}. They correspond to the fact that $f$ commutes with the degeneracy maps in the right-hand square of \eqref{diag:VSnDualMultiplicativity}. In other words, equations \eqref{eq:MultiplicativityMostGeneral} and \eqref{eq:NormalizationMostGeneral} represent, respectively, the commutative diagrams for the only non-trivial face map and that for the degeneracy maps, i.e.
\[\begin{tikzcd}[ampersand replacement=\&]
	{\Lambda^{n+1}_k(\huaV)^{r}} \& {\huaV_n^{d_k r}} \\
	{\R^{n+1}} \& \R
	\arrow["{(f^{d_0 r}, \dots, \widehat{f^{d_kr}},\dots, f^{d_n r})}"', from=1-1, to=2-1]
	\arrow["{m_k}", from=1-1, to=1-2]
	\arrow["{f^{d_k r}}", from=1-2, to=2-2]
	\arrow["{m_k^\R}"', from=2-1, to=2-2]
\end{tikzcd}, \quad
\begin{tikzcd}[ampersand replacement=\&]
	{\huaV_n\otimes \Delta[l]_n} \& {(\huaV_{n-1}\otimes \Delta[l]_{n-1})} \\
	\R \& 0
	\arrow["{f_n}", from=1-1, to=2-1]
	\arrow["{s_i = 0}"', from=2-2, to=2-1]
	\arrow["{f_{n-1} =0}", from=1-2, to=2-2]
	\arrow["{s_i}"', from=1-2, to=1-1]
\end{tikzcd}\]where $m_k^\R$ is given in \eqref{eq:defMultiplicationBnR}. 

The $n$-dual $\huaV^{n*}$ is always a $\VS$ $n$-groupoid. In this case,  by Theorem \ref{thm:FiniteDataVS},  $\check{m}_k=d_k: \huaV^*_{n+1}\cong \Lambda^{n+1}_k(\huaV^{n*}) \to \huaV^{n*}_n$ is the $k$-th face map. Take  \eqref{eq:MultiplicativityMostGeneral} with $l=n+1$ and $r=E_{n+1}=01\dots n(n+1)$, the unique non-degenerate $n+1$-simplex in $\Delta[n+1]$. Notice that $d_i E_{n+1} = \delta_i E_n$ and  $f^{\delta_i E_{n+1}} = (\widecheck{d}_i f)^{E_{n+1}}$ which are precisely the interior components of each face of an $n+1$-simplex $f \in \huaV^{n*}_{n+1}$. 
Solving for each $(d_k f)^{E_n}$ determines the map $\widecheck{m}_k$ in terms of the interiors of the other faces of $f$, which form an $(n+1, k)$-horn. 
Therefore, in terms of the $n$-dual pairing, the equation \eqref{eq:MultiplicativityMostGeneral} for $r=E_{n+1}$ is exactly the multiplicativity condition of $\langle, \rangle$ as a pairing, defined in \eqref{eq:nShiftedPairingMultiplicative}. That is, for any $(n+1, k)$-horns $(f_i)_{0\le i \neq k \le n} \in \Lambda^{n+1}_k(\huaV^{n*})$ and $(v_i)_{0\le i \neq k \le n} \in \Lambda^{n+1}_k(\huaV)$,
\begin{equation}\label{eq:ndual-mult-pairing}
    \langle \widecheck{m}_k((f_i)_{0\le i \neq k \le n}), m_k((v_i)_{0\le i \neq k \le n})\rangle = \sum_{i=0, \, i\neq k }^{n} (-1)^{i-k+1} \langle f_i , v_i \rangle.
\end{equation}
In other words, each multiplication of the $n$-dual can be defined by being the unique multiplication $\Lambda^{n+1}_k(\huaV^{n*}) \to \huaV^{n*}_n$ that makes the $n$-dual pairing multiplicative. This is how the multiplication of the dual $\VB$ groupoid is defined in \cite{Pradines1988} and later references, while \eqref{eq:ndual-mult-pairing} was obtained for the first time by other means in \cite{CosteDazordWeinstein1987} as the multiplication of the cotangent groupoid.

As a final comment, the fact that $\langle, \rangle$ is normalized as a pairing as in \eqref{eq:nShiftedPairingNormalized} is equivalent to the normalization conditions \eqref{eq:NormalizationMostGeneral} for $p = E_{n-1}$. 
In fact, by the definition of the degeneracy maps in equation \eqref{eq:HomSpaceMaps} and the fact that $s_i E_{n-1} = \sigma_i E_n$, $\forall i \in [0,  n-1], f_n \in \huaV^{n*}_n, v\in \huaV_n$,  the $E_n$ component of $\widecheck{s}_i f_n$ is given by $f^{\sigma_i(E_{n-1})}_n$, and 
\begin{equation*}
    \langle \widecheck{s}_i f_n, s_i v \rangle = f_n^{\sigma_iE_{n}}(s_i v) = f_n^{s_iE_{n-1}}(s_i v) = 0.
\end{equation*}

\subsubsection{Annihilators and $(\ker d_i)^*$ }\label{sec:annihilators-and-dual-kernels}

This section is devoted to the study of the solution spaces of the normalization equations \eqref{eq:NormalizationMostGeneral}. 
As a motivating example, in the case of a $\VB$ 1-groupoid, the dual of the \textbf{core} appears in the definition of the dual $\VB$ groupoid. This is an object with multiple isomorphic descriptions (c.f. \cite[\S 3.2.1-3]{GraciaSazMehta2017}). 
First of all, the core can be defined as the kernel of either the source or the target map of a $\VB$ groupoid: $\ker \widetilde{d}_0^1$ is known as the right core and $\ker \widetilde{d}_1^1$ is known as the left core. 
These are isomorphic through the involution mentioned in Example \ref{ex:MooreFillers1Gpd}. Additionally, their duals are isomorphic to the annihilator of the units, $\Ann(\widetilde{s}_0\huaV_0)\subseteq \huaV_1^*$, which in the case of the tangent $\VB$ groupoid is the conormal bundle to the units of the base.
To compute $n$-duals it is convenient to have a systematic understanding of these isomorphisms, so we introduce the following terminology in analogy with this example.

\begin{definition}
Let $\huaV$ be a simplicial vector space. 

For any subset $A \subseteq [m] \in \Delta$ of cardinality $|A|=j$, the associated \textbf{degree $j$ $m$-dimensional core} is $\bigcap_{i \in A} \ker d^m_i \subseteq \huaV_m$.\footnote{When defining generalized horns as in \cite{Joyal2008}, \cite{stefano-thesis}, this is the kernel of a horn projection: $\ker p^m_{[m]/A}$.}

For any subset $B \subseteq [m-1] \in \Delta$ of cardinality $|B|=j$, the associated \textbf{degree $j$ $m$-dimensional degeneracy annihilator} is $O_B := \Ann\left(\sum_{i\in B} s_i \huaV_{m-1}\right) \subseteq \huaV_{m}^*$. 
\end{definition}

We now focus on the two extreme cases needed to compute 1- and 2-duals, which we use in Section \ref{sec:computation-1dual} and Section \ref{sec:computation-2dual}. These are degree 1 $n$-dimensional cores and degree $n$ $n$-dimensional cores. 

Since different annihilators appear separately as solution spaces of the normalization equations \eqref{eq:NormalizationMostGeneral}, we refrain from identifying them as one object. This will further pay off in the computation of the $\VB$ $n$-duals which we carry out in \cite{stefano-thesis, CuecaRonchi-temp}. 

The degree $1$ $n$-dimensional degeneracy annihilators are the spaces that contain the components of an element in the $n$-dual that are normalized with respect to a single degeneracy map, $O_i:=\Ann(s_i\huaV_{n-1})$. We observe that for $0\le i < n$, the following exact sequences canonically split:
\begin{equation}\label{eq:VSnDualSESd_is_i}
% https://q.uiver.app/#q=WzAsNSxbMSwwLCJcXGtlciBkX2kiXSxbMiwwLCJcXGh1YVZfbiJdLFswLDAsIjAiXSxbMywwLCJcXGh1YVZfe24tMX0iXSxbNCwwLCIwIl0sWzAsMSwiIiwxLHsib2Zmc2V0IjotMSwic3R5bGUiOnsidGFpbCI6eyJuYW1lIjoiaG9vayIsInNpZGUiOiJ0b3AifX19XSxbMiwwXSxbMSwzLCJkX2kiLDAseyJvZmZzZXQiOi0xfV0sWzMsMSwic19pIiwwLHsib2Zmc2V0IjotMX1dLFszLDRdLFsxLDAsImlkLSBzX2lkX2kiLDAseyJvZmZzZXQiOi0xfV1d
\begin{tikzcd}[ampersand replacement=\&]
	0 \& {\ker d_i} \& {\huaV_n} \& {\huaV_{n-1}} \& 0
	\arrow[shift left=1, hook, from=1-2, to=1-3]
	\arrow[from=1-1, to=1-2]
	\arrow["{d_i}", shift left=1, from=1-3, to=1-4]
	\arrow["{s_i}", shift left=1, from=1-4, to=1-3]
	\arrow[from=1-4, to=1-5]
	\arrow["{id- s_id_i}", shift left=1, from=1-3, to=1-2],
\end{tikzcd}
\end{equation}
\begin{equation}\label{eq:VSnDualSESd_iplus1s_i}
    % https://q.uiver.app/#q=WzAsNSxbMSwwLCJcXGtlciBkX3tqKzF9Il0sWzIsMCwiXFxodWFWX24iXSxbMCwwLCIwIl0sWzMsMCwiXFxodWFWX3tuLTF9Il0sWzQsMCwiMCJdLFswLDEsIiIsMSx7Im9mZnNldCI6LTEsInN0eWxlIjp7InRhaWwiOnsibmFtZSI6Imhvb2siLCJzaWRlIjoidG9wIn19fV0sWzIsMF0sWzEsMywiZF97aisxfSIsMCx7Im9mZnNldCI6LTF9XSxbMywxLCJzX3tqfSIsMCx7Im9mZnNldCI6LTF9XSxbMyw0XSxbMSwwLCJpZC0gc197an1kX3tqKzF9IiwwLHsib2Zmc2V0IjotMX1dXQ==
\begin{tikzcd}[ampersand replacement=\&] 
	0 \& {\ker d_{i+1}} \& {\huaV_n} \& {\huaV_{n-1}} \& 0
	\arrow[shift left=1, hook, from=1-2, to=1-3]
	\arrow[from=1-1, to=1-2]
	\arrow["{d_{i+1}}", shift left=1, from=1-3, to=1-4]
	\arrow["{s_{i}}", shift left=1, from=1-4, to=1-3]
	\arrow[from=1-4, to=1-5]
	\arrow["{id- s_{i}d_{i+1}}", shift left=1, from=1-3, to=1-2].
\end{tikzcd}
\end{equation}
This implies that $\huaV_n \cong \ker d_i \oplus s_i(\huaV_{n-1}) \cong \ker d_{i+1} \oplus s_i(\huaV_{n-1})$.  As a result, we have the isomorphisms
\begin{equation}
   O_i = \Ann(s_i\huaV_{n-1}) \cong (\ker d_i)^* \cong (\ker d_{i+1})^*, \quad 0\le i \le n. 
\end{equation}
Different choices of isomorphic descriptions of the $O_i$ lead to different explicit descriptions of the $n$-dual. In the dual $\VB$ groupoid example this can be seen by comparing the one obtained in \cite[\S 11.2]{Mackenzie2005} by choosing the right core (kernel of the source $d_0$), and the one obtained in \cite[\S 2]{Pradines1988}, by choosing the annihilator of the units of the $\VB$ groupoid.

In the other extreme case, the degree $n$ $n$-dimensional degeneracy annihilator is the space that contains fully normalized components
\begin{equation*}
    O_{01\dots (n-1)} := \Ann(D_n\huaV) = \bigcap_{i=0}^n O_i,
\end{equation*}
where $D_n\huaV = s_0\huaV_{n-1} + \dots + s_{n-1}\huaV_{n-1}$ is the space generated by all degenerate $n$-simplices in $\huaV$. As observed in Lemma \ref{lem:MooreFillers-degeneracies}, $D_n\huaV$ is isomorphic to the horn space $\Lambda^n_k(\huaV)$ via $\mu^n_k$, for any $0 \le k \le n$. The degree $n$ $n$-dimensional core $\ker p^n_k$ appears in the exact sequence
% https://q.uiver.app/#q=WzAsNSxbMSwwLCJcXGtlciBwXm5fayJdLFsyLDAsIlxcaHVhVl9uIl0sWzAsMCwiMCJdLFszLDAsIlxcTGFtYmRhXm5fayhcXGh1YVZfXFxidWxsZXQpIl0sWzQsMCwiMCJdLFswLDEsIiIsMSx7Im9mZnNldCI6LTEsInN0eWxlIjp7InRhaWwiOnsibmFtZSI6Imhvb2siLCJzaWRlIjoidG9wIn19fV0sWzIsMF0sWzEsMywicF5uX2siLDAseyJvZmZzZXQiOi0xfV0sWzMsMSwiXFxtdV5uX2siLDAseyJvZmZzZXQiOi0xfV0sWzMsNF0sWzEsMCwiXFxnYW1tYV5uX2sgPSBpZC0gXFxtdV5uX2siLDAseyJvZmZzZXQiOi0xfV1d
\begin{equation}\label{eq:VSnDualSESp_kmu_k}
\begin{tikzcd}[ampersand replacement=\&,cramped]
	0 \& {\ker p^n_k} \& {\huaV_n} \& {\Lambda^n_k(\huaV)} \& 0,
	\arrow[from=1-1, to=1-2]
	\arrow[shift left, hook, from=1-2, to=1-3]
	\arrow["{\gamma^n_k = id- \mu^n_k}", shift left, from=1-3, to=1-2]
	\arrow["{p^n_k}", shift left, from=1-3, to=1-4]
	\arrow["{\mu^n_k}", shift left, from=1-4, to=1-3]
	\arrow[from=1-4, to=1-5]
\end{tikzcd}
\end{equation}
where we call the retract $\gamma^n_k:= id-\mu^n_k$ the \textbf{$k$-th (degree $n$) core projection}. These projections also provide isomorphisms between cores for different indices $k$. From the dual sequence we have
\begin{equation*}
    O_{01\dots (n-1)} = \Ann(D_n\huaV) \cong (\ker p^n_k)^*, \qquad \forall 0 \le k \le n,
\end{equation*}
where the isomorphisms are given by the $\gamma^*_k$. 

\begin{remark}\label{rem:general-degree-cores}
    By using the algorithm in the proof of Prop. \ref{thm:MooreHornFillers} we expect to be able to construct fillers for arbitrary generalized horns and obtain sequences such as \eqref{eq:VSnDualSESp_kmu_k}. This would allow to classify all degree $j$ cores.
    We plan to study this in more detail in future work.
\end{remark}

\subsection{The 1-dual}\label{sec:computation-1dual}

Let $\huaV$ be a $\VS$ 1-groupoid $\huaV_1 \rightrightarrows \huaV_0$.
As remarked in Section \ref{sec:annihilators-and-dual-kernels}, the normalized components of each $m$-simplex in the $n$-dual are elements of an annihilator space, which is isomorphic to the dual of each core. For the $1$-dual, only the degree 1 1-dimensional degeneracy annihilator $O_0= \Ann(s_0\huaV_0)$ is relevant. Recall that
\begin{equation*}
    O_0 = \Ann(s_0\huaV_0) \cong (\ker d_0)^* \cong (\ker d_1)^* \subseteq \huaV_1^*.
\end{equation*}
In the literature (\cite{Mackenzie2005}, \cite{GraciaSazMehta2017}), $\ker d_0$ is known as the right core and $\ker d_1$ as the left core. The core projections $\gamma_0: \huaV_1 \to \ker d_0$ and $\gamma_1:\huaV_1 \to \ker d_1$ appearing in \eqref{eq:VSnDualSESp_kmu_k} for $p^1_1=d_0$ and $p^1_0=d_1$ are 
\begin{equation*}
    \gamma_0v = v - 1d_0v, \qquad \gamma_1v = v - 1d_1v. 
\end{equation*}
Their dual maps $\gamma_i^*: (\ker d_i)^* \to O_0$ give the isomorphisms between $O_0$ and $(\ker d_i)^*$. Note that the isomorphism between left and right core is also given by restriction of the appropriate projection, e.g. $\gamma_1|_{\ker d_0}: \ker d_0 \overset{\cong}{\to} \ker d_1$. As in Example \ref{ex:MooreFillers1Gpd} this happens to be the opposite of the groupoid inversion: $\gamma_1 c = - c^{-1} = d_2\mu^2_2(c, 0)$, for any $c \in \ker d_0$. 

We summarize the computation of $\huaV^{1*}$ in the following proposition.

\begin{proposition}[1-dual of a $\VS$ 1-groupoid] \label{prop:VS1dual}
    Let $\huaV$ be a $\VS$ 1-groupoid. Then its 1-dual $\huaV^{1*}$ is
    \begin{equation*}
            \huaV_1^* 
            \rightrightarrows 
            O_0,  
    \end{equation*}
    with face maps given for any $\xi \in \huaV_1^*$ and $v \in \huaV_1$ by
    \begin{equation}
    \begin{aligned}
    \widecheck{d}_0\xi(v) = \gamma_1^*\xi(v) = \xi(v - s_0d_1v) ,\quad 
    \widecheck{d}_1\xi(v) = \gamma_0^*\xi(v) = \xi(v - s_0d_0v) ,
    \end{aligned}
    \end{equation}
    and unit map  $\widecheck{s}_0: O_0 \to \huaV_1^*$ the inclusion. 
    The multiplication $\widecheck{m}_1$ is
    \begin{equation}
        \langle \widecheck{m}_1(\eta,\xi), m_1(v,w) \rangle 
        = \langle \xi \cdot \eta, w \cdot v \rangle 
        = \langle \eta, v\rangle + \langle \xi, w \rangle,
    \end{equation}
    for any composable pairs $(\eta, \xi) \in \Lambda^2_1(\huaV^{1*})$, $(v,w)\in \Lambda^2_1(\huaV)$. 
\end{proposition}
\begin{proof}
According to Section \ref{sec:ndual-overview}, to compute each level of the 1-dual $\huaV^{1*}$ we only have to solve the linear equations given by the multiplicativity  \eqref{eq:MultiplicativityMostGeneral} and normalization \eqref{eq:NormalizationMostGeneral} conditions.  

Beginning with level 0, since the only 1-simplex in $\Delta[0]$ is $00$, we have that an element $f \in \huaV^{1*}_0 = \Hom(\huaV, B^1\R)$ only has one component $f^{00} = f$. This has to satisfy only one normalization condition \eqref{eq:NormalizationMostGeneral} for $i=0$ and $p = 0 \in \Delta[0]_0$ and one multiplicativity condition \eqref{eq:NormalizationMostGeneral} for $r=000 \in \Delta[0]_2$.
The normalization condition reads $f(s_0x) = 0$ for all $x \in \huaV_0$, therefore $f \in O_0$. 
The multiplicativity condition reads $f(w\cdot v) = f(w) + f(v)$, for any $(v,w) \in \Lambda^2_1(\huaV)$. Since $w\cdot v = d_1\mu^2_1(v,w) = v + w - 1d_0w$ as in Example \ref{ex:MooreFillers1Gpd}, this is automatically satisfied by linearity and normalization of $f$. Therefore $\huaV^{1*}_0 = O_0$. 

For level 1, any $\xi\in \huaV^{1*}_1$ has three components, $\xi^{00}$, $\xi^{01}$ and $\xi^{11}$. 
There are two normalization conditions given by \eqref{eq:NormalizationMostGeneral} for each of the two 0-simplices $p = 0,1$ in $\Delta[1]$, and four multiplicativity conditions given by \eqref{eq:MultiplicativityMostGeneral} for each of the four 2-simplices $r=000, 001, 011, 111$ in $\Delta[1]$.
The normalization conditions read $\xi^{00}(s_0x)=0$ and $\xi^{11}(s_0x) =0$ for any $x \in \huaV_0$. Hence $\xi^{00}, \xi^{11} \in O_0$. Note that by the definition of the face maps in \eqref{eq:HomSpaceMaps}, $\widecheck{d}_0\xi=\xi^{11}$ and $\widecheck{d}_1\xi= \xi^{00}$, and the normalization conditions state precisely the fact that these are in $\huaV^{1*}_0$, as should be expected. 
For $i=0,1$, the multiplicativity equations for $r=iii$ read $\xi^{ii}(w\cdot v)=\xi^{ii}(v) + \xi^{ii}(w)$ for and $(v,w)\in\Lambda^2_1(\huaV)$. These are exactly the multiplicativity conditions at level 0 for each of the components $\xi^{ii}$, hence they are automatically satisfied by the same argument as before. 
The only remaining equations are the multiplicativity conditions for $001$ and $011$. We claim that they define the face maps $\widecheck{d}_0\xi=\xi^{11}$ and $\widecheck{d}_1\xi= \xi^{00}$. In fact, these equations read, for any $(v,w)\in \Lambda^2_1(\huaV)$,
\begin{equation*}
    \xi^{00}(w) = \xi^{01}(v) - \xi^{01}(w \cdot v), \qquad 
    \xi^{11}(v) = \xi^{01}(w) - \xi^{01}(w \cdot v).
\end{equation*}
Again, by Example \ref{ex:MooreFillers1Gpd}, $w \cdot v = d_1\mu^2_1(v,w)$ and we get 
\begin{equation*}
\begin{split}
        \xi^{00}(w) &= \xi^{01}(v) - \xi^{01}(v + w - 1d_1v) = \xi^{01}(v - 1d_1v),\\
        \xi^{11}(v) &= \xi^{01}(w) - \xi^{01}(v + w - 1d_0w) = \xi^{01}(w - 1d_0w).
\end{split}
\end{equation*}
Hence $\xi^{01}$ can be seen as the only independent component of $\xi$ and $\huaV^{1*}_1 = \huaV_1^*$.

By the discussion at the end of Section \ref{sec:ndual-overview}, the multiplication can now be equivalently defined without computing Moore fillers by the multiplicativity equation \eqref{eq:MultiplicativityMostGeneral} for $r=012 \in \Delta[2]_2$. Let $\phi \in \huaV^{1*}_2 \cong \Lambda^2_1(\huaV^{1*})$. Written as a horn this is $(\phi^{12}, \phi^{01}) = (\eta, \xi)$, and $\xi \cdot \eta = m_1(\eta, \xi) = \phi^{02}$ so the equation for $012$ is 
\begin{equation*}
    (\xi \cdot \eta) ( w \cdot v) = \eta(v) + \xi (w),
\end{equation*}
as expected. The other multiplications can be defined analogously by renaming the different components of $\phi$.
\end{proof}

\begin{remark}
    As an immediate observation, this result is consistent with the well-known dual $\VB$ groupoid construction of \cite{Pradines1988}.
    if $\huaV \to M$ is a $\VB$ groupoid over the identity groupoid of a manifold $M$, then each fiber of the dual $\VB$ groupoid of $\huaV$ concides with the 1-dual $(\huaV |_m)^{1*}$ at each point $m \in M$. 
    In particular, the restriction of any $\VB$ groupoid $\huaV \to \huaG$ to the units $\huaG_0$ is such a $\VB$ groupoid, and its fiber at any $p\in \huaG_0$ is exactly the $\VS$ 1-dual $(\huaV|_p)^{1*}$ of the fiber of $\huaV$ at $p$. 
\end{remark}

\begin{remark}\label{rem:VS1Dual-DK}
    A straightforward computation shows that the normalized complex of the 1-dual of $\huaV$ is isomorphic to 
    \begin{equation*}
        \huaV_0^* \xrightarrow{-d_1^*} (\ker d_0)^*,
    \end{equation*}
    which is precisely the 1-shifted dual $N(\huaV)^*[-1]$ of the normalized complex of $\huaV$. Therefore, in this case $\huaV^{1*} \cong DK(N(\huaV)^*[-1])$ and the 1-dual pairing is non-degenerate ``on the nose''.
\end{remark}

\subsection{The 2-dual of a \texorpdfstring{$\VS$}{VS} 2-groupoid}\label{sec:computation-2dual}

Let $\huaV$ be a $\VS$ 2-groupoid. Following Section \ref{sec:annihilators-and-dual-kernels}, the solution spaces of the normalization conditions \eqref{eq:NormalizationMostGeneral} appearing in the 2-dual $\huaV^{2*}$ are the degree 1 2-dimensional degeneracy annihilators $O_0 = \Ann(s_0\huaV_1)$ and $O_1 = \Ann(s_1\huaV_1)$, and their intersection $O_{01} = \Ann(D\huaV_2)$, the degree 2 2-dimensional degeneracy annihilator. 

By the dual sequences of \eqref{eq:VSnDualSESd_is_i} and \eqref{eq:VSnDualSESd_iplus1s_i},  we have the isomorphisms
\begin{equation} \label{eq:O0O1}
    O_0 \cong (\ker d_0)^*\cong (\ker d_1)^*, \quad O_1 \cong (\ker d_1)^*\cong (\ker d_2)^*.
\end{equation}
Analogously, the dual of the sequence \eqref{eq:VSnDualSESp_kmu_k} for $0\le k\le 2$ gives isomorphisms
\begin{equation} \label{eq:O01}
    O_{01} \cong (\ker p^2_0)^* \cong (\ker p^2_1)^* \cong (\ker p^2_2)^*,
\end{equation}
which are the dual maps of the core projections
\begin{equation}\label{eq:VS2DualCoreProjections}
    \begin{split}
        \gamma_0 X &= X - s_0d_1X - s_1d_2X + s_0d_2 X : \huaV_2 \to \ker p^2_0,\\
        \gamma_1 X &= X - s_0d_0X - s_1d_2X + 1d_1d_0 X : \huaV_2 \to \ker p^2_1,\\
        \gamma_2 X &= X - s_0d_0X - s_1d_1X + s_1d_0 X : \huaV_2 \to \ker p^2_2,
    \end{split}
\end{equation}
appearing in \eqref{eq:VSnDualSESp_kmu_k}. As in Section \ref{sec:computation-1dual}, the isomorphisms between the three different degree 2 cores are given by the restrictions of these projections. 

\begin{remark}\label{rem:VS2Dual-det-norm-comps}
    These isomorphism translate into the useful principle that an element of $\huaV_2^*$ that satisfies a single normalization condition --- i.e. an element of $O_0$ or $O_1$ --- is determined by its evaluation on elements of any degree 1 core $\ker d^2_i$ for $i = 0,1,2$. 
    In the same way, an element of $\huaV_2^*$ that satisfies both normalization conditions --- i.e. an element of $O_{01}$ --- is determined by its evaluation on elements of any degree 2 core $\ker p^2_i$ for $i = 0,1,2$.
    We will use this repeatedly in the computation of the 2-dual.
\end{remark}

\begin{remark}
    Observe that $\gamma_0^*|_{O_0}=\gamma_1^*|_{O_0}$, while $\gamma_2^*|_{O_0}$ is a different map. Similarly,  $\gamma_1^*|_{O_1}=\gamma_2^*|_{O_1}$, while $\gamma_0^*|_{O_1}$ is a different map.  
\end{remark}

\begin{theorem}[2-dual of a $\VS$ 2-groupoid] \label{thm:VS2Dual}
Let $\huaV$ be a $\VS$ 2-groupoid. Then its 2-dual $\huaV^{2*}$ is
    \begin{equation} \label{eq:2-dual}
                    \huaV_2^* \times_{O_{01}} 
                    \left(O_0 \times_{O_{01}} O_1\right) \times_{O_{01}} O_0 \rightthreearrows O_0 \times_{O_{01}} O_1 \rightrightarrows O_{01},
    \end{equation}
where the fiber products at levels 1 and 2 are 
\begin{equation*}
    \begin{split}
        \huaV^{2*}_1 &= \{(\eta^{001}, \eta^{011}) \in O_0 \oplus O_1 \mid  \eta^{001}(\gamma_1X) = \eta^{011}(\gamma_1X), \quad \forall X \in \huaV_2\},\\
        \huaV^{2*}_2 &= \{(\phi^{012}, \phi^{112}, \phi^{122}, \phi^{001}) \in \huaV_2^* \oplus 
        \left(O_0 \oplus O_1\right) \oplus O_0| \\
        &\qquad \qquad  \phi^{012}(\gamma_0 X) = \phi^{112}(\gamma_1 X) = \phi^{122}(\gamma_1 X) \text{ and } \phi^{012}(\gamma_2 X) = \phi^{001}(\gamma_1 X), \quad \forall X \in \huaV_2\},
    \end{split}
\end{equation*}
and this is equipped with the following face and degeneracy maps\footnote{Notice that since $\eta^{011}\in O_1=\Ann(s_1\huaV_1)$, $\eta^{011}(s_1d_2(X))=0$. Thus we have a simplification (rather than typo) in \eqref{eq:vs-2dual-face-degen-1}. Similarly for other simplifications in \eqref{eq:vs-2dual-face-degen-1} and \eqref{eq:vs-2dual-face-degen-2}. }:
\begin{equation}\label{eq:vs-2dual-face-degen-1}
    \begin{aligned}
        \check{d}_0^1(\eta^{001},\eta^{011})(X)
        &= \eta^{011}(\gamma_0 X) = \eta^{011} (X - s_0d_1 X + s_0d_2 X),
        \\
        \check{d}_1^1(\eta^{001},\eta^{011})(X) 
        &= \eta^{001}(\gamma_2 X)= \eta^{001}(X - s_1d_1 X + s_1d_0 X),
        \\
        \check{s}_0^0(\epsilon) &= (\epsilon, \epsilon)
    \end{aligned}
\end{equation}
for all $(\eta^{001}, \eta^{011}) \in O_0 \times_{O_{01}} O_1$, $\epsilon \in O_{01}$, $X \in \huaV_2$, and
\begin{equation}\label{eq:vs-2dual-face-degen-2}
    \begin{aligned}
        \check{d}^2_0 (\phi^{012}, \phi^{112}, \phi^{122}, \phi^{001}) 
        &= (\phi^{112}, \phi^{122}),
        \\
        \check{d}^2_1 (\phi^{012}, \phi^{112}, \phi^{122}, \phi^{001}) 
        &= (\phi^{012} - \phi^{012}s_0d_0 + \phi^{001}s_1d_2, \phi^{012} - \phi^{012}s_1d_2 + \phi^{122}s_0d_0),
        \\
        \check{d}^2_2 (\phi^{012}, \phi^{112}, \phi^{122}, \phi^{001})
        &= (\phi^{001}, \phi^{012} - \phi^{012}s_1d_1 + \phi^{112}s_1d_0),
        \\
        \check{s}^1_0(\eta^{001}, \eta^{011})
        &= (\eta^{001}, (\eta^{001}, \eta^{011}), \check{d}_1(\eta^{001},\eta^{011})),
        \\
        \check{s}^1_1(\eta^{001}, \eta^{011})
        &= (\eta^{011}, \check{s}_0\check{d}_0(\eta^{001}, \eta^{011}), \eta^{001}),
    \end{aligned}
\end{equation}
for all $(\phi^{012}, \phi^{112}, \phi^{122}, \phi^{001}) \in \huaV^{2*}_2$ and $(\eta^{001}, \eta^{011}) \in \huaV^{2*}_1$. 
The multiplication $\widecheck{m}_1$ is defined by the property that 
\begin{equation}\label{eq:mul-pairing}
    \langle \widecheck{m}_1 (\phi, \phi', \phi''), m_1(W, Y, Z) \rangle =
    \langle \phi \square \phi' \phi'', W \square Y Z \rangle
    = \langle \phi, W \rangle 
    + \langle \phi', Y \rangle 
    - \langle \phi'', Z \rangle,
\end{equation}
for any $(\phi, \phi', \phi'') \in \Lambda^3_1(\huaV^{2*})$ and any $(W, Y, Z) \in \Lambda^3_1(\huaV)$.\footnote{Note that this equation relates only the interior component $(\phi\square\phi'\phi'')^{012}$ of the product to the interior components of the factors. As previously explained in Section \ref{sec:ndual-overview}, the other components can be inferred by the simplicial identities. We write these in \eqref{eq:m_1}.}. 
\end{theorem}

\begin{proof}
Similarly to the computation of the 1-dual, we follow Section \ref{sec:ndual-overview} and compute $\huaV^{2*}$ level by level by solving the linear equations given by the multiplicativity  \eqref{eq:MultiplicativityMostGeneral} and normalization \eqref{eq:NormalizationMostGeneral} conditions. 
This is however more complicated than in the 1-dual case, as there are many instances where multiple equations determine the same variable, so we will also need to check that the equations do not over-determine the solutions, that is the solution space is not empty. 
For ease of reading, we organize the proof into subsections.

\subsubsection{Level 0}\label{sec:level0} 
For $\epsilon \in \huaV^{2*}_0$, the only 2-simplex in $\Delta[0]_2$ is $000$, so we have only one component $\epsilon = \epsilon^{000}$. Since $000 = s_0 00 = s_1 00$, we have two normalization conditions \eqref{eq:NormalizationMostGeneral}, and we have that $\epsilon \in \Ann(s_0\huaV_1 + s_1\huaV_1) = O_{01}$. There is only one multiplicativity condition, which is \eqref{eq:MultiplicativityMostGeneral} for $0000$, which reads
\begin{equation*}
    \epsilon (W \square Y Z) = \epsilon(W) + \epsilon(Y) - \epsilon(Z).
\end{equation*}
But this is already implied  by  \eqref{eq:TriMultOverPointAll} and the normalization condition of $\epsilon$. Therefore $\huaV^{2*}_0=O_{01}$.

\subsubsection{Level 1}\label{sec:level1} 
Any $\eta \in \huaV^{2*}_1$ consists of four components $\eta^{000}, \eta^{001}, \eta^{011}, \eta^{111}$, which are normalized in the following way:
    \begin{equation*}
        \eta^{000}, \eta^{111} \in \Ann(s_0\huaV_1 + s_1\huaV_1)=O_{01},
        \quad \eta^{001} \in \Ann(s_0\huaV_1)=O_0,
        \quad \eta^{011} \in \Ann(s_1\huaV_1)=O_1.
    \end{equation*}

The multiplicativity conditions are
    \begin{equation}\label{eq:VS2DualGenericMult}
        \eta^{d_1r}((W\square Y Z)) = \eta^{d_0r}(W) + \eta^{d_2r}(Y) - \eta^{d_3r}(Z)
    \end{equation}
for all $(W,Y,Z) \in \Lambda^3_1(\huaV)$, and $r = 0000, 0001, 0011, 0111, 1111$. The two equations relative to $r = 0000$ and $r = 1111$ are automatically satisfied by the same argument as in level 0.

For the others, we make use of Remark \ref{rem:VS2Dual-det-norm-comps}, by which $\eta^{000}, \eta^{111} \in O_{01}$ are determined by their evaluation on elements in $\ker p^2_j$ for $0 \le j \le 2$.
For $r=0001$, since for any $Z\in \ker p^2_2$ we get $0\square Z Z = 0$ by \eqref{eq:TriMultOverPointAll},  \eqref{eq:VS2DualGenericMult} reads
    \begin{equation*}
            0 = \eta^{001}(Z) - \eta^{000}(Z) \quad
            \iff \quad \eta^{000}(Z) = \eta^{001}(Z), \quad \forall Z \in \ker p^2_2.
    \end{equation*}
Hence by Remark \ref{rem:VS2Dual-det-norm-comps}, $\eta^{000}$ is entirely determined by $\eta^{001}$.  To see what $\eta^{000}$ is when evaluated on a generic $X\in \huaV_2$ we use the projections in \eqref{eq:VS2DualCoreProjections} and obtain that
    \begin{equation}\label{eq:VS2Dual000}
        \eta^{000}(X) = \eta^{000}(X - s_0d_0X - s_1d_1X + s_1d_0X)= \eta^{001}(\gamma_2 X),
    \end{equation}
by the normalization condition on $\eta^{000}$. It is easy to see that imposing \eqref{eq:VS2Dual000} is equivalent to the multiplicativity condition \eqref{eq:VS2DualGenericMult} for $r=0001$ because substituting \eqref{eq:VS2Dual000} in \eqref{eq:VS2DualGenericMult} does not impose extra conditions on $\eta^{001}$.

A symmetric argument applies to $r=0111$, by using the duality principle in Remark \ref{rem:front-to-back}. Thus we have that 
\begin{equation}\label{eq:VS2Dual111}
    \eta^{111}(X) = \eta^{011}(\gamma_0(X))
\end{equation}
is equivalent to the multiplicativity condition \eqref{eq:VS2DualGenericMult} for $r=0111$ and $\eta^{111}$ is completely determined by $\eta^{011}$. 

Lastly, we take care of \eqref{eq:VS2DualGenericMult} for $r=0011$. For $Y\in \ker p^2_1$, $0\square Y0 = Y$, and this equation reads $\eta^{011}(Y) = \eta^{001}(Y)$. Again, by the normalization condition,  \eqref{eq:VS2DualGenericMult} for $r=0011$ implies that
    \begin{equation}\label{eq:VS2DualLevel1FiberProductCondition}
        \eta^{011}(X - s_0d_0 X) = \eta^{001}(X - s_1d_2 X) \iff \eta^{011}(\gamma_1 X) = \eta^{001}(\gamma_1 X), 
    \end{equation} for any $X \in \huaV_2$. 
Conversely, with \eqref{eq:VS2DualLevel1FiberProductCondition}, using normalization conditions,  \eqref{eq:VS2DualGenericMult} for $r=0011$ becomes automatic:
    \begin{equation*}
        \begin{split}
            \eta^{011}(W + Y - Z - s_0d_1W + s_0d_0Z) 
            &= \eta^{011}(W) + \eta^{001}(Y) - \eta^{001}(Z)\\
\iff            \eta^{011}(Y - s_0d_0Y) + \eta^{011}(- Z + s_0d_0 Z) &= \eta^{001}(Y) - \eta^{001}(Z)\\
\iff            \eta^{001}(Y - \cancel{s_1d_2Y}) + \eta^{001}(- Z + \cancel{s_1d_2 Z}) &= \eta^{001}(Y) - \eta^{001}(Z)
        \end{split}
    \end{equation*}
Thus \eqref{eq:VS2DualLevel1FiberProductCondition} is equivalent to \eqref{eq:VS2DualGenericMult} for $r=0011$.

In summary, each $\eta$ is determined by the pair of components $(\eta^{001},\eta^{011})$ which satisfies \eqref{eq:VS2DualLevel1FiberProductCondition}. That is, $\huaV^{2*}_1$ is the fiber product 
    % https://q.uiver.app/?q=WzAsNCxbMCwwLCJDXzBeKlxcb3BsdXNfeyhcXGtlciBwXjJfMSleKn1DXzFeKiJdLFswLDEsIkNfMF4qIl0sWzEsMSwiKFxca2VyIHBeMl8xKV4qIl0sWzEsMCwiQ18xXioiXSxbMCwxXSxbMCwyLCIiLDAseyJzdHlsZSI6eyJuYW1lIjoiY29ybmVyIn19XSxbMywyLCJpZC1zXzBkXzAiXSxbMCwzXSxbMSwyLCJpZC1zXzFkXzIiLDJdXQ==
    \[\begin{tikzcd}[ampersand replacement=\&]
	{\huaV^{2*}_1 \cong O_0\times_{(\ker p^2_1)^*}O_1} \& {O_1} \\
	{O_0} \& {(\ker p^2_1)^*}
	\arrow[from=1-1, to=2-1]
	\arrow["\lrcorner"{anchor=center, pos=0.125}, draw=none, from=1-1, to=2-2]
    \arrow["{\gamma_1^*=(id-s_0d_0)^*}", from=1-2, to=2-2]
	\arrow[from=1-1, to=1-2]
	\arrow["{\gamma_1^*=(id-s_1d_2)^*}"', from=2-1, to=2-2]
    \end{tikzcd}\]
    with face maps 
    \begin{equation*}
        \begin{split}
            \check{d}_0(\eta^{001},\eta^{011}) &= \eta^{111} = \eta^{011} \circ \gamma_{0}\\
            \check{d}_1(\eta^{001},\eta^{011}) &= \eta^{000} = \eta^{001} \circ \gamma_{2},
        \end{split}
    \end{equation*}
by \eqref{eq:VS2Dual000}, \eqref{eq:VS2Dual111}.  The degeneracy map $\check{s}_0 : \huaV^{2*}_0 \to \huaV^{2*}_1$ is the obvious diagonal map.

\subsubsection{Level 2 --- Solving Equations}\label{sec:level2} An element $\phi\in \huaV^{2*}_2$ consists of 10 components, one for each 2-simplex of $\Delta[2]$. Nine of them are normalized in the sense that
    \begin{equation*}
        \phi^{iii} \in \Ann(s_0\huaV_1 + s_1\huaV_1)=O_{01},
        \quad \phi^{iij} \in \Ann(s_0\huaV_1)=O_0,
        \quad \phi^{ijj} \in \Ann(s_1\huaV_1)=O_1,
    \end{equation*}
for all $0\le i\le j \le 2$. The component $\phi^{012}$ is the only one with no normalization conditions. 
    
There are 15 multiplicativity conditions from \eqref{eq:MultiplicativityMostGeneral}, which we write for $m_1$ as 
\begin{equation}\label{eq:VS2DualGenericMult-c}
    \phi^{d_1r}((W\square Y Z)) = \phi^{d_0r}(W) + \phi^{d_2r}(Y) - \phi^{d_3r}(Z)
\end{equation}
for all $(W,Y,Z) \in \Lambda^3_1(\huaV)$ and for all $r=ijkl$ with $0\le i\le j\le k\le l \le 2$.  

Those for $r = iiii$ with $0\le i\le 2$ are automatically satisfied, by the same argument in \ref{sec:level0};    those for $r= ijjj, iijj, ijjj$ for $0\le i\le j \le 2$ can also be treated as in \ref{sec:level1}, and they impose that each pair $(\phi^{iij},\phi^{ijj})$ forms an element of $\huaV^{2*}_1$.
Notice that each pair $(0001,0002)$, $(0111,1112)$, and $(0222,1222)$ gives two possible ways to determine $000$, $111$ and $222$, respectively. We will check later that they give consistent results by using other multiplicativity conditions in \ref{sec:consist}. Now it only remains to solve the three equations for $r = 0012, 0112, 0122$. 

We begin with $r=0012$. Then \eqref{eq:VS2DualGenericMult-c} reads
    \begin{equation}\label{eq:VS2DualMult0012}
        \phi^{012}((W\square Y Z)) = \phi^{012}(W) + \phi^{002}(Y) - \phi^{001}(Z),
    \end{equation} for all $(W,Y,Z)\in \Lambda^3_1(\huaV)$.
As before, since $\phi^{002} \in O_0$, by picking $Y\in \ker d_0^2$, we can determine $\phi^{002}$ completely (see Remark \ref{rem:VS2Dual-det-norm-comps}). By plugging in $W=0$ and $Y\in \ker d_0^2$,  \eqref{eq:VS2DualMult0012} becomes
    \begin{equation}\label{eq:VS2DualMult0012with0}
        \phi^{012}(0\square YZ) = \phi^{012}(Y - Z + s_1d_1Z) = \phi^{002}(Y) - \phi^{001}(Z),
    \end{equation}
where $Z \in \ker d_0^2$ is any (2,1)-horn filler of $(0, d_2Y)$. Observe that $s_1d_2Y$ is one such horn filler, since $d_0d_2Y=d_1d_0Y = 0$. A general filler is given by $s_1d_2Y + k_1$, with $k_1\in \ker p^2_1$. In other words the space of (2,1)-fillers is an affine space modelled over $\ker p^2_1$. Since $s_1d_1(s_1d_2Y + k_1) = s_1d_2Y + s_1d_1 k_1$, by using \eqref{eq:TriMultOverPointAll}, \eqref{eq:VS2DualMult0012with0} becomes
    \begin{equation*}
        \phi^{002}(Y) = \phi^{012}(Y - k_1 + s_1d_1k_1) + \phi^{001}(s_1d_2Y + k_1).
    \end{equation*}
Thus for any $X\in \huaV_2$, since $X-s_0d_0X \in \ker d_0^2$ and $\phi^{002}\in O_0$, we get
    \begin{equation}\label{eq:VS2Dual002GenericFiller}
        \phi^{002}(X) = \phi^{002}(X-s_0d_0X) = \phi^{012}(X - s_0d_0X - k_1 + s_1d_1k_1) + \phi^{001}(s_1d_2X + k_1), 
    \end{equation}
where, for the last term, we used the fact that $s_1d_2(s_0d_0X) =s_0s_0 d_1d_0X$ and $\phi^{001}\in O_0$.    
Now we need to see what conditions are imposed on $\phi^{012}$ and $\phi^{001}$ by \eqref{eq:VS2DualMult0012} after inserting \eqref{eq:VS2Dual002GenericFiller}.
With \eqref{eq:TriMultOverPointAll}, \eqref{eq:VS2DualMult0012} reads
\begin{equation*}
    \phi^{012}(W + Y - Z - s_0d_1W + s_0d_0Z - s_1d_0Z +s_1d_1Z) 
            = \phi^{012}(W) + \phi^{002}(Y) - \phi^{001}(Z), 
\end{equation*}  for all $(W, Y, Z)\in \Lambda^3_1(\huaV)$. By inserting \eqref{eq:VS2Dual002GenericFiller}
we get
    \begin{equation*} \label{eq:c012-2}
        \begin{split}
            \phi^{012}(Y - Z - s_0d_1W + s_0d_0Z - s_1d_0Z +s_1d_1Z)
            &= \phi^{012}(Y - s_0d_0Y - k_1 + s_1d_1k_1) \\
            &\qquad + \phi^{001}(s_1d_2Y + k_1 - Z),
        \end{split}
    \end{equation*}        
with $k_1$ an arbitrary element in $\ker p^2_1$. Because $d_1W = d_0Y$, $d_2Y = d_2Z$ and $\phi^{001}(\gamma_1 X) = \phi^{001}(X - s_1d_2 X)$,this further simplifies to      
\begin{equation*}          
        \begin{split}    
            \phi^{012}(Z - s_0d_0Z + s_1d_0Z - s_1d_1Z - k_1 + s_1d_1k_1)
            &= \phi^{001}(Z - s_1d_2Z - k_1)\\
          \iff  \phi^{012}(\gamma_2(Z - k_1)) 
            &= \phi^{001}(\gamma_1(Z - k_1)).
        \end{split}
    \end{equation*}
This means that \eqref{eq:VS2DualMult0012} holds if and only if \eqref{eq:VS2Dual002GenericFiller} holds for all $X\in \huaV_2$ and $k_1 \in \ker p^2_1$, and 
    \begin{equation}\label{eq:VS2DualLevel2FiberProductCondition1}
        \phi^{012}(\gamma_2 X) = \phi^{001}(\gamma_1 X), \quad \forall X\in \huaV_2.
    \end{equation}
In fact, to use \eqref{eq:VS2Dual002GenericFiller} to determine $\phi^{002}$, we need to show the right-hand side of \eqref{eq:VS2Dual002GenericFiller} does not depend on the choice of $k_1 \in \ker p^2_1$. However this follows precisely from \eqref{eq:VS2DualLevel2FiberProductCondition1} for $X=k_1$. In summary, \eqref{eq:VS2DualMult0012} is equivalent to \eqref{eq:VS2Dual002GenericFiller} for any $Y$ and any filler. Because this must hold for any filler, \eqref{eq:VS2DualLevel2FiberProductCondition1} must also hold.

The front-to-back symmetric case of $r=0122$ can be treated analogously by Remark \ref{rem:front-to-back}. The multiplicativity condition for $r=0122$ reads, for any $(W, X, Y) \in \Lambda^3_2(\huaV)$,
\begin{equation}\label{eq:VS2DualMult0122with0}
    \phi^{022}(X) = \phi^{012}(W X\square 0) + \phi^{122}(W) = \phi^{012}(-W + X +s_0d_1W) +\phi^{122}(W).  
\end{equation}
This is equivalent to imposing
\begin{equation}\label{eq:VS2DualLevel2FiberProductCondition2}
    \phi^{012}(\gamma_0 X) = \phi^{122}(\gamma_1 X), \quad \forall X \in \huaV_2
\end{equation}
and defining $\phi^{022}$ in terms of $\phi^{012}$ and $\phi^{001}$ by
\begin{equation}\label{eq:VS2Dual022GenericFiller}
    \phi^{022}(X) = \phi^{012}(X - k_1 + s_0d_1k_1) +\phi^{122}(s_0d_0X + k_1), \quad \forall X\in \huaV_2
\end{equation}    
for any choice of $k_1\in \ker p^2_1$. As before, \eqref{eq:VS2DualLevel2FiberProductCondition2} for $X=k_1$ is equivalent to \eqref{eq:VS2Dual022GenericFiller} being independent of the choice of $k_1$. 

Moving on to $r=0112$, the multiplicativity condition reads 
    \begin{equation}\label{eq:VS2DualMult0112}
        \phi^{012}(W \square Y Z) = \phi^{112}(W) + \phi^{012}(Y) - \phi^{011}(Z), \quad \forall (W,Y,Z)\in \Lambda^3_1(\huaV).
    \end{equation} 
Once more, we want to determine $\phi^{011}$ from $\phi^{012}$ and $\phi^{112}$. 
Take $Z \in \ker d^2_2$, then \eqref{eq:VS2DualMult0112} implies that
    \begin{equation*}
        \begin{split}
            \phi^{011}(Z) &= - \phi^{012}(W \square 0 Z) + \phi^{112}(W) \\
            &= \phi^{012}( - W + Z - s_0d_0Z + s_1d_0Z - s_1d_1Z) + \phi^{112}(W),
        \end{split}
    \end{equation*}
for any (2,0)-horn filler $W\in \ker d^2_1$ of $(0, d_0Z)$. Since $s_1d_0 Z - s_0d_0 Z$ is one such horn filler, an arbitrary filler is of the form $s_1d_0 Z - s_0d_0 Z + k_0$ with $k_0 \in \ker p^2_0$. Hence we have
    \begin{equation*}
        \phi^{011}(Z) = \phi^{012}(Z - s_1d_1Z - k_0) + \phi^{112}(s_1d_0Z + k_0), \quad \forall k_0 \in \ker p^2_0,
    \end{equation*} 
since $\phi^{112}\in O_0 $. 
Again for all $X\in \huaV_2$, since  $X-s_1d_2X \in \ker d^2_2$ and $\phi^{011}\in O_1$, we have 
    \begin{equation}\label{eq:VS2Dual011GenericFiller}
        \phi^{011}(X) = \phi^{011}(X - s_1d_2X) = \phi^{012}(X - s_1d_1X - k_0) +\phi^{112}(s_1d_0X + k_0),
    \end{equation}
where, for the last term, we use again the fact that $s_1d_0(s_1d_2X) =s_0s_0 d_0 d_2X$ and $\phi^{112}\in O_0$. 
With \eqref{eq:TriMultOverPointAll}, \eqref{eq:VS2DualMult0112} reads
    \begin{equation*}
        \phi^{012}(W + Y - Z - s_0d_1W + s_0d_0Z - s_1d_0Z + s_1d_1Z) 
        = \phi^{112}(W) + \phi^{012}(Y) - \phi^{011}(Z), 
    \end{equation*}  
for all $(W, Y, Z)\in \Lambda^3_1(\huaV)$.  
By inserting \eqref{eq:VS2Dual011GenericFiller} in this, we get
    \begin{equation*}
            \begin{split}
                \phi^{012}(W - Z - s_0d_1W + s_0d_2W - s_1d_2W + s_1d_1Z) 
                &= \phi^{112}(W) - \phi^{012}(Z - s_1d_1Z - k_0) \\
                &\qquad - \phi^{112}(s_1d_0 Z + k_0),
            \end{split}
    \end{equation*}        
with $k_0$ an arbitrary element in $\ker p^2_0$. Because $d_0Z = d_2W$, and $\phi^{112}\in O_0$, this further simplifies to      
    \begin{equation*}       
        \begin{split} 
            \phi^{012}(W - s_0d_1W + s_0d_2W - s_1d_2W - k_0)
            &= \phi^{112}(W - s_1d_2W - k_0)\\
            \iff  \phi^{012}(\gamma_0(W - k_0)) 
            &= \phi^{112}(\gamma_1(W - k_0)).   
        \end{split}
    \end{equation*}
This means that \eqref{eq:VS2DualMult0112} holds if and only if \eqref{eq:VS2Dual011GenericFiller} holds for all $X\in \huaV_2$ and $k_0 \in \ker p^2_0$ and 
    \begin{equation}\label{eq:VS2DualLevel2FiberProductCondition3}
        \phi^{012}(\gamma_0 X) = \phi^{112}(\gamma_1 X), \quad \forall X\in \huaV_2.
    \end{equation}
The latter is precisely the fact that $\phi^{011}$ is determined by \eqref{eq:VS2Dual011GenericFiller} without depending on the choice of filler $k_0 \in \ker p^2_0$.
In summary, \eqref{eq:VS2DualMult0112} is equivalent to \eqref{eq:VS2DualLevel2FiberProductCondition3} and the fact that $\phi^{011}$ is determined by $\phi^{012}$ and $\phi^{122}$ via \eqref{eq:VS2Dual011GenericFiller}.

As a side note, \eqref{eq:VS2DualLevel2FiberProductCondition3} is also already implied by \eqref{eq:VS2DualLevel2FiberProductCondition2} and the multiplicativity condition for $r=1122$, which, by the discussion in \ref{sec:level1}, is equivalent to $\phi^{112}(\gamma_1 X) = \phi^{122}(\gamma_1 X)$.

\subsubsection{Level 2 --- Consistency}\label{sec:consist} Since some of the components are overdetermined, one needs to check consistency of the following conditions:
\begin{itemize}
\item     It is equivalent to determine $\phi^{000}$ by imposing \eqref{eq:VS2DualGenericMult-c} for $r = 0001$ and  for $r=0002$. This follows from
\begin{equation}\label{eq:001-002}
    \phi^{001}(\gamma_2 X) = \phi^{002}(\gamma_2 X), \quad \forall X \in \huaV_2.
\end{equation}
\item Analogously, it is equivalent to determine $\phi^{111}$ and $\phi^{222}$ in the two possible ways. This follows from 
\begin{equation*}
    \phi^{011}(\gamma_0 X)= \phi^{112}(\gamma_2 X), \quad \phi^{022}(\gamma_0 X) = \phi^{122}(\gamma_0 X), \quad \forall X \in \huaV_2.
\end{equation*}    
\item $\phi^{011}$ satisfies the multiplicativity condition for $r = 0011$ together with $\phi^{001}$,  when determined using the multiplicativity condition for $r=0122$. This follows from 
\begin{equation}\label{eq:001-011}
    \phi^{001} (\gamma_1 X) = \phi^{011} (\gamma_1 X), \quad \forall X \in \huaV_2.
\end{equation}
\item Analogously, $\phi^{002}$ and $\phi^{022}$ that we determined with \eqref{eq:VS2Dual002GenericFiller} and \eqref{eq:VS2Dual022GenericFiller} also satisfy the multiplicativity condition for $r=0022$. This follows from 
\begin{equation}\label{eq:002-022}
    \phi^{002} (\gamma_1 X) = \phi^{022} (\gamma_1 X), \quad \forall X \in \huaV_2.
\end{equation}
\end{itemize}
All the equations above follow from a straightforward calculation. Here we give the proof for the first one, \eqref{eq:001-002}. As the others can be obtained by similar arguments, we leave them to the interested reader.\footnote{
These consistency equations can be used to see more interesting facts, for example, by \eqref{eq:VS2Dual002GenericFiller} and \eqref{eq:VS2Dual022GenericFiller}, \eqref{eq:002-022} is $\phi^{012}(\gamma_1 X) = \phi^{002} (\gamma_1 X) = \phi^{022} (\gamma_1 X)$. Together with \eqref{eq:VS2DualLevel2FiberProductCondition1}, \eqref{eq:001-011}, \eqref{eq:VS2DualLevel2FiberProductCondition2}, and \eqref{eq:VS2DualLevel2FiberProductCondition3}, this can be summarized in the identity $\gamma_i^*\phi^{012} = \gamma_1^*\phi^{s_jd_i012}$ for any $i=0,1,2$, $j=0,1$. This can be used to rewrite $\huaV^{2*}_2$ as a different fiber product, by changing which components of $\phi$ are taken as independent variables, as in Remark \ref{rem:VS2Dual-other-choice-vars}.}

Using \eqref{eq:VS2DualCoreProjections} and normalization conditions, \eqref{eq:VS2DualLevel2FiberProductCondition1} is equivalent to 
\begin{equation*}
    \phi^{012}(X-s_0d_0X)+\phi^{001}(s_1d_2X)=\phi^{012}(s_1d_1 X -s_1d_0X)+\phi^{001}(X), \quad \forall X \in \huaV_2.
\end{equation*}
Then by combining this with \eqref{eq:VS2Dual002GenericFiller} for $k_1=0$, we get
\begin{equation*}
    \phi^{002}(X) = \phi^{012}(s_1d_1X-s_1d_0X) + \phi^{001}(X) , \quad \forall X \in \huaV_2.
\end{equation*}
Since the image of $\gamma_2$ is the intesection of the kernels of $d_1$ and $d_0$, by replacing $X$ with $\gamma_2X$ we obtain \eqref{eq:001-002}. 

\subsubsection{Level 2 --- Structure Maps and Multiplication}
In the above two subsections, we have shown that an element $\phi\in \huaV^{2*}_2$ depends only on the components $\phi^{012}$, $\phi^{112}$, $\phi^{122}$, and $\phi^{001}$. Thus $\huaV^{2*}_2 \subseteq \huaV^*_2 \oplus (O_0)^2 \oplus O_1$. Furthermore, by \eqref{eq:VS2DualLevel2FiberProductCondition1}, \eqref{eq:VS2DualLevel2FiberProductCondition2}, \eqref{eq:VS2DualLevel2FiberProductCondition3}, and \ref{sec:consist}, $\huaV^{2*}_2$ is exactly the fiber product
    \begin{equation*}
        \begin{split}
            \huaV^{2*}_2 &= \huaV^*_2 \times_{(\gamma_0^*, \gamma_2^*), (O_{01})^2, (\gamma_1^*\circ pr_1, \gamma_1^*)} \times \left(\left( O_0 \times_{\gamma_1^*, O_{01}, \gamma_1^*} O_1 \right) \times O_0 \right)\\
            &= \{(\phi^{012}, (\phi^{112}, \phi^{122}), \phi^{001}) \mid \phi^{012}\gamma_2 = \phi^{001}\gamma_1, \quad 
            \phi^{012}\gamma_0 =\phi^{122}\gamma_1= \phi^{112}\gamma_1
           \}.
        \end{split}
    \end{equation*}

    By the definition of the internal hom, the simplicial maps are described by \eqref{eq:HomSpaceMaps}. For example, the face maps are
    \begin{equation*}
        \begin{split}
            \check{d}_0(\phi^{012}, (\phi^{112}, \phi^{122}), \phi^{001}) 
            &= (\phi^{112}, \phi^{122})\\
            \check{d}_1(\phi^{012}, (\phi^{112}, \phi^{122}), \phi^{001}) 
            &= (\phi^{002}, \phi^{022}) \\
            &= \left(\phi^{012} \circ (id - s_0d_0) + \phi^{001}s_1d_2 , 
            \phi^{012} \circ (id - s_1d_2) + \phi^{122}s_0d_0 \right)\\
            \check{d}_2(\phi^{012}, (\phi^{112}, \phi^{122}), \phi^{001}) 
            &= (\phi^{001}, \phi^{011})
            = \left(\phi^{001}, 
            \phi^{012} \circ (id - s_1d_1) + \phi^{112}s_1d_0 \right)
        \end{split}
    \end{equation*}
The degeneracy maps can be similarly computed and we refer to the statement for their explicit description. 

Finally, a (3,1)-horn in $\huaV^{2*}$ can be written in components as $p^3_1\psi$ for a unique 3-simplex $\psi \in \huaV^{2*}_3$:
\begin{equation*}
    p^3_1\psi= ((\psi^{123},(\psi^{223},\psi^{233}),\psi^{112}),(\psi^{013},(\psi^{113},\psi^{133}),\psi^{001}),(\psi^{012},(\psi^{112},\psi^{122}),\psi^{001})).
\end{equation*}
Following the discussion at the end of Section \ref{sec:ndual-overview}, the multiplication $\widecheck{m}_1$ is given by
\begin{equation} \label{eq:m_1}
    \widecheck{m}_1(p^3_1\psi) = \check{d}_1 \psi  = (\psi^{023}, (\psi^{223},\psi^{233}), \psi^{002}),
\end{equation}
where $(\psi^{223},\psi^{233})$ and $\psi^{002}$ are determined as before from the horn data, while $\psi^{023}$ can be obtained by the multiplicativity condition for $r=0123$, that is, for any $(W,Y,Z) \in \Lambda^3_1(\huaV)$, we have
\begin{equation}
    \psi^{023}(W\square YZ) = \psi^{123}(W) + \psi^{013}(Y) - \psi^{012}(Z),
\end{equation}
which is the expected formula. 
\end{proof}

\begin{remark}\label{rem:VS2Dual-other-choice-vars}
    The choice of including $\phi^{112}, \phi^{122}$ and $\phi^{001}$ as ``given data'' of a certain element $\phi \in \huaV^{2*}_2$ is not unique. In fact, there are six possible combinations of choices coming from the fact that the components $001$ and $002$ are related by the multiplicativity condition relative to $r=0012$, $011$ and $112$ are related by the one for $r=0112$, while $122$ and $022$ by the one for $r=0122$. These all determine  isomorphic but different expressions of $\huaV^{2*}_2$ as a fiber product.
\end{remark} 

\begin{remark}\label{rem:VS2Dual-DK} 
    As in the case of the 1-dual, a slightly lengthier computation which appears in \cite[Prop. 2.3.19]{stefano-thesis}, shows that the normalized complex of the 2-dual of $\huaV$ is isomorphic to 
    \begin{equation*}
        \huaV_0^* \oplus (\ker d^1_0)^* \xrightarrow{(pr_2, pr_2 + d_1^* pr_1)} (\ker d^1_0)^* \oplus (\ker d^1_0)^* \xrightarrow{d_2^* pr_2 - d_2^* pr_1} (\ker p^2_2)^*,
    \end{equation*}
    which is not the 1-shifted dual $N(\huaV)^*[-2]$ of the normalized complex of $\huaV$, but is nevertheless homotopy equivalent to it by Theorem \ref{thm:EZ-thm-hom}. Therefore, in this case $\huaV^{2*} \simeq DK(N(\huaV)^*[-2])$ and the 2-dual pairing is only non-degenerate on the homology.
\end{remark}

\appendix

\section{The model \texorpdfstring{$DK(N(\huaV)^*[-n])$}{DK(N(V)*[-n])}}\label{sec:appendix-DK-N-nshifted-dual}

As anticipated in Remark \ref{rem:DK-n-shift-dual-appendix}, we now discuss the other possible model of $n$-dual imported from chain complexes via the Dold-Kan correspondence. For simplicity, throughout this section, whenever we talk about the $n$-dual of $\huaV$ we consider $\huaV$ to be a $\VS$ $n$-groupoid.

First of all, $DK(N(\huaV)^*[-n])$ is the $n$-dual of the $\VS$ $n$-groupoid $\huaV$ for $n=0, 1$ because in these cases the deformation retract in Theorem \ref{thm:EZ-thm-hom} is actually an isomorphism, as in Example \ref{ex:0dual} and Remark \ref{rem:VS1Dual-DK}. 
In these cases the $n$-dual pairing is also non-degenerate on the nose and it coincides with a trivial extension of the canonical evaluation pairing of the $n$-th level $\huaV_n$ with its dual vector space. 

In general cases (for $n\ge 2$), we would like to retain this property of the $n$-dual pairing being a trivial extension of the evaluation pairing of $\huaV_n^*$ to a canonical simplicial pairing $\huaV^{n*} \otimes \huaV \to B^n{\R}$, while weakening the non-degeneracy condition and requiring it to be non-degenerate only up to homotopy. 
This is precisely what we get, as we saw in Section \ref{sec:ndual-overview}. 
On the other hand, using the maps in \eqref{eq:EZ-AW-internal-hom}, one could produce a pairing $DK(AW^*\lambda)$ on $DK(N(\huaV)^*[-n])$, where $\lambda: N(\huaV)^*[-n] \otimes N(\huaV) \to \R[-n]$ is the evaluation pairing on chains. 
But there are two reasons why this is not a good definition:

Firstly, although at level $n$, $DK(N(\huaV)^*[-n])_n \cong \huaV_n^*$, the pairing $DK(AW^*\lambda)$ will not coincide with the evaluation pairing $ev: \huaV_n^* \otimes \huaV_n \to \R$ outside of the cases $n=0,1$. This is because, as we show in Proposition \ref{prop:DK2shiftedDual-counterex}, for $n\ge 2$, $ev$ does not extend trivially to a \textit{simplicial} pairing $DK(N(\huaV)^*[-n]) \otimes \huaV \to B^n(\R)$ even in simple cases. 

The second reason has to do with our main motivation for pursuing this project, which is its extension to the theory of higher vector bundles. 
As shown in \cite{HoyoTrentinaglia2021, HoyoTrentinaglia2023}, the Dold-Kan correspondence for a higher vector bundle $E$ over a simplicial manifold $M$ involves the extra data given by a representation up to homotopy (RUTH) of $M$ on the normalized complex $N$ of $E$. 
While the construction of $N$ is entirely analogous to that for simplicial vector spaces, a choice of a cleavage is needed to obtain the RUTH. 
In addition, a RUTH of a Lie $n$-groupoid $\huaG$ for $n \ge 2$ cannot be canonically dualized to a RUTH of $\huaG$, but only to a RUTH of its opposite $\huaG^{op}$ obtained by simplicial front-to-back symmetry as in Remark \ref{rem:front-to-back}.
This is due to the absence of an inversion map for arrows in Lie $n$-groupoids for $n\ge 2$. Notice that the inverse map gives the canonical isomorphism between a Lie groupoid and its opposite. Even assuming both $DK(N(\huaV)^*[-n])$ and $DK(AW^*\lambda)$ can be somehow defined over $\huaG$ and not its opposite, these must necessarily depend on a choice of cleavage. 
On the contrary, our approach in defining the $n$-dual and its pairing directly on the simplicial side can be extended to avoid both the choice of a cleavage and the problem with dualizing representations, directly resulting in a canonical description of the $n$-dual and its pairing. We will discuss these issues in more detail in \cite{CuecaRonchi-temp}.

We now show that the canonical evaluation pairing of $\huaV_2^*$ cannot be trivially extended to a simplicial pairing $DK(N(\huaV)^*[-2]) \otimes \huaV \to B^2\R$ in general.
To simplify the description of $DK(N(\huaV)^*[-2])$, we claim that given a $\VS$ 2-groupoid $\huaV$, $DK(N(\huaV)^*[-2])$ is isomorphic to the $\VS$ 2-groupoid $K^*(\huaV)$, which we define as 
% https://q.uiver.app/?q=WzAsMyxbMiwwLCJLX3syLDF9XioiXSxbMSwwLCJLX3syLDB9XioiXSxbMCwwLCJWXzJeKiJdLFsxLDAsIiIsMCx7Im9mZnNldCI6LTF9XSxbMSwwLCIiLDIseyJvZmZzZXQiOjF9XSxbMiwxXSxbMiwxLCIiLDAseyJvZmZzZXQiOi0yfV0sWzIsMSwiIiwwLHsib2Zmc2V0IjoyfV1d
\[\begin{tikzcd}
	{\huaV_2^*} & {(\ker d_0^2)^*} & {(\ker p_2^2)^*}
	\arrow[shift left=1, from=1-2, to=1-3]
	\arrow[shift right=1, from=1-2, to=1-3]
	\arrow[from=1-1, to=1-2]
	\arrow[shift left=2, from=1-1, to=1-2]
	\arrow[shift right=2, from=1-1, to=1-2]
\end{tikzcd}\]
with the following simplicial maps:
Given $\phi \in \huaV_2^*$, $\nu \in (\ker d_0^2)^*$, $\theta \in (\ker p_2^2)^*$, $v \in \huaV_2$, $h \in (\ker d_0^2)$, $h\in (\ker p_2^2)$, take
\begin{equation*}
\begin{array}{ccc}
    (\widecheck{d}_0^2\phi) (h) = \phi (h),
    &(\widecheck{d}_1^2\phi) (h) = \phi (h - s_0d_1 h),
    &(\widecheck{d}_2^2\phi) (h) = \phi (h - s_0d_1 h - s_1d_2 h + 1d_1d_1 h),\\
    (\widecheck{d}_0^1\nu) (k) = \nu (k),
    &(\widecheck{d}_1^1\nu) (k) = \nu (k - s_1d_2 k), &
\end{array}
\end{equation*}
and
\begin{equation*}
\begin{array}{cc}
    (\widecheck{s}_0^1\nu)(v) = \nu (v - s_0d_0v),
    &(\widecheck{s}_1^1\nu)(v) = \nu (v - s_0d_0v - s_1d_1 v + 1d_0d_0 v),\\
    (\widecheck{s}_0^0\theta) (h) = \theta (h - s_1d_1 h). &
\end{array}
\end{equation*}
A straightforward computation shows these obey the simplicial identities \eqref{eq:face-degen}. This 2-truncated simplicial vector space is then equipped with the canonical $\VS$ 2-groupoid structure by using the multiplications in Example \ref{ex:MooreFillers2Gpd}. 

\begin{proposition}
    For any $\VS$ 2-groupoid $\huaV$ we have $K^*(\huaV) \cong DK(N(\huaV)^*[-2])$.
\end{proposition}

\begin{proof}
We prove this by showing that the normalized complex of $K^*(\huaV)$ is isomorphic to the 2-shifted dual of the normalized complex of $\huaV$, $N(\huaV)^*[-2]$.    

First of all $N(\huaV)$ is the chain complex 
\begin{equation*}
    \ker p_2^2 \xrightarrow{d_2} \ker d_0^1 \xrightarrow{-d_1} \huaV_0
\end{equation*}
which is concentrated in degrees 0 to 2. Then, its $2$-shifted dual $N(\huaV)^*[-2]$ is the chain complex 
\begin{equation*}
    \huaV_0^* \xrightarrow{d_1^*} (\ker d_0^1)^* \xrightarrow{d_2^*} (\ker p_2^2)^*,
\end{equation*}
which is also concentrated in degrees 0 to 2.

The normalized chain complex of $K^*(\huaV)$ is
\begin{equation*}
    \ker\widecheck{p}^2_2
    \xrightarrow{\widecheck{d}_2} \ker\widecheck{d}_0^1 
    \xrightarrow{-\widecheck{d}_1}(\ker p_2^2)^*.
\end{equation*}

By definition of $\widecheck{d}_0^1$, $\ker\widecheck{d}_0^1 = \Ann(\ker p^2_2)\subset (\ker d_0^2)^*$. 
$\ker d_0^2$ splits as $\ker d_0^2 \cong \ker p^2_2 \oplus s_1(\ker d^1_0)$ via the splitting $h =  (h - s_1d_2h) + s_1d_2h$ for any $h \in \ker d_0^2$. So $\Ann(\ker p^2_2) \cong (\ker d^1_0)^*$, with the explicit isomorphism given by
\begin{equation*}
(\ker d^2_0)^* \supset \Ann(\ker p^2_2) \newrightleftarrows{s_1^*}{d_1^*} (\ker d^1_0)^*
\end{equation*}
This isomorphism intertwines the differentials $-\widecheck{d}^1_1: \Ann(\ker p^2_2) \to (\ker p^2_2)^*$ and $d_2^* : (\ker d_0^1)^* \to (\ker p^2_2)^*$, because
\begin{equation*}
    -(\widecheck{d}^1_1\nu)(k) = - \nu ( k - s_1d_2 k )
    = - \nu( s_1d_1(k - s_1d_2 k)) = \nu(s_1d_2 k), \quad \forall k \in \ker p^2_2.
\end{equation*} 

For $\ker\widecheck{p}^2_2 = \ker \widecheck{d}^2_0 \cap \ker \widecheck{d}^2_1$, we have $\ker\widecheck{d}_0^2 = \Ann(\ker d_0^2)\subseteq \huaV_2^*$. Moreover, for any $\phi \in \ker\widecheck{d}_0^2 \cap \ker \widecheck{d}_1^2$, we have
\begin{equation*}
    (\widecheck{d}_1^2\phi)(h) =  \phi(h) - \phi(s_0d_1 h) = -  \phi(s_0d_1 h) = 0, \quad \forall h \in \ker d_0^2.
\end{equation*}
Then, by surjectivity of $d_1: \ker d_0^2 \to \ker d_0^1$, $\phi (s_0 x) = 0$ for any $x \in \ker d_0^1$. So $\phi \in \Ann(\ker d_0^2) \cap \Ann(s_0(\ker d_0^1)) \subseteq \huaV_2^*$. The other inclusion is obvious, so $\ker\widecheck{p}^2_2 = \Ann(\ker d_0^2) \cap \Ann(s_0(\ker d_0^1))$.  Furthermore, by the splitting $\huaV_2 \cong \ker d_2^0 \oplus s_0 (\ker d_0^1) \oplus 1(\huaV_0)$ given for any $v\in \huaV_2$ by 
\begin{equation*}
    v = (v - s_0d_0 v)  + (s_0d_0 v - 1d_0d_0 v) + 1d_0d_0v, \quad \forall v\in \huaV_2,
\end{equation*}
$\Ann(\ker d_0^2) \cap \Ann(s_0(\ker d_0^1)) \cong \huaV_0^*$. 
The explicit isomorphism is 
\begin{equation*}
\huaV_2^* \supset \Ann(\ker d_0^2) \cap \Ann(s_0(\ker d_0^1)) 
\newrightleftarrows{1^*}{(d_0d_0)^*} \huaV_0^*,
\end{equation*}
which intertwines the differentials $\widecheck{d}_2^2: \Ann(\ker d_0^2) \cap \Ann(s_0(\ker d_0^1)) \to \Ann(\ker p^2_2)$ and $d_1^* : \huaV_0^* \to (\ker d^1_0)^*$, as 
\begin{equation*}
    (\widecheck{d}^2_2 (\theta d_0d_0))(s_1x) 
    = \theta (d_0d_0 (s_1x - s_0x - s_1x + 1d_1x))
    = \theta (d_1x), \quad \forall x \in \ker d^1_0.
\end{equation*}

By combining all of this we obtain the isomorphism of chain complexes
% https://q.uiver.app/#q=WzAsNixbMCwwLCJcXGtlciBcXGNoZWNre3B9XjJfMiJdLFsxLDAsIlxca2VyIFxcY2hlY2t7ZH1eMV8wIl0sWzIsMCwiKFxca2VyIHBeMl8yKV4qIl0sWzAsMSwiVl8wXioiXSxbMSwxLCIoXFxrZXIgZF8wXjEpXioiXSxbMiwxLCIoXFxrZXIgcF4yXzIpXioiXSxbMCwxLCJcXGNoZWNre2R9XjJfMiJdLFsxLDIsIi0gXFxjaGVja3tkfV4xXzEiXSxbMCwzLCItMV4qIiwyLHsib2Zmc2V0IjoxfV0sWzEsNCwic18xXioiLDIseyJvZmZzZXQiOjF9XSxbMywwLCItZF8wXipkXzBeKiIsMix7Im9mZnNldCI6MX1dLFsyLDUsIiIsMCx7ImxldmVsIjoyLCJzdHlsZSI6eyJoZWFkIjp7Im5hbWUiOiJub25lIn19fV0sWzQsNSwiZF8yXioiLDJdLFszLDQsIi1kXzFeKiIsMl0sWzQsMSwiZF8xXioiLDIseyJvZmZzZXQiOjF9XV0=
\[\begin{tikzcd}[ampersand replacement=\&,sep=scriptsize]
	{\ker \widecheck{p}^2_2} \& {\ker \widecheck{d}^1_0} \& {(\ker p^2_2)^*} \\
	{\huaV_0^*} \& {(\ker d_0^1)^*} \& {(\ker p^2_2)^*.}
	\arrow["{\widecheck{d}^2_2}", from=1-1, to=1-2]
	\arrow["{- \widecheck{d}^1_1}", from=1-2, to=1-3]
	\arrow["{1^*}"', shift right, from=1-1, to=2-1]
	\arrow["{s_1^*}"', shift right, from=1-2, to=2-2]
	\arrow["{d_0^*d_0^*}"', shift right, from=2-1, to=1-1]
	\arrow[Rightarrow, no head, from=1-3, to=2-3]
	\arrow["{d_2^*}"', from=2-2, to=2-3]
	\arrow["{d_1^*}"', from=2-1, to=2-2]
	\arrow["{d_1^*}"', shift right, from=2-2, to=1-2]
\end{tikzcd}\]
\end{proof}

\begin{proposition}\label{prop:DK2shiftedDual-counterex}
    Let $\huaV$ be the pair groupoid of $\R$, $\R^2 \rightrightarrows \R$. Then the 2-shifted pairing $\langle \cdot,\cdot \rangle: K^*(\huaV)_2 \otimes \huaV_2 \to \R$ on 2-simplices is not multiplicative.
\end{proposition}

\begin{proof}
Take arbitrary $(\phi,\chi,\psi)\in \Lambda^3_1(K^*(\huaV))$ (each of these is an element in $\huaV_2^*$), and $(W,Y,Z)\in \Lambda^3_1(\huaV)$. The pairing is multiplicative if and only if 
\begin{equation*}
    \langle \phi\square\chi\psi, W\square YZ\rangle
    = \langle \phi, W \rangle + \langle \chi, Y \rangle - \langle \psi, Z \rangle.
\end{equation*}
But now consider the elements $\widecheck{s}_0\phi = (\phi,\widecheck{s}_0\widecheck{d}_1\phi, \widecheck{s}_0\widecheck{d}_2\phi) \in \Lambda^3_1(K^*(\huaV))$, with $\widecheck{d}_1\widecheck{s}_0\phi = \phi$, and $s_1 W = (s_0d_0W,W,s_1d_2W)$ with $d_1s_1W = W$. Then if the pairing is multiplicative, 
\begin{equation*}
    \begin{aligned}
    \langle \phi, W \rangle &= \langle \phi\square(\widecheck{s}_0\widecheck{d}_1\phi)(\widecheck{s}_0\widecheck{d}_2\phi), (s_0d_0W)\square W(s_1d_2W) \rangle\\
    &= \langle \phi, s_0d_0W \rangle + \langle \phi - \phi s_0d_1, W\rangle - \langle \phi -\phi s_0d_1 - \phi s_1d_2 + 1d_1d_1, s_1d_2 W\rangle\\
    \end{aligned}
\end{equation*}
which is true if and only if
\begin{equation*}
    \begin{aligned}
    0  &= \langle \phi, s_0d_0W - s_0d_1 W - s_1d_2 W + s_0\stkout{d_1s_1}d_2 W + s_1\stkout{d_2s_1}d_2 W - 1d_1\stkout{d_1s_1}d_2 W\rangle\\
    &= \langle \phi, s_0\partial W - 1 d_1d_2 W\rangle,
    \end{aligned}
\end{equation*}
for $\partial= d_0 -d_1 +d_2$, the boundary map. In particular, for any $\theta \in \huaV_0^*$, we have $\theta d_0d_0 \in \huaV_2^*$, so that the multiplicativity of the pairing implies
\begin{equation*}
    \theta(d_0\partial W) = \theta(d_0d_2W) = \theta(d_1d_2W). 
\end{equation*}
Since $d_2^2$ is surjective (for any $\huaV$), this implies that $\partial^*\theta = 0$ for any $\theta \in \huaV_0^*$. 
For $\huaV = \R^2 \rightrightarrows \R$, when evaluating this on any $(x,0) \in \R^2$, we get that
\begin{equation*}
    \theta (x) = 0,
\end{equation*}
for any $\theta \in \R^*$ and $x \in \R$, yielding a contradiction.
\end{proof}

\begingroup
\sloppy
\printbibliography
\endgroup

\end{document}